\magnification=\magstep1
\input amstex
\documentstyle{amsppt}
\catcode`\@=11 \loadmathfont{rsfs}
\def\mycal{\mathfont@\rsfs}
\csname rsfs \endcsname \catcode`\@=\active

\vsize=6.5in

\topmatter 
\title A II$_1$ factor approach \\ to the Kadison-Singer problem
\endtitle
\author SORIN POPA \endauthor

\rightheadtext{Kadison-Singer Problem}

\affil     {\it Dedicated to R.V. Kadison and I.M. Singer} \endaffil

\address Math.Dept., UCLA, Los Angeles, CA 90095-1555\endaddress
\email popa\@math.ucla.edu\endemail

\thanks Supported in part by NSF Grant DMS-1101718 and a Simons Fellowship \endthanks

\abstract  We show that the Kadison-Singer problem, asking whether the pure states of the diagonal subalgebra
$\ell^\infty\Bbb N\subset \Cal B(\ell^2\Bbb N)$ have unique state extensions to $\Cal
B(\ell^2\Bbb N)$, is equivalent to a similar statement in II$_1$ factor framework, concerning the ultrapower inclusion 
$D^\omega \subset R^\omega$, where $D$ is the {\it Cartan subalgebra} of the hyperfinite II$_1$ factor $R$ (i.e. a maximal abelian $^*$-subalgebra of $R$ 
whose normalizer generates $R$, e.g. $D=L^\infty([0, 1]^{\Bbb Z}) \subset L^\infty([0,1]^{\Bbb Z}\rtimes \Bbb Z = R$), and $\omega$ 
is a free ultraflter.  
Instead, we prove here that 
if $A$ is any {\it singular} maximal abelian $^*$-subalgebra of $R$ (i.e., whose normalizer consists of the unitary group of $A$, e.g. $A=L(\Bbb Z)\subset 
L^\infty([0,1]^\Bbb Z)\rtimes \Bbb Z=R$), then the inclusion $A^\omega \subset R^\omega$  does satisfy the Kadison-Singer property. 
\endabstract

\endtopmatter

\document

\heading 0. Introduction \endheading

A famous problem posed by Kadison and Singer  in the late 1950s ([KS]) asks whether any pure state on the diagonal 
$\ell^\infty \Bbb N$ of the algebra $\Cal B(\ell^2 \Bbb N)$, of all linear bounded operators on the Hilbert space 
$\ell^2 \Bbb N$, has unique state extension to $\Cal B(\ell^2 \Bbb N)$. We will refer to this 
property of the inclusion of algebras $\ell^\infty \Bbb N \subset \Cal B(\ell^\infty \Bbb N)$ as the 
{\it Kadison-Singer property}. As already pointed out in [KS], 
it is equivalent to the following property for operators on the Hilbert space,  
known as the {\it paving property}: if $x\in \Cal B(\ell^2 \Bbb N)$ has only $0$ on the diagonal,  then for any $\varepsilon > 0$, 
there exists a finite partition of $\Bbb N$ into subsets $Y_1, ..., Y_n$, such that if $p_i\in \ell^\infty \Bbb N$ denotes the characteristic function 
of $Y_i$, viewed as a diagonal operator operator on $\ell^2\Bbb N$, then $\|\Sigma_{i=1}^n p_i x p_i \|\leq \varepsilon \|x\|$. 
It was later shown in [An1, An2] that this is in fact equivalent to the following finite dimensional version 
of the property, known as the {\it uniform paving property}: for any $\varepsilon > 0$, 
there exists $n=n(\varepsilon)$ such that for any $m$ and any $x \in \Cal B (\ell^2_m)$ with $0$ on the diagonal, 
there exists a partition of $\{1, 2, ..., m\}$ into $n$ sets $Y_i$, such that the corresponding diagonal operators $p_i$ 
satisfy $\|\Sigma_i p_i x p_i\|\leq \varepsilon \|x\|$. 

The Kadison-Singer problem has attracted much interest over the years, 
proving to have deep connections to a large number of fields of mathematics, with interesting 
equivalent re-formulations in harmonic analysis, frame theory, discrepancy theory, etc. Several partial results have been obtained so far 
(see e.g. [A1], [A2], [AkA], [BT], [BeHKW], [We], etc), showing for instance that certain classes of operators in $\Cal B(\ell^2 \Bbb N)$ can indeed be paved. 
We refer the reader to [CaFTW] for a beautiful, comprehensive account on this problem, and on its interdisciplinary aspects.

In this paper, we attempt a new approach to the problem, based on a reformulation in II$_1$ factor framework. Recall that a 
II$_1$ {\it factor} is a von Neumann algebra $M$ that is infinite dimensional, has trivial center and a completely additive trace state $\tau$. 
Any maximal  abelian $^*$-subalgebra (MASA) $A$ of a II$_1$ factor $M$ is diffuse (i.e. has no atoms) and if $A$ is also countably generated, 
then $A\simeq L^\infty([0,1])$ ([vN]).

II$_1$ factors can in fact be viewed as non-commutative versions of the function algebra $L^\infty([0,1])$, with the role of the Lebesgue integral 
$\int \cdot \ \text{\rm d} \mu$ played by the positive functional $\tau: M \rightarrow \Bbb C$ which satisfies $\tau(1)=1$ (it is a {\it state}) 
and $\tau (xy)=\tau(yx)$, $\forall x, y \in M$ (it is a {\it trace}). A specific type of (non-commutative) analysis has been developed in this framework,  
often exploiting the interplay between the operator norm and the Hilbert-norm implemented by the trace, as well as ergodicity properties of the $Ad$-action of  
the unitary group of $M$.  
One should note that the algebra $M_{m \times m}(\Bbb C)$, of $m$ by $m$ matrices with complex entries ($\simeq \Cal B(\ell^2_m)$),  
has both a trace state (given by the normalized trace $tr$) and is a factor, but it is finite dimensional. However, inductive limits and 
ultraproducts of these algebras give rise to II$_1$ factors. 

Thus, the most ``basic'' example of  a II$_1$ factor is the {\it hyperfinite} II$_1$ {\it factor} $R$ of Murray and von Neumann,  
defined as the infinite tensor product $(R, \tau)=\overline{\otimes}_k (M_{2 \times 2}(\Bbb C), tr)_k$. By [MvN2],  
$R$ is in fact the unique {\it approximately finite dimensional} II$_1$ factor, and by [C1] it is even the unique {\it amenable} II$_1$ factor.  
So $R$ can be represented in many different ways, for instance as the {\it group measure space} II$_1$ factor $L^\infty(X)\rtimes \Gamma$, 
associated with a free  ergodic measure preserving action of a countable amenable group $\Gamma$ on a probability space $(X, \mu)$. 
In particular, $R=L^\infty([0, 1]^\Bbb Z) \rtimes \Bbb Z$, where $\Bbb Z \curvearrowright X=[0, 1]^\Bbb Z$ is the Bernoulli action. 
When viewed this way, $R$ has $D=L^\infty(X)$ as a natural {\it Cartan subalgebra}, i.e. a MASA $D\subset R$ 
whose normalizer generates $R$. By [CFW], [OW] the Cartan subalgebra of $R$ is in fact unique, up 
to conjugacy by an automorphism of $R$. We may thus represent $D\subset R$ as the infinite tensor product  $\overline{\otimes}_k (D_2)_k\subset 
\overline{\otimes}_k (M_{2 \times 2}(\Bbb C))_k$, where 
$D_2$ is the diagonal subalgebra in $M_{2 \times 2}(\Bbb C)$. 

But  the hyperfinite II$_1$ factor $R$ also has MASAs $A\subset R$ whose normalizer 
is trivial, i.e. the only unitary elements $u\in R$ normalizing $A$,  $uAu^*=A$, are the unitaries in $A$. 
Such MASAs  are called {\it singular} and their existence  was discovered in [D1].  
A typical example of singular MASA is given by the subalgebra $L(\Bbb Z)\subset R$, generated by the canonical unitary implementing the Bernoulli action 
$\Bbb Z \curvearrowright [0, 1]^\Bbb Z$, in the above representation of the hyperfinite factor $R=L^\infty([0, 1]^\Bbb Z)\rtimes \Bbb Z$.

There is an {\it ultraproduct} procedure of constructing II$_1$ factors from 
a free ultrafilter $\omega$ on $\Bbb N$ and a sequence of factors $(M_m, \tau)$, 
with $M_m$ either II$_1$, or finite dimensional with $\text{\rm dim} M_m \nearrow \infty$ 
(see [W], [F]). The initial motivation for our work has been the observation that the Kadison-Singer 
property for $\ell^\infty \Bbb N\subset \Cal B(\ell^2\Bbb N)$, as well as its paving version, are equivalent to the analogue statements for 
the ultrapower inclusions $D^\omega \subset R^\omega$, respectively $\Pi_\omega D_m \subset \Pi_\omega M_{m \times m}(\Bbb C)$.  
Paving here means that if $x \in R^\omega$ (resp. $x\in \Pi_\omega M_{m \times m}(\Bbb C)$) has trace preserving expectation onto $D^\omega$  
(resp. $\Pi_\omega D_m$)  equal to $0$, then for any $\varepsilon > 0$, there exists a partition of $1$ with finitely many projections $p_1, ..., p_n$ in  
$D^\omega$ (resp. $\Pi_\omega D_m$), such that $\|\Sigma_{i=1}^n p_i x p_i \|\leq \varepsilon \|x\|$.  

The operator norm 
of an element $y$ in a II$_1$ factor can be calculated by the formula $\|y\|=\lim_n \tau((y^*y)^{2n})^{1/2n}$. So in order to pave $x$ one needs to 
control the ``higher moments'' $\tau((y^*y)^n)$ for $y=\Sigma_i p_i x p_i$. Our idea here is to approach such calculations  
by using a technique developed in [P6], which consists of building the paving $p_i$ by patching together small, ``infinitesimal''  pieces of projections, 
with ``incremental'' control of the moments.  Ideally, one wants to build the partition $p_i$ so that to be 
``free independent'' with respect to the given $x$, because then  the paving diminishes the operator norm  by $\sqrt{\varepsilon}$ if the mesh of the 
partition is less than $\varepsilon$, due to norm calculations in [V2]. 

In the case of  the Cartan subalgebra $D\subset R$, the independence ``breaks'' after the 3rd moment, 
more precisely we show that given any $x\in R^\omega \ominus D^\omega$, $D^\omega$ contains finite partitions with projections $p_i$ that are $3$-independent to $x$, but if $x$ normalizes 
$D$ then $xux^*u^*=u^*xux^*$, for any $u\in D^\omega$, so $4$-independence may fail in general. 

Nevertheless, our approach does provide ``free paving'' for any ultrapower $A^\omega \subset R^\omega$ of  a singular MASA   
$A \subset R$, in fact for any ultraproduct of singular MASAs in II$_1$ factors (N.B.: by [P3] any II$_1$ factor contains singular MASAs). Note that this result 
provides the first case when the Kadison-Singer property is established for a MASA in an infinite dimensional von Neumann factor.

\proclaim{0.1. Theorem (Kadison-Singer for ultrapowers of singular MASAs)} Let  $A_m \subset M_m$, $m\geq 1$, be a sequence of singular MASAs   
in $\text{\rm II}_1$ factors and denote $\text{\bf A} = \Pi_\omega A_m \subset \Pi_\omega M_m = \text{\bf M}$, their 
ultraproduct, over a free ultrafilter $\omega$ on $\Bbb N$. 
Then $\text{\bf A} \subset \text{\bf M}$ satisfies the Kadison-Singer property, i.e. any pure state on $\text{\bf A}$ 
has a unique state extension to $\text{\bf M}$. Moreover, $\text{\bf A}\subset \text{\bf M}$ has the uniform paving property:  
if $x\in \text{\bf M}$ has $0$-expectation on $\text{\bf A}$, then $\forall \varepsilon > 0$, $\exists$ $p_1, ..., p_n$ partition of $1$ with projections 
in $\text{\bf A}$, with $n \leq C \varepsilon^{-6}$ for some universal constant $C$,  
such that $\|\Sigma_{i=1}^n p_i x p_i \|\leq \varepsilon \|x\|$. 
\endproclaim

As we mentioned before, the way we prove the above result is by showing that given any $x \perp \text{\bf A}$, there exists a 
diffuse abelian subalgebra $B_0\subset \text{\bf A}$ which is free independent to $\{x, x^*\}$, i.e. any alternating word $\Pi_{i=1}^k u_i x_i$, 
with letters $x_i \in \{x, x^*\}$, $u_i \in B_0 \ominus \Bbb C$, has  $0$-trace. The presence of ``asymptotic freeness'' 
in a MASA $A\subset M$ characterizes in fact singularity, and 
for it to be satisfied, asymptotic $4$-independence is actually sufficient ($B_0\subset A^\omega$ is $n$-independent to $X \subset M^\omega \ominus A^\omega$ 
if $\tau(\Pi_{i=1}^k u_i x_i)=0$, $\forall k\leq n$, $u_i \in B_0 \ominus \Bbb C$, $x_i \in X$). In turn, we will show in Section 3 and 5.3.1 that existence 
of asymptotic 3-independence holds in any MASA.

\proclaim{0.2. Theorem (Characterizations of singularity for MASAs)} Let $A$ be a MASA in a $\text{\rm II}_1$ factor $M$. The following are  equivalent: 
\vskip .05in
$1^\circ$  $A$ is singular in $M$; 
\vskip .05in 
$2^\circ$ $A^\omega$ is maximal amenable in  $ M^\omega$; 
\vskip .05in 
$3^\circ$ $A^\omega$ is maximal among the $^*$-subalgebras $P\subset M^\omega$ that contain $A^\omega$ 
and are countably generated both as a left and right $A^\omega$-modules $($i.e., $\exists X\subset P$ countable such that ${\text{\rm sp}}A^\omega X$ 
and $\text{\rm sp} X A^\omega$ are dense in $P)$.  
\vskip .05in
$4^\circ$ Given any countable set $X\subset M^\omega \ominus A^\omega$, there exists $B_0\subset A^\omega$ 
diffuse such that $B_0$, $X$  are free independent relative to $A^\omega$.
\vskip .05in
$5^\circ$ Given any self-adjoint element $x\subset M\ominus A$, there exists $B_0\subset A^\omega$ 
diffuse such that $B_0$, $\{x\}$  are $4$-independent.
\endproclaim

The paper is organized as follows: In Section 1 we recall a classic result from [KS], 
on the equivalence between the unique state extension of pure states from $\ell^\infty \Bbb N$ to $\Cal B(\ell^2 \Bbb N)$ 
and the paving property, as well as other basic facts. In Section 2 we prove the equivalence between the  Kadison-Singer problem 
and several similar statements  in II$_1$ factor framework. In Section 3 we show that given any MASA $A$ in a II$_1$ factor $M$, 
one can pave any finite set $X \subset M\ominus A$ with respect to the $L^2$-norm given by the trace. This result has already been shown in ([P1]; see also A.1 in [P5]), 
but we give it here a different proof (which gives better estimates of the paving size), by showing that $A$ contains  
finite partitions of arbitrary small mesh that are approximately $2$-independent to $X$. In Section 4 we prove Theorem 0.1 (as Corollary 4.3). 
We do this by utilizing the $L^2$-paving from Section 3 and the ``incremental patching method'' from [P6]. In Section 5 
we derive Theorem 0.2 (as Theorem 5.2.1) and obtain several related results, including existence of approximate $3$-independence in arbitrary MASAs 
(see Theorem 5.3.1). We also formulate a conjecture generalizing Kadison-Singer (see 5.5.1). 

While we made an effort to make this paper as self-contained as possible, 
for the most basic facts on von Neumann algebras and II$_1$ factors, we refer the reader to  
the classic books [D2], [KR]. 

This work was completed during my stay at the Jussieu Math Institute in Paris, during the year 2012-2013. 
I want to gratefully acknowledge A. Connes, G. Pisier, G. Skandalis and S. Vassout, for their kind hospitality and support.

\vskip .05in

{\it Note added in the proof}.  The present paper has been posted on the arXiv on March 22nd, as arXiv:13031424. Since then, 
A. Marcus, D. Spielman and N. Srivastava 
have posted the paper  
{\it Interlacing Families II: Mixed Characteristic Polynomials and the Kadison-Singer Problem} (arXiv:1306.3969), 
where they solve the classic Kadison-Singer problem in the affirmative.  
They do this by settling a finite dimensional version of the paving property 
emphasized in [AkA] and [We]. Note that,  due to the 
equivalent re-formulation of the problem established in Theorem 2.2 below, their result also implies that the inclusions $D^\omega \subset R^\omega$ and 
$\Pi_\omega D_m \subset \Pi_\omega M_{m \times m}(\Bbb C)$ have the Kadison-Singer property.

\heading 1. Preliminaries \endheading

We recall in this section a  result   
from [KS], showing that pure states on a maximal abelian von Neumann 
subalgebra $\Cal A$ of a von Neumann algebra 
$\Cal M$ have unique state extensions to $\Cal M$ if and only if  all elements in $\Cal M$ have a certain ``paving property'' relative to $\Cal A$. 
For the reader's convenience, we have included a proof. It is essentially the original one from [KS], 
but explained in more modern terms, and adding the reformulation of paving in terms of ``relative Dixmier property''.  We   
also introduce some necessary terminology and prove some basic related results.  

\vskip .05in \noindent {\bf 1.1. Notation}. Let $\Cal M$ be a von Neumann algebra and $\Cal A \subset \Cal M$ 
a maximal abelian $^*$-subalgebra (hereafter abbreviated MASA) in $\Cal M$. If $x\in \Cal M$ then we denote 
by $C_{\Cal A}(x)$ the norm closure of the convex hull of the set $\{uxu^* \mid u\in \Cal U(\Cal A)\}$. Also, 
given a finite $n$-tuple of unitaries $V=(v_1, ..., v_n)$ in $\Cal A$ and $y\in \Cal M$, we denote $T_V(y)=n^{-1}\Sigma_{i=1}^n v_i y v_i^* \in C_{\Cal A}(y)$. 
Note right away that the commutativity of $\Cal A$ implies  $T_U(T_V(y))=T_V(T_U(y))$ for any two such 
tuples $U, V$. Also, $\|T_U(y)\|\leq \|y\|$ and $T_U(a_1ya_2)=a_1T_U(y)a_2$, 
$\forall a_1, a_2\in \Cal A$ (i.e., the maps $T_U$ are $\Cal A$-bimodular).

\proclaim{1.2. Theorem (Kadison-Singer [KS])} If $\Cal A \subset \Cal M$ is a MASA in a von Neumann algebra $\Cal M$, then the following 
conditions are equivalent: 
\vskip .05in 
\noindent 
$(1.2.1)$ Any pure state on $\Cal A$ has a unique pure state extension to $\Cal M$. 
\vskip .05in 
\noindent 
$(1.2.2)$ $C_{\Cal A}(x) \cap \Cal A \neq \emptyset$, $\forall x\in \Cal M$. 
\vskip .05in 
\noindent
$(1.2.3)$ $C_{\Cal A}(x) \cap \Cal A$ is a single point set $\{E_{\Cal A}(x)\}$, $\forall x\in \Cal M$.  
\vskip .05in 
\noindent
$(1.2.4)$ For all $x\in \Cal M$ and all $\varepsilon > 0$ there exists  
a finite partition of $1$ with projections $q_k \in \Cal A$ such that 
$\text{\rm d}(\Sigma_k q_k x q_k, \Cal A) \leq \varepsilon$. 
\vskip .05in 
\noindent
$(1.2.5)$ For all $x\in \Cal M$, there exists a unique element $E(x) \in \Cal A$ with the property that $\forall \varepsilon > 0$,  
$\exists q_k \in \Cal P(\Cal A)$ a finite partition of $1$ such that 
$\|\Sigma_k q_k x q_k - E(x)\|\leq \varepsilon$.

\vskip .05in 

Moreover, if these conditions are satisfied then $E(x)=E_{\Cal A}(x)$ and the map $E_{\Cal A}$ satisfies the following additional properties: 

\vskip .05in 
\noindent
$(i)$ $\overline{C_{\Cal A}(x)}^w \cap \Cal A =\{E_{\Cal A}(x)\}$, $\forall x\in \Cal M$. 
\vskip .05in 
\noindent
$(ii)$ $E_{\Cal A}$ is the unique conditional expectation of $\Cal M$ onto $\Cal A$. 
\vskip .05in 
\noindent
$(iii)$ Given any pure state $\psi$ on $\Cal A$, $\psi \circ E_{\Cal A}$ is the unique state extension of $\psi$ to $\Cal M$,  
and it is a pure state. 
\endproclaim 
\noindent
{\it Proof}. The implication $(1.2.3) \implies (1.2.2)$ is trivial. If $(1.2.2)$ is satisfied and  $a, b \in C_{\Cal A}(x)$ are distinct, then 
there exist tuples $U$, $V$ such that $\|T_U(x)-a\|\leq \|a-b\|/4$ and $\|T_V(y)-b\|\leq \|a-b\|/4$. But  
$\|T_V(T_U(x))-a\|= \|T_V(T_U(x)-a)\| \leq \|T_U(x)-a\|$ and $\|T_U(T_V(x))-b\|= \|T_U(T_V(x)-b)\| \leq \|T_V(x)-b\|$. Thus we have 
$$
\|a-b\| \leq  \|T_V(T_U(x))-a\| + \|T_V(T_U(x))-T_U(T_V(x))\|+ \|T_U(T_V(x))-b\|
$$
$$
= \|T_V(T_U(x))-a\| + \|T_U(T_V(x))-b\| 
 \leq \|a-b\|/4 + \|a-b\|/4=\|a-b\|/2, 
$$
a contradiction. This proves $(1.2.2) \implies (1.2.3)$. 

A similar argument shows that $(1.2.4)$ and $(1.2.5)$  are equivalent. 

Assuming now $(1.2.3)$, we prove $(1.2.5)$ as well as the properties  $(i)-(iii)$.  
Let $x\in \Cal M$, $\|x\|\leq 1$, and $\varepsilon > 0$. Let  
$u_1, ..., u_n \in \Cal A$ be so that $\|T_U(x)-E_{\Cal A}(x)\|\leq \varepsilon/2$, where $U=(u_1, ..., u_n)$. For each $i=1, ..., n$ let 
$\{e_{ij}\}_j\in \Cal A$ be spectral projections of $u_i$ such that if we denote $v_i=\Sigma_j \lambda_{ij}e_{ij}$, then 
$\|u_i -v_i\| \leq \varepsilon/4.$ Thus, if we denote  $V=(v_1, ..., v_n)$, then $\|T_U(x) - T_V(x)\|\leq \varepsilon/2$ and hence $\|T_V(x)-E_{\Cal A}(x)\|\leq \varepsilon$. 
By taking into account  that if $\{q_k\}_k$ denotes a relabeling of the set of projections $\{e_{ij}\}_{i,j}$, then  $\Sigma_k q_k T_v(x) q_k = \Sigma_k q_k x q_k$,   
we thus get 
$$
\|\Sigma_k q_k x q_k - E_{\Cal A}(x)\| = \|\Sigma_k q_k (T_v(x)-E_{\Cal A}(x)) q_k \| 
\leq \|T_v(x)-E_{\Cal A}(x) \| \leq \varepsilon,   
$$
proving the existence part of $(1.2.5)$. Since any element of the form $\Sigma_j p_j x p_j $, with $p_1, ..., p_n$ a partition with projections in $\Cal A$, 
is of the form $T_{W}(x)\in C_{\Cal A}(x)$, where $W=(w^{j-1})_{j=1}^n$, $w=\Sigma_{i=1}^n \lambda^{(i-1)}p_i$, where $\lambda=\exp(2\pi i/n)$, 
by the uniqueness in $(1.2.3)$ we get the uniqueness 
part in $(1.2.5)$ and that $E(x)=E_{\Cal A}(x)$. We have also implicitly  shown that $(1.2.5)\implies (1.2.2)$. 

Since $\Cal A$ is abelian (thus amenable), there exists a conditional expectation $\Cal E:\Cal M \rightarrow \Cal A$ (obtained by taking 
a Banach limit of appropriate averages $T_U$). By $(1.2.3)$, for any fixed $x\in \Cal M$ and any $\varepsilon > 0$, 
there exists a tuple $V=(v_1, ..., v_n)$ in $\Cal A$ 
such that fixed $\|T_V(x)-E_{\Cal A}(x)\|\leq \varepsilon$. Thus 
$$
\|\Cal E(x)-E_{\Cal A}(x)\| = \|T_V(\Cal E(x))-E_{\Cal A}(x)\| 
$$
$$
=\|\Cal E(T_V(x)-E_{\Cal A}(x))\leq \|T_V(x)-E_{\Cal A}(x)\| \leq \varepsilon.
$$
Since $\varepsilon > 0$ was arbitrary, this shows that $\Cal E(x)=E_{\Cal A}(x)$, $\forall x\in \Cal M$, proving $(ii)$.

Fix now $x_0\in \Cal M$ and let $y_0 \in \overline{C_{\Cal A}(x)}^w \cap \Cal A$ and $\{U_\iota\}_{\iota \in I}$ be a net of tuples of unitaries in $\Cal A$ such that 
the weak limit of $\{T_{U_\iota}(x_0)\}_\iota$ is equal to $y_0$. Let $\text{\rm Lim}_\iota$ be a Banach limit over $\iota$ and for each $x\in \Cal M$ denote  
$\Phi (x) =\text{\rm Lim}_\iota T_{U_\iota}(x)$. Then $\Phi: \Cal M\rightarrow \Cal M$ is linear, positive, $\Cal A$-bimodular, $\Phi(a)=a$, $\forall a\in \Cal A$,  
and  $\Phi(x_0)=y_0$. But then $\Cal E(x)=E_{\Cal A}(\Phi(x))$ is a conditional expectation of $\Cal M$ onto $\Cal A$ satisfying $\Cal E(x_0)=y_0$. By $(ii)$, 
this forces $y_0=E_{\Cal A}(x_0)$, proving $(i)$. 

Let now $\psi$ be a pure state on $\Cal A$. By Gelfand-Naimark, $\psi$ is given by the evaluation at some point in the spectrum $\Omega$ of $\Cal A$  
(thus $\Omega$ is a hyperstonian compact space and $\Cal A=C(\Omega)$). 
In particular,  $\psi$ is multiplicative and takes only the values $0, 1$ 
on the set of projections $\Cal P(\Cal A)$, with $\psi(1)=1$. This implies that any state extension  $\varphi$ of $\psi$ to $\Cal M$ has 
$\Cal A$ in its centralizer $\Cal M_\varphi$. Indeed, because if $\psi(p)=0$ for some 
$p\in \Cal A$, then by the Cauchy-Schwartz inequality for $\varphi$ we have 
$|\varphi (px)| \leq \varphi(p)^{1/2} \varphi(x^*x)^{1/2}=0$,  $|\varphi (xp)| \leq \varphi(xx^*)^{1/2}\varphi(p)^{1/2} =0$, $\forall x\in \Cal M$. 
Since for any projection $p\in \Cal A$ we either have $\psi(p)=0$ or $\psi(p)=1$ and $1\in \Cal M_\varphi$, this shows that $\Cal P(\Cal A)\subset \Cal M_\varphi$, thus 
all $\Cal A$ is contained in $\Cal M_\varphi$. Hence, $\varphi$ is constant on  $C_{\Cal A}(x)$, which contains $E_{\Cal A}(x)$ by $(1.2.3)$, 
implying that $\psi(x)=\psi(E_{\Cal A}(x))$.  This proves $(1.2.3) \implies (1.2.1)$ and $(1.2.3) \implies (iii)$. 

We have shown so far that $(1.2.2)-(1.2.5)$ are equivalent and that they imply $(1.2.1)$ and $(i)-(iii)$. 
To prove the remaining implication $(1.2.1) \implies (1.2.2)$,  let $b\in \Cal M_h$ and fix a point $t\in \Omega$ in the spectrum of $\Cal A$. 
Letting $c_0=\inf \{ a(t) \mid a \in \Cal A_h, a \geq b\}$, $c_1= \sup \{a(t) \mid a\in \Cal A_h, a \leq b\}$, we first show that 
condition $(1.2.1)$ implies $c_0=c_1$.  
For if not, then the maps $\psi_i: \Cal A + \Bbb C b\rightarrow \Bbb C$ defined by $\psi_i(y+ \alpha b)=y(t)+ \alpha c_i$, $i=0, 1$, $y\in \Cal A$, $\alpha\in \Bbb C$,  
are well defined, linear and positive; thus $\|\psi_i\|=1$ and by Hahn-Banach each $\psi_i$ can be extended to a norm-1 linear functional  
$\varphi_i:\Cal M \rightarrow \Bbb C$; we have thus obtained two states $\varphi_0, \varphi_1$ on $\Cal M$,  
which extend the pure state $t$ and are distinct (because $\varphi_0(b)\neq \varphi_1(b)$), contradicting $(1.2.1)$. Let now $\varepsilon > 0$ 
and for each $t\in \Omega$ denote $c_t=\inf \{a(t) \mid a\in \Cal A_h, a\geq b\} = \sup \{a(t) \mid a\in \Cal A_h, a \leq b\}$. Let $a^{\pm}_t\in \Cal A_h$
be such that $a_t^+ \geq b \geq a_t^-$ and $c_t + \varepsilon/2 > a^+_t(t),  a^-_t(t) > c_t - \varepsilon/2.$ By the continuity of $a_t^{\pm}\in \Cal A=C(\Omega)$ 
as a function on $\Omega$, 
there exists an open-closed neighborhood $\Omega_t$ of $t$ in $\Omega$ such that $c_t + \varepsilon/2 > a^+_t(t'), a^-_t(t') > c_t - \varepsilon/2,$ 
$\forall t'\in \Omega_t$. Thus, if we denote by $p_t\in C(\Omega)$ the characteristic function of $\Omega_t$, then $p_t$ is a 
projection in $\Cal A$ satisfying 
$$
(c_t + \varepsilon/2) p_t \geq a_t^+ p_t \geq p_t b p_t \geq a_t^-p_t \geq (c_t - \varepsilon/2)p_2. 
$$
In particular, $\|p_tbp_t-c_t p_t\|\leq \varepsilon$. Since $\Omega$ is compact, there exist $t_1, ..., t_n\in \Omega$ such that 
$\cup_i \Omega_{t_i}= \Omega$. If we now take $q_1$ to be the characteristic function of $\Omega_{t_1}$ and for each 
$j\geq 2$, $p_j$ to be the characteristic function of $\Omega_j \setminus \cup_{i=1}^{j-1} \Omega_i$, viewed as a projection in $\Cal A$, 
it follows that $\|\Sigma_j q_j b q_j - \Sigma_j c_{t_j}q_j\|\leq \varepsilon$ with $\Sigma_j c_{t_j}q_j \in C_{\Cal A}(b)$. 

\hfill 
$\square$ 
\vskip .05in 
\noindent
{\bf 1.3. Remark.} The above proof actually shows that properties $(1.2.1)-(1.2.4)$ are equivalent 
for any given element $x\in \Cal M$. More precisely, we have proved the following ``local'' statement: Let $\Cal A$ be a MASA 
in a von Neumann algebra $\Cal M$ and let $x\in \Cal M$. The following properties are equivalent:  
\vskip .05in 
\noindent 
$(1.3.1)$ Any two state extensions on $\Cal M$ of a pure state on $\Cal A$ 
coincide at $x$. 
\vskip .05in 
\noindent 
$(1.3.2)$ $C_{\Cal A}(x) \cap \Cal A \neq \emptyset$. 
\vskip .05in 
\noindent
$(1.3.3)$ $C_{\Cal A}(x) \cap \Cal A$ is a single point set $\{E_{\Cal A}(x)\}$.  
\vskip .05in 
\noindent
$(1.3.4)$ For all $\varepsilon > 0$, there exists a finite partition of $1$ with projections $q_k \in \Cal A$ such that 
$\text{\rm d}(\Sigma_k q_k x q_k, \Cal A) \leq \varepsilon$. 
\vskip .05in 
\noindent
$(1.3.5)$ There exists a unique element $E(x) \in \Cal A$ such that for all $\varepsilon > 0$,  
there exists a finite partition of $1$ with projections $q_k \in \Cal A$ such that 
$\|\Sigma_k q_k x q_k - E(x)\|\leq \varepsilon$. 
\vskip .05in
Moreover, if these conditions are satisfied for $x$, then $E(x)=E_{\Cal A}(x)$ and the following additional properties hold true:  

\vskip .05in 
\noindent
$(i)$ $\overline{C_{\Cal A}(x)}^w \cap \Cal A =\{E_{\Cal A}(x)\}$. 
\vskip .05in 
\noindent
$(ii)$ Any conditional expectation $\Cal E$ of $\Cal M$ onto $\Cal A$ 
(which always exist because $\Cal A$ is abelian)  satisfies $\Cal E(x)=E_{\Cal A}(x)$. 
\vskip .05in 
\noindent
$(iii)$ Any extension of a pure state $\psi$ on $\Cal A$ to a state $\varphi$ on $\Cal M$, 
satisfies  $\varphi(x)=\psi(E_{\Cal A}(x))$.

\vskip .05in 
\noindent
{\bf 1.4. Definitions.} Let $\Cal M$ be a von Neumann algebra and $\Cal A \subset \Cal M$  a MASA in $\Cal M$. 
We will use the following terminology: 
\vskip .05in 
\noindent
$(1.4.1)$  $\Cal A\subset \Cal M$ satisfies the {\it Kadison-Singer} (abbreviated {\it KS}) {\it  property} 
if   $(1.2.1)$ is satisfied. Condition 
$(1.2.4)$ is referred to as the {\it paving property} for $\Cal A\subset \Cal M$ (the term was coined in [A2]). Also, condition  
$(1.2.3)$ is called the {\it relative Dixmier property} for $\Cal A\subset \Cal M$, 
because of its relation to a phenomenon first emphasized in [D2] (the ``Dixmier averaging by unitaries''). Note that by Theorem 1.2 
these three properties for $\Cal A\subset \Cal M$ are actually equivalent, and they imply  $1.2 (i)-(iii)$ as well.    
\vskip .05in 
\noindent
$(1.4.2)$ An element $x\in \Cal M$ {\it can be paved} (over $\Cal A$) if  condition $(1.3.5)$  is satisfied. 
A set $X\subset \Cal M$ can be paved if each $x\in X$ can be paved.  If $x\in \Cal M$ can be paved and $\varepsilon > 0$, then we 
denote by $\text{\rm n}(\Cal A \subset \Cal M; x, \varepsilon)$ 
(or simply $\text{\rm n}(x,\varepsilon)$ if no confusion is possible) the smallest number $n$ for which there exists a partition of $1$ 
with $n$ projections $p_1, ..., p_n \in \Cal A$, such that 
$\|\Sigma_{i=1}^n p_i x p_i - E_{\Cal A}(x)\|\leq \varepsilon \|x - E_{\Cal A}(x)\|$, where $E_{\Cal A}(x)$ is given by Remark 1.3. 

More generally, if $\Cal E: \Cal M \rightarrow \Cal A$ is a conditional expectation, $x\in \Cal M$ and $\varepsilon > 0$, 
then we say that $x$ can be $\varepsilon$-{\it paved with respect to $\Cal E$}, if there exists a finite partition with projections 
$p_1, ..., p_n\in \Cal A$ such that $\|\Sigma_i p_i x p_i -\Cal E(x)\| \leq \varepsilon \|x-\Cal E(x)\|$ and we denote by 
$\text{\rm n}(\Cal E; x, \varepsilon)$ the smallest number of such projections. If there exists no such finite partition, then we let  
$\text{\rm n}(\Cal E; x, \varepsilon)=\infty.$ One should note that if $\Cal E': \Cal M \rightarrow \Cal A$ is another expectation, then 
$\|\Cal E'(x)-\Cal E(x)\|\leq \varepsilon \|x-\Cal E(x)\|$ and $\|\Sigma_i p_i x p_i -\Cal E'(x)\|\leq 2\varepsilon \|x-\Cal E'(x)\|$, 
so we have $\text{\rm n}(\Cal E; x, \varepsilon) \geq \text{\rm n}(\Cal E'; x, 2 \varepsilon)$. Thus, taking one expectation or another 
doesn't really change the nature of the function $\text{\rm n}( \cdot ; x, \varepsilon)$ and they are all ``comparable'' to 
$\text{\rm n}(\text{\rm d}; x, \varepsilon)$, which is by definition the smallest $n$ for which there exists a partition of $1$ with 
projections $p_1, ...., p_n \in \Cal A$ such that $\text{\rm d}(\Sigma_i p_i xp_i, \Cal A) \leq \varepsilon \text{\rm d}(x, \Cal A)$. In case it is clear from the context 
what expectations we take, then $\Cal E$ will not be mentioned, and we just use the notation $\text{\rm n}(\Cal A \subset \Cal M; x, \varepsilon)$. 
Also, in case there exists a normal conditional expectation of $\Cal M$ 
onto $\Cal A$ (e.g. when $\Cal M=M$ is a finite von Neumann algebra), then $\varepsilon$-pavings are always considered with respect to this expectation. 
\vskip .05in 
\noindent
$(1.4.3)$  A set $X\subset \Cal M$  has the {\it uniform paving property} (over $\Cal A$) if  it can be paved and if $\text{\rm n}(\Cal A\subset \Cal M; X, \varepsilon) 
\overset \text{\rm def} \to = \sup \{\text{\rm n}(x, \varepsilon) \mid x\in X\} $ is finite, $\forall \varepsilon > 0$.  If this holds true for $X=\Cal M$, we  
say that $\Cal A \subset \Cal M$ has the uniform paving property and  use the notation $\text{\rm n}(\Cal A\subset \Cal M; \varepsilon)$ 
for $\text{\rm n}(\Cal A\subset \Cal M; \Cal M, \varepsilon)$. We call  
this function the {\it paving size} of $\Cal A\subset \Cal M$. We will be interested in the order of magnitude of the (decreasing) functions 
$\text{\rm n}(\Cal A\subset \Cal M ; x, \varepsilon)$, 
i.e. up to the equivalence relation $f(\varepsilon) \sim g(\varepsilon)$ for functions $f, g$ requiring the existence of positive constants $ 0 < c < C < \infty$ such that 
$c \leq f(\varepsilon)/g(\varepsilon) \leq C$, $\forall \varepsilon>0$. As we will see below, the uniform paving property appears naturally 
in this context, being often equivalent to the usual paving property (notably in the case $\Cal D \subset \Cal B$), a fact first pointed out by Anderson in [A1], [A2]. 
\vskip .05in 

\noindent 
$(1.4.4)$ Let $\Cal D=\ell^\infty \Bbb N$ be the diagonal MASA in the algebra $\Cal B = \Cal B(\ell^2\Bbb N)$ 
of all linear bounded operators on the Hilbert space $\ell^2\Bbb N$. It is easy to see that the conditional expectation $\Cal B \ni  (\alpha_{jk})_{j, k \in \Bbb N} 
\mapsto (\alpha_{kk})_k \in \Cal D$ is the 
unique conditional expectation of $\Cal B$ onto $\Cal D$ and that it is normal. We use the terminology ``{\it the classic Kadison-Singer problem}'' 
for the question of whether $\Cal D \subset \Cal B$ has the KS property. 
By Theorem 1.2, this property is equivalent to the paving property for $\Cal D \subset \Cal B$. 
The terminology ``{\it Kadison-Singer conjecture}'' is sometimes used for the statement predicting that the KS property does hold true 
for this inclusion, despite the fact that, in their paper,  Kadison and Singer expressed the belief that 
the property doesn't actually hold true for $\Cal D\subset \Cal B$....   

\vskip .05in 
The next result summarizes some well known paving properties, notably J. Anderson's observations    
that uniform paving for $\Cal D \subset \Cal B$ is equivalent to   
paving  and that the classic Kadison-Singer problem is equivalent to  $D_k\subset M_{k \times k}(\Bbb C)$ having uniformly bounded 
paving size (see [A1], [A2]).  

\proclaim{1.5. Proposition} $0^\circ$ Let $\Cal A$ be a MASA in the von Neumann algebra $\Cal M$, $p\in \Cal P(\Cal A)$ a projection 
and $\Cal A \subset \Cal N \subset \Cal M$ an intermediate von Neumann algebra. If $x\in p\Cal Mp$ $($respectively $x\in \Cal N)$ 
then $\text{\rm n}(\Cal Ap \subset p\Cal Mp; x, \varepsilon) = \text{\rm n}(\Cal A \subset \Cal M; x, \varepsilon)$ 
$($resp.  $\text{\rm n}(\Cal A \subset \Cal N; x, \varepsilon) = \text{\rm n}(\Cal A \subset \Cal M; x, \varepsilon))$.

\vskip .05in 
$1^\circ$ Let $\{\Cal A_i \subset \Cal M_i\}_i$ be a family of MASAs in von Neumann algebras and 
for each $i$ let $x_i \in \Cal M_i$, 
$\|x_i\|\leq 1$. Denote $\Cal A = \oplus_i \Cal A_i$, $\Cal M= \oplus_i \Cal M_i$, $x=\oplus_i x_i \in \Cal M$.  
Then $\text{\rm n}(\Cal A \subset \Cal M; x, \varepsilon) = \sup_i \text{\rm n}(\Cal A_i \subset \Cal M_i; x_i, \varepsilon)$ 
and $\text{\rm n}(\Cal A \subset \Cal M; \varepsilon) = \sup_i \text{\rm n} (\Cal A_i \subset \Cal M_i; \varepsilon)$, $\forall \varepsilon > 0$. 

\vskip .05in 
$2^\circ$ If a MASA $\Cal A$ in a von Neumann algebra  $\Cal M$ has the property that there exists 
a sequence of mutually orthogonal projections $p_n \in \Cal A$ with embeddings 
$\theta_n: \Cal M \hookrightarrow p_n\Cal Mp_n$ such that $\theta_n(\Cal A)=\Cal Ap_n$, $\forall n$, then 
$\Cal A\subset \Cal M$ 
has the paving property iff it has the uniform paving property. 

\vskip .05in 
$3^\circ$ The diagonal MASA $\Cal D = \ell^\infty \Bbb N$ in the algebra of all linear bounded operators $\Cal B = \Cal B(\ell^2\Bbb N)$ 
on the Hilbert space $\ell^2\Bbb N$, has the paving property iff it has the uniform paving property. Moreover,  
$\text{\rm n}(\Cal D \subset \Cal B;  \varepsilon)= \sup_k \text{\rm n}(D_k \subset M_{k \times k}(\Bbb C); \varepsilon)$, $\forall \varepsilon > 0$. 
\vskip .05in 
$4^\circ$  If $\Cal A$ is a MASA in a von Neumann algebra $\Cal M$, with $\Cal E:\Cal M \rightarrow \Cal A$ an expectation, and $1> \varepsilon > 0$, then we have 
$\sup \{\text{\rm n}(\Cal E; x, \varepsilon^2) \mid x\in \Cal M\} \leq (\sup \{\text{\rm n}(\Cal E; y, \varepsilon) \mid y\in \Cal M \})^2$. Thus, in order for $\Cal A \subset \Cal M$ to have the 
uniform paving property, it is sufficient that for some $\varepsilon < 1$ we have $\sup \{\text{\rm n}(\Cal E; y, \varepsilon) \mid y\in \Cal M \} < \infty$. 
\endproclaim
\noindent
{\it Proof}. Parts $0^\circ$ and $1^\circ$ are trivial and $2^\circ$ is an immediate consequence of $1^\circ$. Then $2^\circ$ implies the equivalence in the first part of $3^\circ$. 

To establish the formula in $3^\circ$, note that if a sequence of projections $q_n$ in $ \ell^\infty\Bbb N$ is convergent in the weak operator  topology 
to some element $q\in \ell^\infty \Bbb N$, then $q$ is itself a projection and $q_n$ converges to $q$  in the strong operator topology 
as well.  

The inequality $\text{\rm n}(\Cal D \subset \Cal B;  \varepsilon)\geq  \sup_k \text{\rm n}(D_k \subset M_{k \times k}(\Bbb C); \varepsilon)$ 
is trivial because the right hand side is equal to $\text{\rm n}(\oplus_k D_k \subset \oplus_k M_{k \times k}(\Bbb C); \varepsilon)$ 
and one can embed $\oplus_k M_{k \times k}(\Bbb C)$ into $\Cal B$ in a way that takes $\oplus_k D_k$ onto $\Cal D$. 

For the inequality $\leq$ 
let $T\in \Cal B$ be so that $\|T\|\leq 1$ and $T$ has $0$ on the diagonal. Let  
$P_k \in \Cal D=\ell^\infty \Bbb N$ be the projection onto the first $k$ coordinates. Let $\{p_{k,j}\}_j$ be a partition of $P_k$ 
into $n = \sup_k \text{\rm n}(D_k \subset M_{k \times k}(\Bbb C); \varepsilon)$ projections such that $\|\Sigma_j p_{k,j} T p_{k,j} \|\leq \varepsilon$. 
Let $k_1 < k_2 < .... $ be a subsequence such that $\{p_{k_m,j}\}_m$ is weakly convergent for each $j=1, 2, ..., n$ and denote by $q_j$ the 
corresponding weak limit. By the above observation, for each $j$,  $\{p_{k_m,j}\}_m$ converges in fact in the strong operator topology and $q_j$ is a projection. 
Also, since $\Sigma_j p_{k_m,j} = P_{k_m}$ 
and $P_{k_m}$ $so$-converges to $1_{\Cal B}$, it follows that $\Sigma_j q_j = 1$ as well. Thus, $\{q_j\}_j$ is a partition of $1$ by $n$ projections and since 
$\{\Sigma_j p_{k_m, j} T p_{k_m,j} \}_m$ is $so$-convergent to $\Sigma_j q_j T q_j$ and the operator norm is inferior semicontinuous with respect to the 
$so$-convergence, it follow that $\|\Sigma_j q_j T q_j \| \leq \limsup_m \|\Sigma_j p_{k_m, j} T p_{k_m,j}\| \leq \varepsilon.$

Finally, part $4^\circ$ is immediate from the definitions. 
\hfill 
$\square$

\vskip .05in 
\noindent
{\bf 1.6. Remark.} While the classic Kadison-Singer problem is still open, one should point out 
that a large number of beautiful paving results have been obtained over the years, showing that the equivalent 
conditions 1.3 are satisfied for many classes of operators $x\in\Cal B(\ell^2 \Bbb Z)$. Thus,  it is shown in [A2] 
that if $x$ is in the C$^*$-algebra for the reduced C$^*$-algebra of the group $\Bbb Z$,  $C^*_r(\Bbb Z)$, i.e. in the operator norm-closure 
of the span of the range of the left regular representation $\lambda$ of the group $\Bbb Z$, then $x$ can be paved. In [BeHKW] it is shown that 
matrices with non-negative entries in $M_{n \times n}(\Bbb C)$ can be paved,  while in  [BT] it is shown that if an element $x$ in the weak closure $L(\Bbb Z)\simeq L^\infty(\Bbb T)$ of $C^*_r(\Bbb Z)$ in $\Cal B(\ell^2 \Bbb Z)$  has Fourier coefficients satisfying certain growth properties, then $x$ can be paved. Also, a number of results 
have been obtained in [AkA], [A1], [A2], [CaFTW], [We], etc, showing that in order to solve the paving 
conjecture, it is sufficient to be able to pave certain particular  classes of elements (e.g. projections with small diagonal entries in [AkA]). 

\heading 2. Kadison-Singer in II$_1$ factor framework
\endheading

We prove in this section that the KS property for the inclusion of the diagonal MASA $\Cal D=\ell^\infty \Bbb N$ into the algebra $\Cal B = \Cal B(\ell^2 \Bbb N)$,  
of all linear bounded operators on the Hilbert space $\ell^2 \Bbb N$, is equivalent to the KS property of  
MASAs in II$_1$ factors obtained as ultraproducts of certain Cartan inclusions. Also, we use a  dilation trick to prove  that in order to pave arbitrary elements 
in an ultraproduct of inclusions of MASAs, it is sufficient to pave projections that expect on scalars. 

From now on, we fix once for all an (arbitrary) 
free ultrafilter $\omega$ on $\Bbb N$.  All finite von Neumann algebras that we consider are assumed equipped with a faithful normal trace state, 
generically denoted by $\tau$ (unless otherwise specified). 

If $M_n$, $n\geq 1$, is a sequence of finite von Neumann algebras then, we denote by $\Pi_\omega M_n$ their $\omega$-ultraproduct, i.e., 
the finite von Neumann algebra obtained as the quotient of $\oplus_n M_n$ by its ideal 
$\Cal I_\omega=\{ (x_n) \mid \lim_\omega \tau(x_n^*x_n) = 0\}$, endowed with the trace $\tau(y) =\lim_\omega \tau(y_n)$, where $(y_n)_n \in \oplus_n M_n$ 
is in the class $y \in \oplus_n M_n/\Cal I_\omega$ ([W]). Recall that if $M_n$ are factors and $\dim M_n \rightarrow \infty$, 
then $\Pi_\omega M_n$ is a II$_1$ factor ([W], [F]) and that if $A_n\subset M_n$ are MASAs, $n\geq 1$, then $\Pi_\omega A_n$ is a MASA in $\Pi_\omega M_n$ 
(see e.g. [P1]). If $A\subset M$ is a MASA in a finite von Neumann algebra, then $A^\omega \subset M^\omega$ denotes its $\omega$-ultrapower, 
i.e. the ultraproduct of infinitely many copies of $A\subset M$. Note that $M$ naturally embeds into $M^\omega$, as 
the von Neumann subalgebra of constant sequences. 

Recall that a  {\it Cartan subalgebra} $A$ in a finite von Neumann algebra $M$ is a MASA in $M$ 
whose normalizer  $\Cal N_M(A)=\{u\in \Cal U(M)\mid uAu^*=A\}$ generates $M$, i.e. $\Cal N(A)''=M$ (see [FM]).
\vskip .05in
\noindent 
{\bf 2.1. Notations} $(a)$ We denote the  
inclusion $\Pi_\omega D_n \subset \Pi_\omega M_{n \times n}(\Bbb C)$ by $\text{\bf D}(\omega) \subset \text{\bf M}(\omega)$, 
or simply $\text{\bf D}\subset \text{\bf M}$. Note that given any sequence of MASAs 
$A_n\subset M_{n\times n} (\Bbb C)$, the von Neumann algebra $\Pi_n A_n \in \text{\bf M}$ is unitary conjugate to $\text{\bf D}$ in $\text{\bf M}$. 
One should point out that $\text{\bf D}$ is not a Cartan subalgebra in $\text{\bf M}$, in fact $\text{\bf M}$ has no Cartan subalgebras  (cf. [P2]; see 
$2.3.2^\circ$ below). Also, $\text{\bf D}, \text{\bf M}$ are non-separable (cf. [F]) and $\text{\bf M}$ is non-amenable (because it contains $L(\Bbb F_2)$, 
as first noticed in [Wa]). 

$(b)$ We represent the hyperfinite II$_1$ factor $R$ as the infinite tensor product $\overline{\otimes}_n (M_{2 \times 2}(\Bbb C), tr)_n$, where $tr$ is the 
normalized trace on $M_{2 \times 2}(\Bbb C)$. Also, we denote by $D\subset R$ the Cartan subalgebra obtained as the infinite tensor product 
of the diagonals $D_2 \subset M_{2\times 2}(\Bbb C)$. Recall that  any other Cartan subalgebra  $A\subset R$ is conjugate to $D$ by an 
automorphism of $R$ (cf [CFW]).  Thus, if $D^\omega \subset R^\omega$ is the $\omega$-ultrapower of $D\subset R$, then 
any ultraproduct $\Pi_\omega A_n \subset R^\omega$, with $A_n\subset R$ Cartan subalgebras,  is conjugate to $D^\omega$ by an automorphism 
$\theta = (\theta_n)_n$ of $R^\omega$, where $\theta_n\in \text{\rm Aut}(R)$ is so that $\theta_n(A_n)=D$.  We denote by $\text{\bf R}\subset R^\omega$ 
the von Neumann algebra $D^\omega \vee R$, generated by $D^\omega$ and $R$, or equivalently by $D^\omega$ and $\Cal N_R(D)$.  

If $\Gamma \subset \Cal N_R(D)$ is any countable subgroup generating the hyperfinite equivalence relation $\Cal R$ associated with $D\subset R$ (cf [FM]), then 
$\text{\bf R}$ is generated by $D^\omega$ and $\Gamma$. Moreover, if $\Gamma$ acts freely on $D$, then $\Gamma$ acts freely on $D^\omega$ 
as well and so we can view $\text{\bf R}$ as the crossed product $D^\omega \rtimes \Gamma$. Finally, note that $\text{\bf R}$ is an amenable 
II$_1$ von Neumann algebra, but not a factor, in fact any sequence $(a_n)_n \in D^\omega$ with $a_n \in 1 \otimes_{j \geq k_n} (D_2)_j$ 
for some $k_n \rightarrow \infty$, lies in $R'\cap D^\omega = \text{\bf R}'\cap D^\omega=\Cal Z(\text{\bf R})$ (the center of $\text{\bf R}$). 

Note that, while $D^\omega$ is Cartan in $\text{\bf R}$, $D^\omega$ is not Cartan in $R^\omega$, in fact $R^\omega$ has no Cartan subalgebras (by [P2]; 
see $2.3.2^\circ$ below).

\proclaim{2.2. Theorem} $\Cal D \subset  \Cal B$ has the KS property $($equivalently, the paving property$)$ 
if and only if $\text{\bf D} \subset \text{\bf M}$ $($resp. $D^\omega \subset R^\omega$, resp. $D^\omega \subset \text{\bf R})$ 
has this property. Moreover,  all these inclusions have the same paving size $($whether finite or infinite$)$: 
$$
\text{\rm n}(\Cal D \subset \Cal B; \varepsilon)=\text{\rm n}(\text{\bf D}\subset \text{\bf M}; \varepsilon) 
=\text{\rm n}(D^\omega \subset R^\omega; \varepsilon)=\text{\rm n}(D^\omega \subset \text{\bf R}; \varepsilon). \tag 2.2.1
$$

They also have the same paving size as the Cartan subalgebra inclusions $\text{\bf D}\subset \Cal N_{\text{\bf M}}(\text{\bf D})''$ and 
$D^\omega \subset \Cal N_{R^\omega}(D^\omega)''$.

\endproclaim
\noindent 
{\it Proof}. Consider the inclusion $A_0=\oplus_{n=1}^\infty D_{n}\subset \oplus_{n=1}^\infty  M_{n\times n}(\Bbb C) = M_0$ and note that 
by $1.5.1^\circ$ and $1.5.3^\circ$ we have 
$$
\text{\rm n}(\Cal D \subset \Cal B; \varepsilon)=\sup_n \text{\rm n}(D_n \subset M_{n\times n}(\Bbb C); \varepsilon)
=\text{\rm n}(A_0 \subset M_0; \varepsilon).
$$ 

Embed now $A_0$ into $D$ and then extend this to an embedding of $M_0$ into $R$ so that the matrix units of
each direct summand $M_{n \times n}(\Bbb C)$ are in the normalizing groupoid of $D\subset R$ 
(this is possible because $D$ is Cartan in $R$; in fact, semiregularity is sufficient). Note that this implies $M_0$ and $D$ make a commuting square, 
i.e. $E_{M_0}E_D=E_D E_{M_0}=E_{A_0}$. Also, we trivially have 

$$
\text{\rm n}(D^\omega \subset R^\omega; \varepsilon)\geq \text{\rm n}(D^\omega \subset \text{\bf R}; \varepsilon)\geq 
\text{\rm n}(D^\omega \subset \text{\bf R}; R, \varepsilon)\geq \text{\rm n}(D^\omega \subset \text{\bf R}; M_0, \varepsilon).  
$$

Let $x=(x_n)_n \in M_0\subset R$ be so that $E_{D^\omega}(x)=E_{A_0}(x)=0$ 
and note that $\text{\rm n}(D^\omega \subset \text{\bf R}; x, \varepsilon) = \sup_n \text{\rm n}(D^\omega \subset \text{\bf R}; x_n, \varepsilon)$.   
Let $s_n\in D\subset D^\omega$ denote the support projection of $M_{n\times n}(\Bbb C)$ in $A_0\subset D$.  
Each $D^\omega s_n \subset (M_0 \vee D^\omega)s_n$ is 
of the form $C(\Omega) \subset C(\Omega)\otimes M_{n \times n}(\Bbb C)$. Thus, 
if $p_1, ..., p_m \in \Cal P(D^\omega)$ is a partition of $1$ such that $\|\Sigma_i p_i x p_i \|\leq \varepsilon$,   
then the evaluation at any point $t\in \Omega$ of $p_i s_n$, $1\leq i \leq m$,   
gives a partition of $1_{M_{n \times n}(\Bbb C)}$ with $m$ projections $q_i$ in $D_n$ 
such that $\| \Sigma_{i=1}^m q_i x_n q_i\|\leq \varepsilon$. This shows that 
$\text{\rm n}(D^\omega \subset \text{\bf R}; M_0, \varepsilon)\geq \text{\rm n}(A_0\subset M_0; \varepsilon)$. Since the latter 
is equal to $\text{\rm n}(\Cal D \subset \Cal B; \varepsilon)$ and to 
$\sup_n \text{\rm n}(D_n \subset M_{n\times n}(\Bbb C); \varepsilon)$, in order to 
end the proof of the fact that $\Cal D\subset \Cal B$, $D^\omega \subset R^\omega$, $D^\omega \subset \text{\bf R}$ 
(as well as $A_0\subset M_0$) have the same paving size, it is sufficient to show that $\sup_n \text{\rm n}(D_n \subset M_{n\times n}(\Bbb C); \varepsilon)  
\geq \text{\rm n}(D^\omega \subset R^\omega; \varepsilon)$.  

To this end, let $x=(x_n)_n \in R^\omega \ominus D^\omega$ and note that one can take each $x_n$ 
to belong to $R\ominus D$ and such that $\|x_n \|\leq \|x\| + c_n$, 
$\forall n$, for some $c_n \rightarrow 0$. Moreover, since there exists an increasing sequence of 
$2^k \times 2^k$ matrix subalgebras $M_{2^k}\subset R$ with diagonal subalgebra $D_{2^k}\subset D$ making a commuting square with 
$D\subset R$ such that $\overline{\cup_k M_{2^k}}^w=R$ and $\overline{\cup_k D_{2^k}}^w=D$, we may replace each $x_n$ by $E_{M_{2^{k_n}}}(x_n)$, 
and thus assume 
$x_n\in M_{2^{k_n}}\ominus D_{2^{k_n}}$, $\forall n$. Let $\{q^n_j\}_j \subset D_{2^{k_n}}$ be a partition of $1$ with 
$K=\sup_n \text{\rm n}(D_n \subset M_{n\times n}(\Bbb C); \varepsilon) $ projections such that $\|\Sigma_{j=1}^K q^n_j x_n q^n_j \| 
\leq \varepsilon \|x_n \|$. If we denote by $q_j=(q^n_j)_n \in D^\omega$, 
it follows that $\|\Sigma_{j=1}^K q_j x q_j\| \leq \varepsilon \|x\|$. This proves the desired inequality.

Let now $\varepsilon>0$ and assume 
$m=\text{\rm n}(\Cal D \subset \Cal B; \varepsilon)=\sup_n \text{\rm n}(D_n \subset M_{n\times n}(\Bbb C); \varepsilon)$ is finite. 
Any $x \in \text{\bf M}$ with $\|x\|\leq 1$ and $E_{\text{\bf D}}(x)=0$ can be represented by a sequence $x=(x_n)_n$ with $x_n\in M_{n \times n}(\Bbb C)$ 
such that $\|x_n\|\leq 1+c_n$, $E_{D_n}(x_n)=0$, for some $c_n \rightarrow 0$. For each $n$ there exists a partition of $1$ with projections $p^n_j, 1\leq j \leq m$, such that 
$\|\Sigma_{j=1}^m p^n_j x_n p^n_j \|\leq \varepsilon (1+c_n)$. But then $p_j=(p^n_j)_n\in \text{\bf D}$ gives a partition of $1$ satisfying 
$\|\Sigma_{j=1}^m p_j x p_j\|\leq \varepsilon$. Thus, $\sup_n \text{\rm n}(D_n \subset M_{n\times n}(\Bbb C); \varepsilon)
\geq \text{\rm n}(\text{\bf D}\subset \text{\bf M}; \varepsilon)$. 

Conversely, assume $m=\text{\rm n}(\text{\bf D}\subset \text{\bf M}; \varepsilon)$ is finite.  
Let $x$ be an element in $M_{k \times k}(\Bbb C)$, for some $k\geq 1$, with $\|x\|\leq 1$, $E_{D_k}(x)=0$. For each $n$ larger than $k$, 
embed $M_{k\times k}(\Bbb C)$ into $M_{n\times n}(\Bbb C)$ by first letting $n=k d_n + r_n$, with $d_n, r_n \in \Bbb N$, $r_n < k$, 
then letting $s_j^n\in D_n$ 
be mutually orthogonal projections of trace $k/n$, and then identifying $D_k\subset M_{k \times k}(\Bbb C)$ with $D_n s^n_j \subset 
s_j^n M_{n \times n}(\Bbb C)s_j^n$ via some isomorphism $\theta^n_j$, for each $j=1, ..., d_n$, 
and mapping diagonally 

$$
M_{k\times k}(\Bbb C) \ni y \mapsto \theta^n(y)\overset \text{\rm def} \to = \Sigma_j \theta^n_j(y) \in 
\Sigma_{j=1}^{d_n}s^n_jM_{n \times n}(\Bbb C)s^n_j \subset M_{n\times n}(\Bbb C).  
$$  
Then consider the embedding $\theta: M_{k \times k}(\Bbb C) \rightarrow \text{\bf M}$, by $\theta(y)=
(\theta^n(y))_n$. Let $p_1, ..., p_m\in \Cal P(\text{\bf D})$ be a partition of $1$ such that $\|\Sigma_i p_i \theta(x) p_i\|\leq \varepsilon$. 
One can then choose representing sequences $p_i=(p^n_i)_n$, with $p^n_i \in \Cal P(D_n)$, $1\leq i \leq m$, a partition of $1$ for each $n$. 
We claim that for any $\delta > 0$ there exist $n$ and $j \in \{1, ..., d_n\}$, such that $\| \Sigma_i p^n_is^n_j \theta_j(x) p^n_js^n_j\| <  \varepsilon +\delta$. 

Indeed, for if not then for every  $n$ and $j=1, ..., d_n$ the spectral projection of $|\Sigma_i p^n_is^n_j \theta_j(x) p^n_js^n_j|$ corresponding to the 
interval $[\varepsilon +\delta, 1]$ is non-zero, thus having trace at least $1/n$. Since 
$|\Sigma_i p^n_i \theta^n(x) p^n_i|=\Sigma_{j=1}^{d_n} |\Sigma_i p^n_is^n_j \theta^n_j(x) p^n_is^n_j|$, the spectral projection corresponding to $[\varepsilon +\delta, 1]$ 
of $|\Sigma_i p^n_i \theta^n(x) p^n_i|$ has trace $\geq d_n/n$. But this implies that the spectral projection corresponding to $[\varepsilon +\delta/2, 1]$ of 
$|\Sigma_i p_i \theta(x) p_i |=(|\Sigma_i p^n_i \theta^n(x) p^n_i|)_n\in \text{\bf M}$ has trace $\geq \lim_n d_n/n=1/k$. Thus 
$\|\Sigma_i p_i \theta(x) p_i\|\geq \varepsilon + \delta/2$, a contradiction. 

If we now choose some  $n$ and $j\in \{1, ..., d_n\}$ satisfying $\|\Sigma_i p^n_is^n_j \theta^n_j(x) p^n_js^n_j\| <  \varepsilon +\delta$, 
and let  $q_i, 1\leq i \leq m$, be the  pre-image in $D_k$ of the partition $\{p^n_is^n_j\}_i$ via $\theta^n$, then $\|\Sigma_{i=1}^m q_i x q_i\| < \varepsilon + \delta$. 
Letting $\delta=1/n$, $n = 1, 2, ...$, we obtain a sequence of partitions of $1$ by projections $\{q_1(n), ..., q_m(n)\}_n$ in $D_k$ satisfying  
$\|\Sigma_{i=1}^m q_i(n) x q_i(n)\| < \varepsilon + 1/n$. But the unit ball of $D_k$ is compact in the operator norm, so by taking the limit over some subsequence, 
we get a partition of $1$ with projections $q_1, ..., q_m \in D_k$ with $\|\Sigma_i q_i x q_i\|\leq \varepsilon$.  Thus, $m\geq \text{\rm n}(D_k \subset M_k; \varepsilon)$, 
and since $k$ was arbitrary, $m \geq \sup_k \text{\rm n}(D_k \subset M_k; \varepsilon)$. 

Finally, since $D^\omega \subset  \text{\bf R} \subset \Cal N(D^\omega)'' \subset R^\omega$, we have 

$$
\text{\rm n}(D^\omega \subset \text{\bf R} ; \varepsilon) 
\leq \text{\rm n}(D^\omega \subset \Cal N_{\text{\bf M}}(\text{\bf D})''; \varepsilon) \leq \text{\rm n} (D^\omega \subset R^\omega; \varepsilon), 
$$  
and since the first and last terms are equal, they must all be equal. Similarly, $\text{\bf D}\subset \Cal N(\text{\bf D})''\subset 
\text{\bf M}$ implies $\text{\rm n}(\text{\bf D}\subset \Cal N(\text{\bf D})''; \varepsilon) \leq \text{\rm n}(\text{\bf D}\subset \text{\bf M}; \varepsilon)$, 
while arguments above show that $\sup_n \text{\rm n} (D_n \subset M_{n \times n}(\Bbb C); \varepsilon) \leq \text{\rm n}(\text{\bf D}\subset \Cal N(\text{\bf D})''; \varepsilon)$, 
with the first of these terms equal to $\text{\rm n}(\text{\bf D}\subset \text{\bf M}; \varepsilon)$.

\hfill 
$\square$ 

If $M$ is a finite von Neumann algebra with its faithful normal trace state $\tau$, then we denote by $\|x\|_2=\tau(x^*x)^{1/2}$, $x\in M$, 
the $L^2$ (or Hilbert) norm given by the trace. We denote by $L^2M$ the Hilbert space obtained by completing $M$ in this $L^2$ norm 
and view $M$ in its {\it standard representation}, as left multiplication representation on $L^2M$.  
We also use the notation $L^1M$ for the completion of $M$ in the norm $\|x\|_1=\tau(|x|)$. We view the elements 
in $L^2M$ (resp. $L^1M$) 
as square summable (resp. summable) operators affiliated with $M\subset \Cal B(L^2M)$, in the usual way. All 
self-adjoint elements affiliated with $M$ (in particular elements in $L^2M$, $L^1M$) have spectral decomposition belonging to $M$ and they can be 
multiplied. In particular, we have $L^2 M \cdot L^2 M = L^1M$. 

A finite von Neumann algebra $M$ with a normal faithful trace is {\it separable} if it is separable with respect 
to the $\| \cdot \|_2$-norm given by the trace. This condition is easily seen to be equivalent to $M$  
being countably generated. A von Neumann algebra is {\it diffuse} if it has no minimal (non zero) 
projection. Any abelian von Neumann algebra $A$ which is diffuse and separable is isomorphic to $L^\infty([0, 1])$ (or to $L^\infty (\Bbb T)$). Moreover, if $A$ 
is endowed with a faithful normal state $\tau$, then the isomorphism $A \simeq L^\infty([0,1])$ can be taken so that to carry $\tau$ 
onto the integral $\int \cdot \ \text{\rm d}\mu$, where $\mu$ is the Lebesgue measure on $[0, 1]$. 

It is well known that all separable diffuse abelian von Neumann subalgebras in an utraproduct II$_1$ factor are unitary conjugate 
(see e.g. [P2]). We will show below that any II$_1$  factor 
$M$ that has this property will automatically have several other properties, like absence of Cartan subalgebras (already 
noticed in [P2]) and the fact that in order to pave arbitrary elements over a MASA in $M$, it is sufficient to pave projections that expect on scalar 
multiples of $1$. (Note that in Anderson's formulation of the KS problem as the uniform paving property in $M_{k \times k}(\Bbb C)$, $k \nearrow \infty$, 
the reduction of the problem to paving special elements, such as projections with constant diagonal, has been subject of much study, 
see [AkA], [A2], [CaFTW], etc.)  

\proclaim{2.3. Proposition} $1^\circ$. Assume a   $\text{\rm II}_1$ factor $M$ has the property 
that given any projection $p\in M$, any two separable diffuse abelian von Neumann subalgebras of $pMp$ are unitary conjugate. Then $M$ 
satisfies the following properties 
\vskip .05in 
$(a)$ Given any MASA $A$ in $M$,  there exists a diffuse abelian von Neumann subalgebra $B_0\subset M$ perpendicular to $A$. $($Recall from $\text{\rm [P2]}$ 
that two von Neumann subalgebras $B_1, B_2$ of a finite von Neumann algebra $M$ are 
said to be perpendicular if $\tau(b_1b_2)=0$, $\forall b_i\in B_i$ with $\tau(b_i)=0$, $i=1,2.)$ 
\vskip .05in 
$(b)$ $M$ has no separable MASAs and no Cartan subalgebras.  
\vskip .05in 
$(c)$ If $A$ is a  MASA in $M$, then $A\subset M$ 
has the paving property iff any projection in $M$ that expects on a scalar multiple of a projection in $A$ can be paved. 
Moreover, if $\Cal P_0$ denotes the set of such projections, then the paving size $\text{\rm n}(\varepsilon)$ of $A\subset M$ satisfies 
$ \text{\rm n}(\varepsilon) \leq \text{\rm n}(A\subset M; \Cal P_0, \varepsilon/50)^2(\varepsilon^{-1}+1)^2$. 

\vskip .05in 
$2^\circ$ If $\{M_n\}_n$ is a sequence of  finite von Neumann factors with 
$\text{\rm dim} M_n \rightarrow \infty$, then the ultraproduct $\text{\rm II}_1$ 
factor $M=\Pi_\omega M_n$ satisfies the assumption in part $1^\circ$, i.e., 
given any projection $p\in M$, any two separable diffuse abelian von Neumann subalgebras  of $pMp$ are unitary conjugate.  
Thus, $M=\Pi_\omega M_n$ satisfies properties  $(a), (b), (c)$ as well. 
\endproclaim
\noindent 
{\it Proof}. $2^\circ$  is well known (see e.g. Lemma 7.1 in [P2]). 

$1^\circ$ Part $(a)$ is shown in the proof of Theorem 7.3 in [P2], but let us recall the argument here for completeness. 
Let $D\subset A$ be a separable diffuse von Neumann subalgebra. Since any two separable diffue abelian subalgebras 
in $M$ are unitary conjugate and since $M$ contains copies of the hyperfinite II$_1$ factor,  we may assume $D$ is the 
Cartan subalgebra of such a subfactor $R \subset M$, 
represented as $D=D_2^{\otimes \infty}\subset M_{2 \times 2}(\Bbb C)^{\otimes \infty} =R$.  
Let $D_2^0\subset M_{2 \times 2}(\Bbb C)$ be a maximal abelian subalgebra of $M_{2 \times 2}(\Bbb C)$ that is perpendicular to $D_2$ 
and denote $D^0={D^0_2}^{\otimes \infty} \subset R$. Then $D \perp D^0$ and since both $D, D^0$ are MASAs in $R$, 
we have $E_{D'\cap M}(D^0)=E_{D'\cap R}(D^0)=E_{D}(D^0)=\Bbb C$, i.e. $D^0 \perp D'\cap M \supset A$, proving $(a)$. 

To prove $(b)$, let $A$ be a MASA in $M$. If $A$ is separable, then 
it has a diffuse proper subalgebra, $A_0\subset A$, 
which cannot be unitary conjugate to $A$ because it is not a MASA. Moreover, by part $(a)$, there exist separable diffuse abelian subalgebras    
$D, D^0$ in $M$ such that $D\subset A$ and $D^0\perp A$. Let $u\in \Cal U(M)$ be so that $uDu^*=D^0$. Then $u$ is perpendicular 
to the normalizer of $A$ in $M$. Indeed, because for any $v\in \Cal N_M(A)$ and any partition $p_i\in D$ of mesh $\leq \varepsilon$, 
we have 

$$
|\tau(uv)|^2=|\tau(\Sigma_i p_i uv p_i)|^2  
\leq \|\Sigma_i p_i u v p_i\|^2_2 = \Sigma_i \tau(u^*p_iuvp_iv^*) = \Sigma_i \tau(p_i)^2 
\leq \varepsilon. 
$$ 
Since $\varepsilon >0$ was arbitrary, $\tau(uv)=0$. Thus $u \perp \Cal N_M(A)''$. 

To prove  $(c)$, note first that  for any MASA in a II$_1$ factor $M$ and any $x\in M$, the paving size over $A$ of any element $x\in M$ 
behaves well with respect to scalar translations and multiplications:  

\vskip .05in
\noindent  
$(1)$ $\text{\rm n}(x, \varepsilon)= \text{\rm n}(x + \alpha 1, \varepsilon)$,  $\text{\rm n}(\alpha x,  \varepsilon )= \text{\rm n}(x, \varepsilon/|\alpha|)$. 

\vskip .05in 
Now note that if $y\in M$ is a $\delta$-perturbation of $x\in M$, then any $\varepsilon$-paving of $y$ gives a  $\varepsilon+\delta$ paving of $x$, 
more precisely:  

\vskip .05in 
\noindent
$(2)$ If $\|x-y\|\leq 2^{-1}\delta(1+\varepsilon)^{-1} \|x-E_A(x)\|$ and $\|\Sigma_i p_i(y-E_A(y))p_i\| \leq \varepsilon \|y-E_A(y)\|$, 
then 

$$
\|\Sigma_i p_i (x-E_A(x))p_i\|  
$$
$$
\leq \|(x-y) -E_A(x-y)\| +\|\Sigma_i p_i (y-E_A(y))p_i\| 
$$
$$
\leq  \|(x-y) -E_A(x-y)\| + 
\varepsilon \|y-E_A(y)\| 
$$ 
$$
\leq (1+\varepsilon) \|(x-y) -E_A(x-y)\| + \varepsilon \|x-E_A(x)\| 
$$
$$
\leq (\varepsilon + \delta)\|x-E_A(x)\|
$$

\noindent
$(3)$ If $x_1, x_2\in M$ can be paved, then $x_1+x_2$ can be paved. More specifically, if 
$x=y_1 + i y_2$, $y_i=y_i^*$, is the decomposition of $x$ into its real and imaginary parts, 
then $\|y_i\|\leq \|x\| $ and so if $p_i$, $q_j \in \Cal P(A)$ are partitions of $1$ such that $\|\Sigma_j p^i_j y_i p_j^i -E_A(y_i)\| 
\leq \varepsilon \|y_i - E_A(y_i)\|$, $i=1,2$, then the partition $\{p_k\}_k = \{p^1_ip^2_j\}_{i, j}$ satisfies 

$$
\|\Sigma_k p_k x p_k - E_A(x)\| \leq \|\Sigma_k p_k y_1 p_k -E_A(y_1)\| +\|\Sigma_k p_k y_2 p_k -E_A(y_2)\| 
$$
$$
\leq \|\Sigma_j p^1_j y_1 p_j^1 -E_A(y_1)\|  + \|\Sigma_j p^2_j y_2 p_j^2 -E_A(y_2)\| 
$$
$$
\leq \varepsilon \|y_1 - E_A(y_1)\| + \varepsilon \|y_1 - E_A(y_1)\| \leq 2\varepsilon \|x-E_A(x)\|.
$$
Thus, $\text{\rm n}(x, 2\varepsilon) \leq \text{\rm n}(\Re x, \varepsilon)\text{\rm n}(\Im x, \varepsilon)$.

Let now $x=x^*\in M$ with $E_A(x)=0$ and $\|x\|=1$. Note that $y_0= 12^{-1}(x + 5\|x\|)$ satisfies $ 1/3 \leq y_0 \leq 1/2$, 
thus $1/3 \leq E_A(y_0) \leq 1/2$. Denote by $e_k$ the spectral projection of $E_A(y_0)$ corresponding to the interval $[1/3 + (k-1)\varepsilon/6, 
1/3 + k\varepsilon/6)$ and note that there are $\leq \varepsilon^{-1} + 1$ many such non-zero projections. If we let $t_k=1/3 + (k-1)\varepsilon/6$, 
$a_k = E_A(y_0)e_k$, then $ t_ke_k \geq a_k \geq t_{k-1}e_k$. Thus, if we denote  $b= \Sigma_k t_{k-1}^{-1} a_k$ then $1 \leq b \leq 1+ \varepsilon/2$ and  
$(1+\varepsilon/2)^{-1/2} \leq  b^{-1/2} \leq 1 $.  Notice now that $y = b^{-1/2} y_0 b^{-1/2}$ satisfies $E_A(y)=\Sigma_k t_k e_k$, $1/3 - \varepsilon/12 \leq y \leq 1/2$ and 
$\|y_0-y\| \leq 2\|1-b^{-1/2}\| \|y_0\|\leq \|1-b^{-1/2}\| \leq \varepsilon/4.$ 

We now split each $e_k$ into the sum of four projections $e_k^j\in Ae_k$, of equal trace $\tau(e^j_k)=\tau(e_k)/4$, 
$e_k = \Sigma_{j=1}^4 e_k^j$, and denote $y_{kj}=e^j_k y e^j_k$. We still have 
$\|e^j_k y_0 e^j_k - y_{kj}\|\leq \varepsilon/4$, with $(1/3 - \varepsilon/12)e^j_k \leq y_{kj} \leq 1/2 e^j_k$ and $E_A(y_{kj})=t_k e^j_k$.  
Let $p_k^j\in A(1-e^j_k)$ be a projection of trace $\tau(e^j_k-y^j_k)/t_k$. To see that this is possible, we need to have 
$\tau(e^j_k-y^j_k)/t_k \leq \tau(1-e^j_k)$, which is easily seen to hold true due to our choices. Note also that $\tau(p^j_k)  \geq \tau(e^j_k)$. 
Take $B_0$ to be a separable diffuse abelian 
von Neumann subalgebra of $p^j_kMp^j_k$ which is perpendicular to $Ap^j_k$. Let $f^j_k$ be a projection in $B_0$ such that $\tau(f^j_k)=\tau(e^j_k)$. 
Let  $v\in M$ be a partial isometry such that $v^*v=e^j_k$ and $vv^*=f^j_k$, which due to the assumption that any two separable diffuse 
abelian von Neumann subalgebras of $f^j_kMf^j_k$ are unitary conjugate, can be chosen so that in addition we have $v(e^j_k - y_{kj})v^*\in B_0$. 

Denote $g_{kj}= y_{kj} + v(e^j_k - y_{kj})v^* + v(y_{kj}(e^j_k - y_{kj}))^{1/2} + (y_{kj}(e^j_k - y_{kj}))^{1/2}v^*$. It is easy to check that $g_{kj}^2=g_{kj}$,  
i.e., $g_{kj}\in \Cal P(M)$. Moreover  

$$
E_A(g_{kj})= E_A(p^j_k g_{kj} p^j_k + e^j_k g_{kj} e^j_k) = E_A(v(e^j_k - y_{kj})v^*) + E_A(y_{kj}) 
$$
$$
=\tau(v(e^j_k - y_{kj})v^*)/\tau(p^j_k)p^j_k  + t_ke^j_k=t_k(p^j_k + e^j_k).
$$ 

Thus, by our assumption, each one of the projections $g_{kj} \in (p^j_k+e^j_k)M(p^j_k + e^j_k)$ can be paved over $A(p^j_k + e^j_k)$. 
Thus $y_{kj}=e^j_kg_{kj}e^j_k$ can be paved, so  $\Sigma_{k,j} e^j_k y e^j_k$ can be paved, implying that $y$ can be paved. 
By $(2)$, it follows that $y_0$ can be paved, and by $(1)$, $x$ can be paved as well. Thus, any selfadjoint element can be paved, so by 
$(3)$ any element in $M$ can be paved. 

Moreover, if we keep track of the number of projections necessary in the above pavings we see that in order to 
$\varepsilon/2$-pave a selfadjoint element,  $\text{\rm n}(\Cal P_0, \varepsilon/50) (\varepsilon^{-1}+1)$  many projections are sufficient. 
By $(3)$, we get that $\text{\rm n}(\varepsilon) \leq \text{\rm n}(\Cal P_0, \varepsilon/50)^2 (\varepsilon^{-1}+1)^2$. 

\hfill 
$\square$

\vskip .05in 
\noindent
{\bf 2.4. Remarks}. $1^\circ$ Let $A\subset M$ be a MASA in a II$_1$ factor and $x\in M\ominus A$, $\|x\|\leq 1$. 
If we view $x$ as an element in $M^\omega$, 
then its $\varepsilon$-paving number over $A^\omega$ is $\leq n$  iff for any
$\delta > 0$, there exists a partition of 1 with projections $p_1,
..., p_n\in \Cal P(A)$ with the property that the spectral projection
of $|\Sigma_i p_i x p_i|$ corresponding to the interval
$(\varepsilon, \infty)$ has trace $\leq \delta$. 

\vskip .05in 
$2^\circ$ We already mentioned in 2.1 $(b)$ several properties of the algebra $\text{\bf R}$: 
given any representation of $R$ as crossed product $D\rtimes \Gamma$, for some free action of a (necessarily countable amenable) group $\Gamma$ on $D$, 
$\Gamma$ acts freely on $D^\omega$ as well and we have $\text{\bf R}=D^\omega \rtimes \Gamma$; thus, $\text{\bf R}$ is amenable and 
has $D^\omega$ as a Cartan subalgebra, but it has large, non-separable center. In addition, note that due to Rohlin's theorem and F\o lner's 
condition for amenable groups, any two free actions of the same amenable group 
$\Gamma \curvearrowright D^\omega$ are conjugate by a unitary element in $\Cal N_{R^\omega}(D^\omega)$. 
Moreover, the 1-cohomology for $D^\omega \subset \text{\bf R}$ vanishes, so if $\Gamma, \Lambda \subset \Cal N_R(D)$ are two countable amenable groups 
of unitaries that implement free actions on $D$ and $\Delta: \Gamma \simeq \Lambda$ is a group isomorphism, then there 
exists $u\in \Cal N_{R^\omega}(D^\omega)$ such that $u u_g u^* =\Delta(u_g)$, $\forall u_g \in \Gamma$. In particular, if $u_1, u_2\in \Cal N_R(D)$ 
act freely on $D$, then there exists $u\in \Cal N_{R^\omega}(D^\omega)$ such that $uu_1u^*=u_2$. Note that in fact  all these properties hold also for 
countable amenable subgroups $\Gamma \subset \Cal N_{\text{\bf M}}(\text{\bf D})$ acting freely on $\text{\bf D}$.

\vskip .05in 
$3^\circ$ Recall that in ($4.1. (iii)$ and $4.3.3^\circ$ of [P7]) it was shown that if $B$ is a separable amenable von Neumann 
subalgebra in a II$_1$ von Neumann algebra $M$ such that the Pimsner-Popa index  [PiP] of 
the inclusions $pBp\subset pMp$ 
is infinite for any projection $p\in B$, $p\neq 0$,  then there exist non-normal conditional expectations of $M$ onto $B$, and thus  
$E_B$ is not the unique conditional expectation of $M$ onto $B$. 
In particular, if $A$ is a separable MASA in a II$_1$ von Neumann algebra $M$, then $E_A$ is not the unique conditional 
expectation of $M$ onto $A$, and thus $A\subset M$ cannot have the KS property, nor the paving property.   
We recall the argument in [P7], emphasizing a simplification that occurs in the case of a MASA. 

First one constructs a singular state $\varphi$ on $M$ 
with $\varphi_{|B} = \tau$ as follows:  Like in [P7], the hypothesis implies there exists $b_n\in L^1M_+$ satisfying  
$E_B(b_n)=1$ and $\tau(s(b_n))\leq 2^{-n}$, where for a positive element $b\geq 0$ in $M$, $s(b)$ denotes the support 
projection of $b$. (In the case $B=A$ is a MASA in $M$, the argument becomes much simpler, as one can take $b_n = 2^{n}q_n\in M$, 
with $q_n$ the following projection: let $e_{kk}\in A$ be a partition of $1$ with $2^n$ mutually equivalent projections 
and complement it to a set  $e_{jk}\in M$ of matrix units, then define $q_n = 2^{-n} \Sigma_{j,k} e_{jk}$, which is  a projection in $M$ with 
$E_A(q_n)=\Sigma_k e_{kk} E_A(q_n) e_{kk}=E_A(\Sigma_k e_{kk}q_n e_{kk})=2^{-n}$.) 
Then $\varphi_n=\tau(\cdot\, b_n)$ defines a normal state on $M$ which, since $E_B(b_n)=1$, satisfies  $\tau(yb_n)=\tau(y), 
\forall y\in B$. Take $\varphi$ to be a state on $M$ obtained as a weak-limit of $\varphi_n$.  
Then we still have $\varphi_{|B}=\tau_{|B}$ while $\varphi$ is singular on $M$ (because for each fixed $n$ we have $\varphi(1-\vee_{m\geq n} s(b_n))=0$  
and $\lim_n (1-\vee_{m\geq n} s(b_n))=1$). 

Next, since $B$ is amenable, by Connes' Theorem we can 
find a countable amena-\newline ble subgroup $\Cal U_0\subset \Cal U(B)$ such that $\Cal U_0''=B$. For each $x\in M$, put $\psi(x)=\int \varphi(u x u^*) \text{\rm d} u$, 
where the integral is in the Banach limit sense, over an invariant mean on the countable amenable group $\Cal U_0$. Then $\psi$ defines a state on $M$ 
which is in the $\sigma(M^*, M)$-closure of a countable set of singular states on $M$. By [Ak], it follows that $\psi$ is singular as well. Also, by its definition, 
$\psi$ has the span of $\Cal U_0$ in its centralizer and $\psi_{|B}=\tau_{|B}$. If now $a\in B$ is arbitrary and $a_n \in \text{\rm sp} \Cal U_0$, $\|a_n\|\leq \|x\|$ 
are so that $\|x-b_n\|_2\rightarrow 0$, then by Cauchy-Schwartz inequality for $\psi$, for all $x \in M$ we have 
$$
|\psi (ax) - \psi(a_nx)|\leq \psi((a-a_n)(a-a_n))^{1/2} \psi(x^*x)^{1/2}
$$
$$
=\tau((a-a_n)(a-a_n)^*)^{1/2} \psi(x^*x)^{1/2}=\|a-a_n\|_2  \psi(x^*x)^{1/2} \rightarrow 0.
$$
Similarly, $|\psi (xa) - \psi(xa_n)| \rightarrow 0$. Since $\psi(a_nx)=\psi(xa_n)$, $\forall n$, this shows that $\psi(ax)=\psi(xa)$, i.e. $B$ is in the centralizer of $\psi$. 
Taking $E:M \rightarrow B$ to be the unique conditional expectation satisfying $\psi(E(x)a)=\psi(xa)$, $\forall x\in M$, $a\in B$, we have constructed this way 
a singular (thus non-normal) conditional expectation of $M$ onto $B$.

\heading 3. Paving in the $L^2$-norm     \endheading

Given a MASA $A$ in a finite von Neumann algebra $M$ 
and $x$ an element in $M$ with $E_A(x)=0$, our strategy for estimating the norm of its paving $y=\Sigma_i p_i x p_i$, 
for $p_i\in \Cal P(A)$ a partition of $1$, will be to calculate 
the moments $\tau((y^*y)^n)$ and use the well known formula $\|y\|^2=\lim_n \tau((y^*y)^n)^{1/n}$. More generally, in order to have $\|y\| \leq c$, 
we need to prove that $\tau((y^*y)^n) \leq c^{2n}$, for large enough $n$. One way of controlling these moments is to construct the partitions $\{p_i\}_i \subset A$ 
so that to have ``high order of independence'' with respect to $x$. 

We will use the following terminology in this respect: Two sets $V, W\subset M\ominus \Bbb C$ are 
{\it $n$-independent} if any alternating word $x_1y_1.... x_ky_k$, with $k\leq n$ and $x_1 \in V\cup \{1\}$, 
$x_2, .., x_k\in V$, $y_1, ..., y_{k-1}\in W$, $y_k \in W \cup \{1\}$, has trace $0$ (unless $k=1$ and $x_1=y_1=1$).  An algebra $B_0\subset M$ is 
$n$-independent to $V$ if $V$ and $B_0 \ominus \Bbb C$ are $n$-independent. Note that $1$-independence amounts to what one usually 
calls  $\tau$-independence.  

More generally, if $P\subset M$ is a von Neumann subalgebra, then two sets $V \subset M\ominus P$, $W\subset M \ominus \Bbb C1$ are $n$-{\it independent 
relative to} $P$ if $E_P( \Pi_{i=1}^k x_iy_i)=0$, for all $1\leq k \leq n$, all $x_1\in V \cup \{1\}$, $x_i\in V$, $y_k\in W\cup\{1\}$, $y_i\in W$.

In this section we prove a general fact about independence, showing that given any MASA $A$ in a finite von Neumann algebra $M$, we can find partitions of $1$ 
in $A$ that are ``asymptotically $2$-independent'' with respect to any given countable set of elements in $M\ominus A$. In other words,  
given any countable set $X\subset M\ominus A$, there exists a diffuse abelian subalgebra $B_0\subset A^\omega$ such that any word 
with at most 4 alternating letters in $X$ and $B_0\ominus \Bbb C$ has trace $0$. In particular, $\|p_i x p_j\|_2 = \|p_i\|_2 \|p_j\|_2 \|x\|_2$, 
for any $p_i, p_j\in \Cal P(A)$, $x\in X$, so any partition $p_i\in B_0$ with small mesh  gives  $L^2$-pavings of $x\in X$, 
simultaneously for all $x\in X$: $\|\Sigma_i p_i x p_i\|_2^2 = 
\|x\|_2^2 \Sigma_i \tau(p_i)^2 \leq \max_i \{\tau(p_i)\}_i \|x\|_2^2$. We construct such partitions recursively, but another method, 
where moments are controlled through incremental patching, can be used instead (see Section 5.3).

\proclaim{3.1. Lemma}  $1^\circ$ If $\xi \in L^2M$,  $u\in \Cal U(M)$ are so that $u^2=1$, $ | \tau(\xi^*u\xi u^*)| \leq  c\|\xi\|_2^2$, for some $1\geq c>0$, 
then $p_1=(1+u)/2$, $p_2=(1-u)/2$  is a partition of $1$ with projections satisfying $\|p_1\xi p_1 + p_2 \xi p_2\|_2  \leq (1 + c)/\sqrt{2}\|\xi\|_2$. 
If in addition $|\tau(u\xi^*\xi)| \leq c \|\xi\|_2^2$and $|\tau(u)|\leq c/2$,  then $\|p_i \xi p_i\|_2 \leq (1+2c)/\sqrt{2} \|\xi\|_2 \|p_i\|_2$, $i=1, 2$. 
\vskip .05in 
$2^\circ$ If $\xi \in L^2M$  and $u\in \Cal U(M)$ satisfy $ \tau(\xi^*u\xi u^*) \leq  c\|\xi\|_2^2$, for some $c\leq 2^{-7}$, 
and $n \geq 2^7$, 
then the spectral projections $\{e_k\}_{1\leq k \leq n}$ of $u$ defined by $e_k=e_{[e^{2\pi i (k-1)/n}, e^{2\pi i k/n})}(u)$, 
give a partition of $1$ and satisfy $\|\Sigma_k  e_k \xi e_k\|_2 
\leq 3/4 \|\xi\|_2$. 
\endproclaim
\noindent{\it Proof}. $1^\circ$ We have $\| \xi + u\xi u^*\|_2^2= 2 \|\xi\|_2^2 + 2 \Re \tau(\xi^* u \xi u^*) \leq (2 +2c)\|\xi\|_2^2.$ 
Since $p_1 \xi p_1 + p_2 \xi p_2 = 2^{-1} (\xi + u\xi u^*)$, we get $\|p_1\xi p_1 + p_2 \xi p_2\|^2_2 \leq (1+c)/2 \|\xi\|_2^2$ 
and the first part of the statement follows. If $|\tau(u\xi^*\xi)| \leq c \|\xi\|_2^2$ and $|\tau(u)|\leq c/2$ as well, then $|\tau(p_1)-1/2|
=|\tau(u)|/2 \leq c/4$. Thus $\tau(p_1) \geq 1/2 - c/4 \geq 1/4$ and also $1/2 \leq \tau(p_1) + c/4\leq \tau(p_1)(1+c)/4$, so we get: 
$$
\|p_1 \xi p_1\|_2^2 = \tau((1+u)\xi^*(1+u)\xi)/4 =\tau(\xi^*\xi)/4 + \tau(u\xi^*\xi)/2 + \tau(u\xi^*u\xi)/4
$$
$$
\leq (1/4 + 3c/4) \|\xi\|_2^2\leq (1 + 3c)/2\|\xi\|_2^2 (1+c) \tau(p_1)  \leq (1 + 3c)^2/2\|\xi\|_2^2 \|p_1\|_2^2. 
$$
Similarly, by using that $p_2=(1-u)/2$, we obtain 
$$
\|p_2 \xi p_2\|_2^2 \leq (1 + 3c)/2\|\xi\|_2^2 (1+c)\tau(p_2) \leq   (1 + 2c)^2/2\|\xi\|_2^2 \|p_2\|_2^2.
$$

$2^\circ$ We may clearly assume $\|\xi\|_2=1$. Note that if we denote $\lambda_k = e^{2\pi i k/n}$, 
then $\|u-\Sigma_k \lambda_k e_k\| \leq |2\pi i /n - 1| < 2\pi/n$. 
Since the elements $\{e_j \xi e_k\}_{1\leq j, k \leq n}$ are mutually orthogonal in the Hilbert space $L^2M$, by using first Pythagora's 
theorem and then the inequality  $|\lambda_j^*\lambda_k-1|\leq 2$, 
$\forall  j, k$, we get: 

$$
4-4\|\Sigma_k e_k \xi e_k\|_2^2
= 4\|\xi\|_2^2-4\|\Sigma_k e_k \xi e_k\|_2^2 
= 4\|\Sigma_{j\neq k} e_j \xi e_k\|_2^2  
$$
$$
\geq \|\Sigma_{j\neq k} (\lambda_j^*\lambda_k -1) e_j \xi e_k\|_2^2
=\|(\Sigma_j \lambda_j^* e_j)\xi(\Sigma_k \lambda_k e_k)-\xi\|_2^2 
$$
$$
\geq \|u^*\xi u - \xi\|_2^2 -4 \|u-\Sigma_k \lambda_k e_k\| 
$$
$$
= 2 - 2\Re \tau(\xi^*u^*\xi u) - 4\|u-\Sigma_k \lambda_k e_k\| \geq 2 - 2c - 8\pi/n. 
$$

If we now choose $c<2^{-7}$ and $n\geq 2^7$, then from the first and last term of the above estimates we get 

$$
\|\Sigma_k e_k \xi e_k\|_2^2 \leq 1/2 +c/2 + 2\pi/n  \leq 9/16. 
$$

\hfill 
$\square$

\vskip .05in 
\noindent
{\bf 3.2. Remark}. For the following lemmas, 
it will be useful to recall that a unitary representation $\pi$ of a group $G$ on a Hilbert space $\Cal H$ is {\it weak mixing}  
if any of the following equivalent conditions is satisfied: 

\vskip .05in 
\noindent
$(3.2.1)$ Given any finite subset $F\subset \Cal H$ and any $\varepsilon > 0$, there exists $g\in G$ such that $|\langle \pi(g)(\xi), \eta \rangle| \leq 
\varepsilon$, $\forall \xi, \eta \in F$;  

\vskip .05in 
\noindent
$(3.2.2)$ $\pi$ has no non-zero finite dimensional 
invariant subspaces $\Cal H_0\subset \Cal H$; 
\vskip .05in 
\noindent
$(3.2.3)$ The representation $\pi \otimes \overline{\pi}$ of $G$ on $\Cal H \otimes \overline{\Cal H}$ is ergodic, i.e. it has 
no fixed non-zero vectors. 
\vskip .05in 
\noindent
$(3.2.4)$ For any finite dimensional subspace $\Cal H_0 \subset \Cal H$ and 
any $\varepsilon > 0$, there exists $g\in G$ such that, if $p_0$ denotes the orthogonal projection of $\Cal H$ onto 
$\Cal H_0$ and we still denote by $\pi$ the representation of $G$ on the space of Hilbert-Schmidt operators $HS(\Cal H)\simeq 
\Cal H \otimes \overline{\Cal H}$, then $Tr(\pi(g)(p_0) p_0)\leq \varepsilon Tr(p_0)$. 

\proclaim{3.3. Lemma} If $B \subset M$ is a diffuse von Neumann subalgebra, then  
the action $\text{\rm Ad}$ of its unitary group $ \Cal U(B)$  on $L^2(M\ominus (B' \cap M))$  
is weak mixing. Moreover, if $B$ is abelian, then the restriction of this $\text{\rm Ad}$-action to the subgroup of period $2$ elements, 
$\Cal U^0(B) \overset \text{\rm def} \to = \{u\in \Cal U(B)\mid u^2=1\}$, is still weak mixing. 
\endproclaim
\noindent
{\it Proof}. If the action is not weak mixing, then there exists a non-zero finite dimensional subspace $\Cal H_0\subset L^2(M\ominus (B' \cap M))$ 
satisfying $u\Cal H_0u^*=\Cal H_0$, $\forall u\in \Cal U(B)$. In particular, if $A\subset B$ is a diffuse abelian subalgebra of $B$ 
and $\Cal U^0=\Cal U^0(A)$ denotes the group of unitaries of period two in $A$ as in part $2^\circ$, 
then $\Cal H_0$ is invariant to the representation $\xi \mapsto \text{\rm Ad}(u)(\xi)=u\xi u^*$ of $\Cal U^0$ 
on $\Cal H_0$. Since the image of this representation is an abelian subgroup $\Cal V^0$ of $\Cal U(\Cal H_0)$, it can be diagonalized. Thus,  
$\Cal H_0=\Sigma_j \Bbb C \xi_j$, with $\xi_1, \xi_2, ...., \xi_n$ an orthonormal basis of $\Cal H_0$ such that Ad$(\Cal U^0)(\xi_j)\subset 
\Bbb T \xi_j$, $\forall j$, and since any element in $\Cal V^0$ has period $2$, we actually have Ad$(\Cal U^0)(\xi_j)\subset 
\{\pm \xi_j\}$, $\forall j$. But the group $(\Cal U^0, \| \ \|_2)$ is Polish and contractible.  
This can be seen by taking a $\|\cdot\|_2$-continuous path $\{p_t \mid 0\leq t \leq 1\}\subset  \Cal P(A)$ with $\tau(p_t)=t$, $p_t \leq p_{t'}$ iff $t\leq t'$, 
then defining the continuous path of group morphisms 
$\Phi_t: \Cal U^0(A) \rightarrow \Cal U^0(A)$, by $\Phi_t(u)=p_t + u(1-p_t)$, which satisfies $\Phi_0(u)=u$, $\Phi_1(u)=1$, $\forall u\in \Cal U^0(A)$. 
Since  $\Cal U^0$ is contractible and the representation Ad is continuous, and since  the one dimensional 
representation $u\mapsto \text{\rm Ad}(u)$ lies in $\{\pm 1\}$, this representation must be trivial, i.e.  
$u\xi_ju^*=\xi_j$, $\forall u\in \Cal U^0$, $\forall j$. Hence,  $u\xi =\xi u$ for all $\xi\in \Cal H_0$ and for all $u\in \Cal U(A)$ 
(because $\Cal U^0$ generates $A$ as a von Neumann algebra). Since any $u\in \Cal U(B)$ 
lies in a diffuse abelian von Neumann subalgebra $A\subset B$, it follows that $u\xi =\xi u$,  for all $u\in \Cal U(B)$ and all $\xi\in \Cal H_0$.
This means that $\Cal H_0\subset L^2(B'\cap M)$, while at the same time $\Cal H_0\perp B'\cap M$, implying that $\Cal H_0=0$.

\hfill $\square$

We have seen in the previous lemma that if $A$ is a diffuse abelian von Neumann subalgebra in $M$, then the Ad-action 
of the group $\Cal U^0(A)$ (of period $2$ unitaries in $A$) on $L^2(M\ominus A'\cap M)$ is weak mixing. 
We will show in the next lemma that one can choose the corresponding mixing elements in $\Cal U^0(A)$ so that to be 
approximately $\tau$-independent  with respect to any given finite set in $M$ and to have approximately $0$-trace. Proving this in the abelian case 
is quite straightforward. But due to its possible independent interest, 
we will actually prove this type of result for arbitrary diffuse von Neumann subalgebras $B\subset M$, 
a fact that will make the argument a bit more lengthy.

\proclaim{3.4. Lemma} Let $M$ be a finite von Neumann algebra and $B\subset M$ a diffuse von Neumann 
subalgebra. Given any finite dimensional subspaces $X\subset L^2(M\ominus (B \vee (B'\cap M)))$, $Y\subset L^1M$, and any $\delta> 0$, 
there exists a period $2$ unitary element $u\in B$ such that $|\tau (u\xi^*_1 u^*\xi_2)|\leq \delta \|\xi_1\|_2\|\xi_2\|_2$, $|\tau(u \eta)|\leq 
\delta \|\eta\|_1$, $\forall \xi \in X$, $\eta \in Y$. 
\endproclaim
\noindent
{\it Proof}. We first prove that given any finite dimensional subspace $\Cal H_0 \subset L^2(M\ominus B\vee (B'\cap M))$, there 
exists a diffuse abelian von Neumann subalgebra $A\subset B$ such that $E_{A'\cap M}(\Cal H_0)=0$. 
Since $\Cal H_0$ is perpendicular to $B\vee (B'\cap M))$,  it is also perpendicular to 
$B'\cap M$, so by Lemma 3.3 there exists $u\in \Cal U(B)$ such that  $|\tau(\xi^* u \xi u^*)|\leq c \|\xi\|_2^2$, $\forall \xi\in \Cal H_0$. 
By Lemma 3.1, 
if $c = 2^{-7}$ and $n_1=2^7$, then there exists a partition of $1$ with $n_1$ projections $\{e^1_j\}_j$ in $B$ such that $\|\Sigma_j e^1_j \xi e^1_j \|_2 
\leq 3/4 \|\xi\|_2$, $\xi \in \Cal H_0$. Since $\Cal H_0 \perp B\vee B'\cap M$, we also have $\Sigma_j e^1_j \Cal H_0 e^1_j \perp \Sigma_j e^1_j(B\vee B'\cap M)e^1_j 
= B_1 \vee (B_1'\cap M)$, where $B_1=\Sigma_j e^1_j B e^1_j$ (cf. Lemma 2.1 in [P1]). We can thus 
continue recursively, replacing $\Cal H_0$ by $\Sigma_j e^1_j \Cal H_0 e^1_j$ and $B$ by $B_1$, to get a partition of 1 with $n_2$ projections  $\{e^2_k\}_k 
\subset B_1,$ which refines $\{e^1_j\}_j$ 
and satisfies $\|\Sigma_j e^2_k \xi e^2_k \|_2 \leq 3/4 \|\Sigma_k e^1_j \xi e_j^1 \|_2$. Thus,  

$$
\|\Sigma_k e^2_k \xi e^2_k\|_2 \leq 3/4 \| \Sigma_j e^1_j \xi e^1_j \|_2\leq (3/4)^2 \|\Sigma_k e^2_k \xi e_k^2 \|_2, \forall \xi \in \Cal H_0.
$$
By iterating this procedure we get a sequence of finer and finer 
partitions of $1$, $\{e^m_j\}_{j=1}^{n_m}\subset B$, where $n_m=2^{7m}$,  such that 
$\|\Sigma_i e^m_j \xi e^m_j\|_2 \leq (3/4)^m \|\xi\|_2 \leq \varepsilon \|\xi\|_2$, $\forall \xi \in \Cal H_0$ (Note that  
the number $n=n_m$ of projections necessary to get $(3/4)^m  \leq \varepsilon$ 
is majorized by $2^{7\ln (1/\varepsilon)/\ln (4/3) + 1}$, thus $n \leq 2^7 (1/\varepsilon)^{7\ln 2/\ln (4/3)}\approx 2^7 (1/\varepsilon)^{17.02}$, 
so the order of magnitude of the size of the partition is $\varepsilon^{-17.02}$). 
 
If we  define $A$ to be the von Neumann algebra generated by $\{ e^m_j \mid 1\leq j \leq n_m, m\geq 1\}$, it follows that $E_{A'\cap M}(\xi)=0$, $\forall \xi \in \Cal H_0$.  
Consider now the group $\Cal U^0=\Cal U^0(A)$ and  note that it is Polish with respect to the topology implemented by $\| \ \|_2$. Also, 
$(\Cal U^0, \| \  \|_2)$ is connected, in fact even contractible (due to the same construction as in the above proof of Lemma 3.3). 
Consider  the Hilbert space $\Cal K = HS( L^2(\overline{\text{\rm sp}}(A F A))) \oplus HS(L^2M)$ and the unitary representation $\pi$ 
of $\Cal U^0$   on $\Cal K$ given by $u \mapsto \text{\rm Ad}(L_u R_u) \oplus \text{\rm Ad}(L_u)$, $\forall u\in \Cal U^0$, where 
$HS(\Cal H)$ denotes the space of Hilbert-Schmidt operators on the Hilbert space $\Cal H$ (i.e. $HS(\Cal H)=\{x\in \Cal B(\Cal H) \mid Tr(x^*x) < \infty\}$), 
and $L_u$ (resp. $R_u$) are the operators of left (resp. right) multiplication by $u\in \Cal U^0\subset A$. Thus, 
if we identify in the usual way $HS(\Cal H) \simeq \Cal H \overline{\otimes} \Cal H^*$, then for each $\xi_1, \xi_2 \in \overline{\text{\rm sp}}(A \Cal H_0 A)$, 
$\eta_1, \eta_2 \in L^2M$ we have 

$$
\pi(u) ( \xi_1 \otimes \xi_2^* \oplus \eta_1 \otimes \eta_2^*)
=(u\xi_1u^* \otimes u\xi_2^*u^* \oplus  u\eta_1 \otimes \eta_2^*u^*). 
$$

Let $p=p_1 \oplus p_2 \in \Cal K$, where $p_1$ is the orthogonal projection of $\overline{\text{\rm sp}}(A \Cal H_0 A)$ onto $\Cal H_0$  
and $p_2$ is the orthogonal projection of $L^2M$ onto this same space. Let $K_p=\overline{\text{\rm co}}^w\{ \pi(u)(p) \mid u \in \Cal U^0\}$. 
Note that $\pi(u)(K_p)=K_p$ and $\|\pi(u)(x)\|_{\Cal K}=\|b\|_{\Cal K}$, $\forall x \in \Cal K$. Note also that all elements $x$ in $K_p$ are of the form $x=x_1 \oplus x_2$, 
with $x_1\in HS(\overline{\text{\rm sp}}(A F A))\subset \Cal B(\overline{\text{\rm sp}}(A F A))$, $x_2 \in HS(L^2M)\subset \Cal B(L^2M)$ 
positive operators when viewed as acting on the corresponding Hilbert space. 

Since $K_p$ is convex, weakly closed and bounded in $\Cal K$, there exists   a unique 
element $b=b_1 \oplus b_2 \in K_p$ of minimal Hilbert norm $\| \ \|_{\Cal K}$. Since $\pi(u)(b) \in K_p$ and has the same norm as $b$, it 
follows that $\pi(u)(b)=b$, $\forall u \in \Cal U_0$. Thus, $\text{\rm Ad}(L_uR_u)(b_1)=b_1$, $\text{\rm Ad}(L_u)(b_2)=b_2$, $\forall u\in \Cal U^0$. 
This means that, as (positive) operators on the corresponding Hilbert space, $b_1$ commutes with $L_uR_u$ and $b_2$ with $L_u$, $\forall u\in \Cal U^0$. 
Thus, the spectral decomposition of $b_1$ (resp $b_2$) commutes with these unitaries. Since $b_1, b_2$ are Hilbert-Schimdt, they are 
in particular compact, so any spectral projection corresponding to $(t, \infty)$ for $t>0$ is finite dimensional. It follows that if $b_1\neq 0$ (resp. $b_2\neq 0$) 
then there exists a non-zero finite dimensional subspace $\Cal H_1 \subset \overline{\text{\rm sp}}(A F A)$  such that $u\Cal H_1u^*=\Cal H_1$ 
(resp. $\Cal H_2 \subset L^2M$ such that $u\Cal H_2 = \Cal H_2$), $\forall u \in \Cal H_2$. 

Let us first notice that this implies $\Cal H_2=0$. Indeed, because $\Cal U^0\Cal H_2=\Cal H_2$ implies $ A \Cal H_2=\Cal H_2$, which contradicts the 
finite dimensionality of $\Cal H_2$, unless $\Cal H_2=0$. To see that $\Cal H_1=0$ as well, 
note that $\Cal H_1 \ni \xi \mapsto \text{\rm Ad}(u)(\xi) \in \Cal H_1$, $\forall u\in \Cal U^0$ defines a continuous unitary representation of 
the abelian Polish group $\Cal U^0$ on $\Cal H_1$. 
Since the image of this representation is an abelian subgroup of $\Cal U(\Cal H_0)$, it can be diagonalized. Thus,  
$\Cal H_0=\Sigma_j \Bbb C \xi_j$, with $\xi_1, \xi_2, ...., \xi_n$ an orthonormal basis of $\Cal H_0$ such that Ad$(\Cal U^0)(\xi_j)\subset 
\Bbb T \xi_j$, $\forall j$. Since any element in $\Cal U^0$ has period $2$, we actually have Ad$(\Cal U^0)(\xi_j)\subset 
\{\pm \xi_j\}$, $\forall j$. But as we have noticed above, the Polish group $(\Cal U^0, \| \ \|_2)$ is connected (even contractible), 
implying that 
$u\xi_ju^*=\xi_j$, or equivalently $u\xi = \xi u$, $\forall u\in \Cal U^0$, $\forall j$. Hence,  $a\xi =\xi a$ for all $\xi \in \Cal H_1$ and all $a\in A$ 
(because $\Cal U^0$ generates $A$ as a von Neumann algebra), i.e. $\Cal H_1 \subset L^2(A'\cap M)$. But $E_{A'\cap M}(F)=0$ 
implies $E_{A'\cap M}(A F A)=0$ and thus $E_{A'\cap M}(\overline{\text{\rm sp}}(A F A))=0$, so in particular $E_{A'\cap M}(\Cal H_1)=0$. 
We have thus proved that $\Cal H_1 \subset L^2(A'\cap M)$ and $\Cal H_1 \perp L^2(A'\cap M)$, showing that $\Cal H_1=0$. 

This implies $b=0$ and thus $0\in K_p$. Hence, for any $\delta > 0$ there exists $u\in \Cal U^0$ such 
that $Tr(\pi(u)(p)p)) < \delta$. Indeed, for if there exists $\delta_0 >0$ such that $Tr(\pi(u)(p)p)) \geq \delta_0$, $\forall u\in \Cal U^0$, 
then $Tr(x p)\geq \delta_0$ for all $x\in K_p$, in particular for $x=0\in K_p$, giving $0 \geq \delta_0$, a contradiction.

But $Tr(\pi(u)(p)p) < \delta$ implies  that we have both $u \Cal H_0 \perp_\delta \Cal H_0$ and $u\Cal H_0 u^*\perp_\delta \Cal H_0$, in particular 
$|\tau (u\xi_1^* u^*\xi_2)| <  \delta \|\xi_1\|_2\|\xi_2\|_2$, $|\tau(u \eta)| <  
\delta \|\eta\|_1$, for all non-zero elements 
$\xi_1, \xi_2 \in X$, $\eta \in Y'{Y'}^* + \Bbb C 1$. Thus, if we embed $Y (\subset L^1M)$ in some $Y'{Y'}^*$ 
for some appropriate finite subspace $Y' \subset L^2M$, then all required conditions are satisfied.   

\hfill 
$\square$

\proclaim{3.5. Lemma} Let $M$ be a finite von Neumann algebra and $\Cal H_1 \subset L^2M$, $\Cal H_2\subset L^1M$ 
finite dimensional spaces. Given any $\delta>0$, there exists $\delta'>0$ such that if $x\in M$ satisfies $\|x\|\leq 1$ and 
$\|x\|_2 \leq \delta'$, then $\|x\xi\|_2 \leq \delta \|\xi\|_2$, $\|x\eta\|_1 \leq \delta \|\eta\|_1$, $\forall \xi \in \Cal H_1$, 
$\eta \in \Cal H_2$. 
\endproclaim
\noindent
{\it Proof}. The first part follows from the fact that norm $\| \ \|_2$ implements the $so$-topology on the unit ball of $M$ while the product with a compact 
operator (such as the orthogonal  projection of $L^2M$ onto the finite dimensional space $\Cal H_1$) turns $so$-convergence 
into operator norm convergence. To prove the second part, note that it is sufficient to show that 
given any $\delta>0$ and any finite set $\{\eta_i\}_i \subset \Cal H_2$ which is ``$\delta/2$-dense'' in the $L^1$-unit ball of $\Cal H_2$, 
there exists $\delta'>0$ such that if $x\in M$ satisfies $\|x\|\leq 1$, $\|x\|_2 \leq \delta'$, then $\|x\eta_i\|_1 \leq \delta/2$. In turn, this fact  is an 
immediate consequence of the first part, the Cauchy-Schwartz inequality and the fact  that 
any $\eta\in L^1M$, $\|\eta\|_1=1$ can be decomposed as a product $\xi_1\xi_2$ with $\xi_i\in L^2M$, $\|\xi_1\|_2=\|\xi_2\|_2=1$.  

\hfill 
$\square$  

\proclaim{3.6. Theorem} Let $n\geq 1$ be an integer. Given any  finite von Neumann algebra $M$, any diffuse von Neumann subalgebra  
$B\subset M$, any  finite sets $X \subset  L^2(M \ominus B\vee (B'\cap M))$, $Y \subset L^1M$   
and  any  $\alpha> 0$, there exists a finite dimensional von Neumann subalgebra $C \subset B$ generated by $2^n$ minimal 
projections of trace $2^{-n}$ such that
\vskip .05in 
\noindent
$(a)$ $|\tau(a_1 \xi_1 a_2 \xi_2)|\leq \alpha \|a_1\|_2 \|a_2\|_2$,  $\forall \xi_1, \xi_2 \in X$, $\forall a_1, a_2 \in C \ominus \Bbb C$. 
\vskip .05in 
\noindent
$(b)$ $|\tau(\eta a)| \leq \alpha \|a\| $, $\forall a\in C \ominus \Bbb C$, $\forall \eta \in Y\cup XX^*$.   

\vskip .05in
In particular, if $q_1, ..., q_{2^n}\in C$ are the minimal projections in $C$, then 

\vskip .05in 
\noindent
$(a')$ $| \|q_i \xi q_j \|_2^2 - \|\xi\|_2^2 \| \tau(q_i)\tau(q_j) | \leq 3\cdot 2^{-n}\alpha $, $\forall i,j$, $\forall \xi\in X$;

\vskip .05in 
\noindent
$(b')$ $|\tau(\eta q_i)-\tau(\eta)\tau(q_i))| \leq  \alpha$, $\forall i$, $\forall \eta \in Y\cup XX^*$.  
\vskip .05in 
\noindent
$(c')$  $\|q_i \xi q_i\|_2  \leq (2^{-n/2} \|\xi\|_2 + 2 \alpha^{1/2}) \|q_i\|_2$; 
$\| \Sigma_i q_i \xi q_i \|_2^2 \leq  2^{-n} \|\xi\|_2^2 + 3 \alpha$, $\forall i$, $\forall \xi\in X$. 
\vskip .05in 
\noindent
$(d')$  $\|q_i \xi q_i\|_1 \leq   (2^{-n/2} \|\xi\|_2 + 2 \alpha^{1/2})  \tau(q_i)$, $\forall i$, $\forall \xi \in X$.

\endproclaim
\noindent
{\it Proof}. Note that, without any loss of generality, we may assume $X=X^*$, $Y=Y^*$, $\|\xi\|_2=1$, $\|\eta\|_1=1$, $\forall \xi\in X$, 
$\forall \eta\in Y$. 

We prove the statement by induction over $n\geq 0$. If $n=0$ then  
$C=\Bbb C 1$ and the conditions are trivially satisfied. 
Assume we have proved the statement up to some $n$. Thus, there exists an abelian  
$2^n$-dimensional $^*$-subalgebra $C^0\subset B$ generated by minimal projections $q^0_1, ..., q^0_{2^n}\in B$ 
of trace $2^{-n}$ such that for all $a, a_1, a_2 \in C^0\ominus \Bbb C$, $\xi_1, \xi_2 \in X$, $\eta \in Y \cup XX^* \cup\{1\} $ we have 

$$
|\tau(a_1 \xi_1 a_2 \xi_2)|\leq \alpha' \|a_1\|_2 \|a_2\|_2, 
|\tau(\eta a)| \leq \alpha' \|a\| , \tag 1 
$$
where $\alpha'=2^{-n-2}\alpha$. 

Denote $B_0=\Sigma_i q_i B q_i$ and let $X_0$, respectively $Y_0$, be the linear span of the finite set $\Sigma_{i,j} q^0_i X q^0_j$,  
respectively $Y\cup X_0X_0^* \cup \{q^0_i \mid 1\leq i \leq 2^n\}$. Note that 
the condition $X \perp B \vee B'\cap M$ implies $X_0\perp B_0 \vee B_0'\cap M$. Indeed, by Lemma 2.1 in [P1], 
we have  $B_0 \vee B_0'\cap M= \Sigma_i q^0_i (B \vee B'\cap M)q^0_i$ and so if $x\in X$, $y\in B\vee B'\cap M$, then 

$$
\tau(q^0_j x q^0_k \Sigma_i q^0_i y q^0_i)= \delta_{jk} \tau(q^0_j xq^0_j y q^0_j)=\delta_{jk} \tau(x q^0_j yq^0_j)=0,
$$
the latter equality being due to the fact that $q^0_jyq^0_j \in B\vee B'\cap M$.

Let $\delta=2^{-n-2}\alpha$. By Lemma 3.5, there exists $1\geq \delta'>0$ such that if $x\in M$ satisfies $\|x\|\leq 1$ and $\|x\|_2\leq \delta'$, 
then $\|x\xi\|_2 \leq 3^{-1} \delta \|\xi\|_2$, $\|x\eta\|_1 \leq 3^{-1}\delta \|\eta\|_1$, $\forall \xi \in X_0, \eta\in Y_0$. 
By Lemma 3.4, there exists $v\in \Cal U^0(B_0)$ such that 

$$
|\tau(v\xi_1^*v^*\xi_2)| \leq 3^{-1} \delta \|\xi_1 \|_2\|\xi_2\|_2, |\tau(\eta v)|\leq {\delta'}^2 \|\eta\|_1, \forall \xi_1, \xi_2 \in X_0, \eta \in Y_0. \tag 2 
$$
Since $v$ is a period $2$ unitary commuting with all $q^0_i$,  
and $|\tau(vq^0_i)|\leq 2^{-n}{\delta'}^2$ (by last part of $(2)$ and the fact that  $q^0_i\in Y_0$), using the fact that $B_0$ is diffuse we can split each  
projection $q^0_i$ into the sum of two projections $q^0_i=q_{2i-1}+q_{2i}$ of trace $2^{-n-1}$ such that 
the period $2$ unitary element $u=\Sigma_i (q_{2i-1} - q_{2i})$ satisfies $\|u-v\|_2\leq \delta'$. 
Thus, $u$ satisfies $\|(v-u)\xi\|_2 \leq 3^{-1} \delta \|\xi\|_2$, $\|(v-u)\eta\|_1 \leq 3^{-1}\delta \|\eta\|_1$, $\forall \xi \in X_0, \eta\in Y_0$. 
Combining with $(2)$ and using the Cauchy-Schwartz inequality, we get: 

$$
|\tau(u\xi_1^*u^*\xi_2)| \leq |\tau(v\xi_1^*u^*\xi_2)|+ \|(u-v)\xi_1^*\|_2\|u^*\xi_2\|_2 \tag 3
$$
$$
\leq |\tau(v\xi_1^* v^* \xi_2)|+ \|v\xi_1^*\|_2 \|(u-v)\xi_2\|_2 + \|(u-v)\xi_1^*\|_2\|u^*\xi_2\|_2 
$$
$$
\leq  |\tau(v\xi_1^* v^* \xi_2)| +  2 \delta \|\xi_1\|_2\|\xi_2\|_2/3 \leq \delta \|\xi_1\|_2 \|\xi_2\|_2.
$$
Moreover, since $\delta'\leq 1/3 $ we have ${\delta'}^2 \leq \delta/3'$  
and thus
$$
|\tau(u\eta)| \leq |\tau((u-v)\eta)| + |\tau(v\eta)|  \leq \delta \|\eta\|_1, \forall \eta \in Y_0. \tag 3'
$$

Denote $C$ the linear span of $\{q_j \mid 1\leq j \leq 2^{n+1}\}$ and note that $C=C^0+ uC^0$, $C^0 \perp uC^0$. 
Let $x_i=a_i + ub_i\in C$, with $a_i, b_i \in C^0$. Thus, $\|x_i\|^2_2= \|a_i\|_2^2 + \|b_i\|_2^2$, $i=1, 2$. 
Take $\xi_1, \xi_2\in X$. Then we have 
$$
|\tau(x_1 \xi_1 x_2 \xi_2)| \tag 4
$$
$$
=|\tau(a_1 \xi_1a_2\xi_2)|+ |\tau(b_1u\xi_1a_2\xi_2)| + |\tau(a_1\xi_1 ub_2\xi_2)| + |\tau(b_1u\xi_1 b_2u\xi_2)|
$$
$$
=|\tau(a_1 \xi_1a_2\xi_2)|+ |\tau(u(\xi_1a_2\xi_2b_1))| + |\tau(u(b_2\xi_2a_1\xi_1))| + |\tau(u(\xi_1 b_2)u(\xi_2b_1))|
$$

By $(1)$, for the first term on the last line in $(4)$, we have the estimate 
$$
|\tau(a_1 \xi_1a_2\xi_2)|\leq \alpha' \|a_1\|_2\|a_2\|_2\leq \alpha'\|x_1\|_2\|x_2\|_2. \tag 5
$$

Since $\xi_1a_2\xi_2b_1$ and $b_2\xi_2a_1\xi_1$ belong to $Y_0$, by $(3')$ it follows that for the second term 
on the last line of $(4)$ we have the estimate: 
$$
|\tau(u(\xi_1a_2\xi_2b_1))|\leq \alpha' \|\xi_1a_2\xi_2b_1\|_1 \leq \alpha'\|\xi_1a_2\|_2 \|\xi_2 b_1\|_2  \tag 6
$$
$$
\leq \alpha'\|a_2\| \|b_1\| 
\leq 2^n \alpha' \|b_1\|_2 \|a_2\|_2 \leq 2^n \alpha' \|x_1\|_2 \|x_2\|_2, 
$$
where for the last row we have used the fact that for $a\in C^0$ we have 
$\|a\|\leq 2^{n/2} \|a\|_2$. Similarly, for the third term of the last line in $(4)$ we have
$$
|\tau(u(b_2\xi_2a_1\xi_1))| \leq  2^n \alpha' \|x_1\|_2 \|x_2\|_2. \tag 7
$$

Finally, for the fourth term of the last line in $(4)$, by $(3)$ and the fact that $\xi_1b_1, \xi_2b_2 \in X_0$,  we get 
$$
|\tau(u(\xi_1 b_2)u(\xi_2b_1))| \leq \alpha' \|\xi_1b_2\|_2 \|\xi_2b_2\|_2 \tag 8 
$$
$$
\leq \alpha' \|b_1\| \|b_2\| \leq 2^n \alpha' \|b_1\|_2 \|b_2\|_2 \leq 2^n \alpha' \|x_1\|_2 \|x_2\|_2. 
$$

By combining $(4)-(8)$, we thus obtain for all $\xi_1, \xi_2\in X$ and $x_1, x_2 \in C$ 

$$
|\tau(x_1 \xi_1 x_2 \xi_2)| \leq 4 \cdot 2^n \alpha' \|x_1\|_2 \|x_2\|_2 \leq \alpha \|x_1\|_2 \|x_2\|_2 \tag 9
$$

Similarly, if $\eta \in Y$ and $x=a + bu \in C$, then $\eta b \in Y_0$, $\|\eta b\|_1 \leq 2^n \|b\| \leq 2^n\|x\|$, and by the second part of $(1)$ and $(3')$ we get
$$
|\tau(x \eta)|\leq |\tau(a\eta)| + |\tau(u(\eta b))|   \tag 10 
$$
$$
\leq \alpha' \|a\| + \delta \|\eta b\|_1 \leq \alpha'\|x\| + 2^n \delta \|x\| \leq \alpha \|x\|,
$$
showing that $C$ satisfies conditions  $(a)$ and $(b)$  of the statement. 

If we now assume $(a)$ and $(b)$ are satisfied and combine them with  the identity 

$$
\|q_j \xi q_i\|_2^2 = 
\tau(q_i \xi^*q_j \xi)=\tau((q_i - \tau(q_i))\xi^* q_j \xi)+ \tau(q_i) \tau(q_j \xi \xi^*)
$$
$$
= \tau((q_i - \tau(q_i))\xi^*(q_j-\tau(q_j))\xi)
+ \tau(q_j) \tau((q_i - \tau(q_i))\xi^*\xi) 
$$
$$
+ \tau(q_i) \tau((q_j - \tau(q_j))\xi \xi^*) + \tau(q_i)\tau(q_j)\tau(\xi^*\xi)
$$
then we get: 
$$
| \|q_j \xi q_i\|_2^2  - \tau(q_i)\tau(q_j) \|\xi\|_2^2 | \leq 2^{-n} \alpha +  2 \cdot 2^{-n} \alpha = 3 \cdot 2^{-n}\alpha. 
$$
This proves that $(a)$ and $(b)$ imply $(a')$, while $(b')$ is trivial from $(b)$ and $(c')$ from $(a')$. 
Finally, $(d')$ follows from the first part of $(c')$, via the Cauchy-Schwartz inequality: 
$$
\|q_i\xi q_i\|_1 = \sup \{|\tau(q_i\xi q_i x) | \mid x\in M, \|x\|\leq 1\} 
$$
$$
\leq \|q_i \xi q_i \|_2 \sup \{\|q_i x\|_2 \mid x\in M, \|x\|\leq 1\} = \|q_i \xi q_i\|_2 \|q_i\|_2.
$$

\hfill $\square$

\proclaim{3.7. Corollary} Let $M$ be a $\text{\rm II}_1$ von Neumann algebra and $A\subset M$ a MASA in $M$. Let $X\subset M\ominus A$,  
$Y\subset M$ be finite sets, $n\geq 1$ an integer and $\alpha > 0$. There 
exists a partition of $1$ with projections $q_1, ..., q_{2^n}\in A$ of trace $2^{-n}$ such that if $C$ denotes the algebra generated by $\{q_i\}_i$ 
then for all $x\in X$, $y\in Y$ and 
$i=1, 2, ..., 2^n$ we have: 
\vskip .05in 
\noindent
$(a)$ $|\tau(a_1 x_1 a_2 x_2)| \leq \alpha \|a_1\|_2 \|a_2\|_2$, $\forall x_1, x_2\in X$, $\forall a_1, a_2 \in C\ominus \Bbb C$
\vskip .05in 
\noindent
$(b)$ $|\tau(y q_i)-\tau(y)\tau(q_i))| \leq  \alpha$, $\forall i$, $\forall y \in Y\cup XX^*$.  

\vskip .05in 
Moreover, we have for all $x\in X$ and $1\leq i,j \leq 2^n$: 

\vskip .05in 
\noindent
$(c)$ $| \|q_i x q_j \|_2^2 - \|x\|_2^2 \tau(q_i)\tau(q_j) | \leq 3\cdot 2^{-n}\alpha $.  

\vskip .05in 
\noindent
$(d)$  $\|q_ixq_i\|_1 \leq (2^{-n/2}\|x\|_2+ 2\alpha^{1/2})\tau(q_i)$.

\endproclaim
\noindent
{\it Proof}. Since $A$ is a MASA  we have $A'\cap M=A$ and since $M$ is II$_1$, $A$ must be diffuse. Thus, Theorem 3.6 applies. 

\hfill 
$\square$

\proclaim{3.8. Lemma} Let $M$ be a finite von Neumann algebra and $A\subset M$ a MASA. Given any separable von Neumann subalgebra 
$Q_0\subset M$, there exists a separable von Neumann algebra $Q\subset M$ that contains $Q_0$, such that 
 $E_A(Q)=A\cap Q$ $($i.e. $Q$ and $A \subset  M$ make a commuting square, as in Sec. $1.2$ in $\text{\rm [P4]})$ and $A_0 =A\cap Q$ is maximal abelian in $Q$. 
\endproclaim
\noindent
{\it Proof}. First note that by Theorem 3.6, given any  countable set $X=\{x_n\}_n\subset M$ there exists a countably generated abelian von Neumann subalgebra 
$B_1\subset A$ such that $E_{B_1'\cap M}(x_n)\subset B_1$, $\forall n$. Indeed, this is obtained by taking $B_1$ to be generated by 
the set $\{ E_A(x_n)\}_n$ and partitions $\{p^n_{j,m}\}_j \subset A$, $n, m \geq 1$, satisfying $\|\Sigma_j p^n_{j,m}( x_n - E_A(x_n))p^n_{j,m}\|_2 
\leq 2^{-m}$ (which exist by Theorem 3.6). If we then take $X$ to be a 
$\| \ \|_2$-dense subset in the unit ball of $Q_0$, it follows that the von Neumann algebra $Q_1$ generated by $B_1$ and $Q_0$ 
satisfies:

\vskip .05in 
\noindent
$(1)$ $B_1\subset A$ is separable and satisfies $E_{B_1'\cap M}(Q_0)\subset B_1$; $Q_1$ is generated by $Q_0, B_1$ and is separable; 
\vskip .05in 

Using this first part, it follows that we can 
construct recursively an increasing sequence of inclusions of separable von Neumann algebra $B_n \subset Q_n$, $n\geq 1$, 
satisfying the properties: 
\vskip .05in  
\noindent
$(2)$ $B_n \subset A$, $E_{B_n'\cap M}(Q_{n-1})\subset B_n$ 
and $Q_n$ is the von Neumann algebra generated by $B_n$ and $Q_{n-1}$. 

\vskip .05in

If we now define $A_0=\overline{\cup_n B_n}^w$ and $Q=\overline{\cup_n Q_n}^w$, then all required conditions are clearly satisfied. 

\hfill 
$\square$

\proclaim{3.9. Theorem} Let $M_n$ be a sequence of finite 
factors with $\text{\rm dim} M_n \rightarrow \infty$ and for each $n$, let $A_n\subset M_n$ be a MASA. 
Denote by $\text{\bf A}=\Pi_\omega A_n \subset \Pi_\omega M_n$. 
Let  $Q \subset \Pi_\omega M_n$ be an arbitrary separable von Neumann subalgebra such that $E_{\text{\bf A}}(Q) = \text{\bf A}\cap Q$, i.e. $Q$ and 
$\text{\bf A}\subset \Pi_\omega M_n$ make a commuting square, and denote $B_1=\text{\bf A}\cap Q$. There exists a diffuse von Neumann 
subalgebra $B_0\subset \text{\bf A}$ such that $B_0$ is $2$-independent to $Q\ominus B_1$ and $\tau$-independent to $B_1$, more precisely: 
\vskip .05in  
\noindent
$(a)$ $\tau(x_1a_1x_2a_2)=0$, $\forall x_1\in Q$, $x_2 \in Q\ominus B_1$, $a_1 \in B_0 \ominus \Bbb C$, $a_2 \in B_0$. 
\vskip .05in  
\noindent
$(b)$ $\tau(xa)=\tau(x)\tau(a)$, $\forall x\in Q$, $a\in B_0$, i.e.  
$(B_0 \vee B_1, \tau)\simeq (B_0, \tau) \overline{\otimes} (B_1, \tau)$. 
\vskip .05in 
\noindent
$(c)$ $\|e x f\|_2^2=\tau(e)\tau(f) \|x\|_2^2$, $\forall x\in Q\ominus B_1$, $e, f\in \Cal P(B_0)$.

\endproclaim
\noindent
{\it Proof}. Let $y_0=1, y_1, ... \in  Q$, be $\| \ \|_2$-dense in the unit ball of $Q$ and denote $x_n=y_n - E_{\text{\bf A}}(y_n)$. By 
Corollary 3.7, for each $n\geq 1$, there exists a $2^n$-dimensional $^*$-subalgebra $C_n \subset \text{\bf A}$ generated by minimal 
projections of trace $2^{-n}$ such that if we denote by $u_n\in C_n$ a unitary element with the properties $u_n^{2^n}=1$, $\tau(u_n^p)=0$, 
$1\leq k \leq 2^n-1$, then the following inequalities hold true for all $1\leq k, l \leq n$, $1 \leq p, q\leq 2^n-1$: 
$$
|\tau(u_n^px_ku_n^qx_l)| < 1/n; |\tau(y_ku_n^p)| < 1/n. \tag 1
$$
This implies that if $u_n=(u_{n,m})_m$, $x_n=(x_{n,m})_m$, $y_n=(y_{n,m})_m$ are so that $u_{n,m}\in \Cal U(A_m)$, $\|y_{n,m}\|\leq \|y_n\|$ 
and $x_{n,m}=y_{n,m}-E_{A_m}(y_{n,m})$, then for each $n$ there exists a neighborhood $\Cal V_n$ of $\omega$ such that if $m\in \Cal V_n$ 
then for all $1\leq k, l \leq n$, $1 \leq p, q\leq 2^n-1$ we have 
$$
|\tau(u_{n,m}^px_{k,m}u_{n,m}^qx_{l,m})| < 1/n; |\tau(y_{k,m}u_{n,m}^p)| < 1/n. \tag 2 
$$
Define $u=(u_m)_m \in \text{\bf A}$ by letting $u_m:=u_{n,m}$, $\forall m\in \Cal V_n \setminus \Cal V_{n-1}$. 
Conditions $(2)$ then imply that $u$ is a Haar unitary element in $\text{\bf A}$ and that $B_0:=\{u^n\}_n'' \subset \text{\bf A} $  
satisfies the independence conditions $(a)$, $(b)$ (and thus $(c)$ as well).

\hfill 
$\square$

\heading 4. Asymptotic freeness and paving over singular MASAs 
\endheading

Recall from [D1]  that a MASA $A$ in a finite von Neumann algebra $M$ is {\it singular} if 
the only unitaries in $M$ that normalize $A$ are the unitaries in $A$, i.e. $\Cal N_M(A)=\Cal U(A)$. It is easy to see 
that the existence of such a MASA in a finite von Neumann algebra $M$ implies $M$ is necessarily of type II$_1$  
(unless $M=A$). For concrete examples of singular MASAs in II$_1$ factors, see [D1], [P2] and Section 5.1 below. Note that by [P3], any separable 
II$_1$ factor has singular MASAs. The prototype singular MASA in the hyperfinite II$_1$ factor $R$ is the abelian von 
Neumann algebra $L(\Bbb Z)$ generated by the canonical unitary implementing the Bernoulli action $\Bbb Z \curvearrowright X = [0,1]^\Bbb Z$, 
in the representation of $R$ given by the Murray-von Neumann group measure space construction [MvN1], $R=L^\infty(X)\rtimes \Bbb Z$.

The main result of this section 
shows that if $A\subset M$ is a singular MASA in a finite von Neumann algebra, then the associated ultrapower MASA inclusion $A^\omega \subset M^\omega$ 
satisfies the paving property, and thus the KS property as well. 
In fact, we prove that given any countable set $X=X^*\subset  M^\omega$ perpendicular to $A^\omega$, there exists   
a diffuse subalgebra $B_0$ of $A^\omega$ which is {\it free independent} to $X$, relative to $A^\omega$, i.e., any alternating word in 
$X, B_0\ominus \Bbb C$ has $0$-expectation onto $A^\omega$ (in other words,  
$X, B_0$ are $n$-independent relative to $A^\omega$, $\forall n$). In particular, due to calculations of norms in [V2], this  implies 
that any $x\in X$ which is a selfadjoint element with two point spectrum, has the property that 
any partition of small mesh with projections in $A_0$, provides a paving of $x$. As we saw in 2.3.2$^\circ$, 
this is sufficient to ensure that ANY $x\in M^\omega$ with $E_{A^\omega}(x)=0$ can be paved with finite partitions in $A^\omega$,
and thus $A^\omega \subset M^\omega$ satisfies the Kadison-Singer property.  

More precisely, we prove the following:  

\proclaim{4.1. Theorem} Let $\Cal S=\{A_n\subset M_n\}_n$ be a sequence of singular MASAs in finite von Neumann algebras    
and denote $\text{\bf M}= \text{\bf M}(\Cal S, \omega) = \Pi_\omega M_n$, $\text{\bf A}=\text{\bf A}(\Cal S, \omega)=   \Pi_\omega A_n$. Then we have: 

\vskip .05in
\noindent
$(a)$ If $X\subset \text{\bf M}\ominus {\text{\bf A}}$, $Y\subset \text{\bf A}$  are countable sets, then there exists a diffuse von Neumann subalgebra 
$B_0\subset \text{\bf A}$  such that $B_0$ is $\tau$-independent to $Y$ and free independent to $X$ relative to $\text{\bf A}$, i.e.,  
$E_{\text{\bf A}}(x_0 \overset k \to{\underset i = 1 \to \Pi} y_ix_i )=0$, for all $k \geq 1$ and all 
$x_0 \in X \cup \{1\}$, $x_i \in X$, 
$y_i \in B_0\ominus \Bbb C$, $1\leq i \leq k$. 

\vskip .05in 
\noindent
$(b)$ Let $B\subset \text{\bf M}$ be a countably generated von Neumann subalgebra such that $E_{\text{\bf A}}(B)$ $=\text{\bf A}\cap B$, i.e. $B$ and $\text{\bf A}$ 
make a commuting square.  
Then there exists a diffuse von Neumann subalgebra $B_0\subset \text{\bf A}$ such that if we denote by $B_1=\text{\bf A} \cap B$, 
then $B_0$ and $B_1$ are in tensor product situation and $B\vee B_0 = B *_{B_1} (B_0 \otimes B_1)$. In particular, if $B \perp \text{\bf A}$ then $B \vee B_0 = B * B_0$. 
\endproclaim

One should note that by Lemma 3.8, any separable subalgebra $Q \subset \text{\bf M}$ is contained in a larger von Neumann 
subalgebra $B\subset \text{\bf M}$ satisfying the commuting square condition in part $(b)$ of 4.1. 

The above theorem implies that given any countable set  $X\perp \text{\bf A}$, one can find partitions of arbitrarily small mesh in $\text{\bf A}$ that 
are free with respect to $X$. For special type of elements $x\in \text{\bf M}$ with $0$-expectation on $\text{\bf A}$, 
such as unitaries or selfadjoint elements with $2$-points spectrum, any such ``free paving'' 
diminishes the operator norm, due to Kesten-type phenomena [Ke] and Voiculescu's calculations of spectra 
for products of free-independent variables [V2]:

\proclaim{4.2. Corollary} Let  $A_n \subset M_n$, $\text{\bf A}$, $\text{\bf M}$ be as above. Then we have:

\vskip .05in
\noindent
$1^\circ$ If $u\in \text{\bf M}$ is a unitary  element such that $E_{\text{\bf A}}(u)=0$, then for  any 
$n \geq 1$, 
there exists a  partition of $1$ with $n$ projections $q_1, ..., q_n \in \text{\bf A}$ 
such that $\|\Sigma_{i=1}^n q_i u q_i\| \leq  (\sqrt{n-1}+1)/n$.
\vskip .05in
\noindent
$2^\circ$. If $e$ is a projection in $\text{\bf M}$ 
such that $E_{\text{\bf A}}(e)=\tau(e)1$ and $\tau(e)\leq 1/2$, then for any 
$n \geq \tau(e)^{-1}$  
there exists a  partition of $1$ with $n$ projections $q_1, ..., q_n \in \text{\bf A}$ 
such that $\|\Sigma_{i=1}^n q_i e q_i - \tau(e) 1\| \leq  2/\sqrt{n}$. Also, there exists $p\in \text{\bf A}$ of trace $\tau(p)=1/2$ such that 
$\|pep + (1-p)e(1-p)\|\leq (\tau(e)(1-\tau(e))^{1/2}  + 1/2.$
\endproclaim

As we saw in Proposition 2.3, the paving of projections that expect on scalars on ultrapowers of MASAs, is in fact sufficient to ensure paving of any element, 
so from 4.2 above we deduce:  

\proclaim{4.3. Corollary (Kadison-Singer for ultraproduct of singular MASAs)} Let \newline $A_n \subset M_n$, $\text{\bf A}$, $\text{\bf M}$ be as above. 
Then the inclusion $\text{\bf A} \subset \text{\bf M}$ satisfies the KS property. Thus, any pure state on $\text{\bf A}$ 
has a unique state extension to $\text{\bf M}$ and $E_{\text{\bf A}}$ is the unique conditional expectation of $\text{\bf M}$ onto $\text{\bf A}$. Moreover,  
$\text{\bf A} \subset \text{\bf M}$ has the uniform paving property, with paving size $\text{\rm n}(\varepsilon)$ 
majorized by a scalar multiple of $ \varepsilon^{-6}$. 

\endproclaim

The proof of Theorem 4.1 will  follow quite closely the type of arguments 
that we have developed in [P6]. We will thus use extensively the notations, terminology 
and technical lemmas from that paper, which we recall here in details, for the reader's convenience.

\noindent
\vskip .05in \noindent {\bf
4.4. Notation}. Let $\Cal M$ be a von Neumann algebra.  If $v \in \Cal M$ is a
partial isometry with $v^*v = vv^*$, $X \subset \Cal M$ is a subset and
$k \leq n$ are nonnegative integers then denote $X_v^{0,n} \overset
\text{\rm def} \to = X$ and $X_v^{k,n} \overset \text{\rm def} \to =
\{ x_0 \overset k \to{\underset i = 1 \to \Pi} v_i x_i \mid x_i
\in X, \ 1 \leq i \leq k - 1, \ x_0, \ x_k \in X \cup \{ 1 \}, 
v_i\in \{v^j \mid 1\leq |j|\leq n\}  \}$.

\proclaim{4.5.  Lemma}  Let $A$ be a singular MASA in the finite von Neumann algebra $M$.  
Let $\varepsilon > 0$,  $n\geq 1$ an integer, $f\in A$ a non-zero projection, and $F = F^*\subset M\ominus A$, $Y\subset A$ finite sets.  
There exists a partial isometry $v\in Af$ such that $\tau(vv^*) > \tau(f)/2$ and 
$\| E_A (x) \|_{_1} \leq \varepsilon$, $\forall x \in \overset
n \to{\underset k = 1 \to \cup} F_v^{k,n}$, and $|\tau(yv^k)|\leq  \varepsilon \tau(vv^*)$, $\forall y\in Y$, $1\leq |k|\leq n$. 
\endproclaim
\noindent 
{\it Proof}.  It is clearly sufficient to prove the statement in case
$\|x\|, \|y\| \leq 1$, $\forall x\in  F$, $y\in Y$. Let $\delta > 0$. Denote $\varepsilon_0 = \delta, \
\varepsilon_k = 2^k \varepsilon_{k-1}, \ k \geq 1$. Let ${\mycal
W} = \{ v \in Af \mid vv^*\in \Cal P(A), \|E_A (x) \|_1 \leq \varepsilon_k \tau(v^*v), |\tau(yv^k)|\leq \varepsilon \tau(vv^*), \forall 1 \leq k \leq n, 
x \in F_v^{k,n},  y\in Y \}$.  Endow ${\mycal W}$ with the order $\leq$ in
which $w_1 \leq w_2$ iff $w_1 = w_2 w_1^*w_1$. $({\mycal W} , \leq)$
is then clearly inductively ordered.  Let $v$ be a maximal element
in ${\mycal W}$. Assume $\tau(v^*v) \leq \tau(f)/2$ and denote $p = f - v^*v$.  If $w$ is a
partial isometry in $Ap$  and  $u = v + w$, 
then for $x = x_0 \overset k \to{\underset i = 1 \to \Pi} u_i x_i
\in F_u^{k,n}$ we have
$$
x = x_0 \Pi_{i=1}^k v_i x_i + \Sigma_\ell \Sigma_i z_{0,i} \Pi_{j=1}^\ell
w_{i_j} z_{j,i}, \tag 1
$$
where the second sum is taken over all $\ell = 1, 2, \dots , k$ and all $i
= (i_1, \dots , i_\ell)$, with $1 \leq
i_1 < \cdots < i_\ell \leq k$, and where $w_{i_j} = w^s$ whenever $v_{j} =
v^s$,  $z_{0,i} = x_0 v_1 x_1 \cdots x_{i_1-1} p$, $ z_{j,i} = p x_{i_j}
v_{i_j+1} \cdots v_{i_{j+1}} x_{i_{j+1}}p$, for $1 \leq j < \ell$, 
and $z_{\ell,i} = p x_{i_\ell} v_{i_\ell + 1}$ $ \cdots v_k x_k$ .

\vskip .05in

By applying part $(d)$ of Corollary 3.7 to the finite set $X$ of all elements of the form $pzp-E_{Ap}(pzp) \in pMp\ominus Ap$, 
where $z$ is of the form $z_{j,i}$, for some $i = (i_1,
\dots , i_\ell)$, $1 \leq j \leq \ell - 1$, $\ell \geq 2$, as well as to the set $Y$ of elements $|E_{Ap}(pzp)|$ for such $z$, 
it follows that $\forall \alpha > 0$,  $\exists q\in \Cal P(Ap)$ such that
$$
\| q z q - E_{Ap} (pzp) q \|_{_{1,pMp}} < \alpha \tau_{pMp} (q), \tag 2
$$ 
$$
\|E_{Ap}(pzp)q\|_{_{1,pMp}} \leq 
(1+\alpha) \|E_{Ap}(pzp)\|_{_{1,pMp}} \tau(q)
$$

Since for $y_1, y_2, y \in M$ with $\|y_1\| \leq 1, \|y_2\| \leq 1$ we have $\|E_A(y_1 y y_2)\|_1 \leq \|y_1 y y_2\|_1\leq \|y\|_1$, it follows that 
for any $l \geq 2$ 
we have: 

$$
\|E_A (z_{0,i} w_{i_1} z_{1, i} w_{i_2} z_{2,i} \dots w_{i_\ell}
z_{\ell, i}) \|_{_1} \tag 3 
$$
$$
\leq \| w_{i_1} z_{1,i} w_{i_2} \|_{_1} 
=  \| q z_{1,i} q \|_1 = \| q z_{1,i} q \|_{1,pMp} \tau(p) , 
$$
which by applying consecutively to $z=z_{1,i}$ the two inequalities in $(2)$,  is further majorised by 
$$
\leq  ( \|E_{Ap} (z_{1,i}) q \|_{1,pMp} + \alpha \tau_{pMp}(q))  \tau(p)   \tag 4 
$$
$$
\leq (1+\alpha) \| E_{Ap}(z_{1, i}) \|_{1,pMp} \tau_{pMp}(q)  \tau(p) + \alpha \tau_{pMp}(q) ) \tau(p)
$$
$$
=(1+\alpha) (\| E_{A}(z_{1, i}) \|_1 \tau(p)^{-1}) (\tau (q) \tau(p)^{-1})  \tau(p) +\alpha \tau(q) 
$$
$$
=(1+\alpha) \|E_A(p x_{i_1}
v_{i_1+1} \cdots v_{i_2} x_{i_2}p)\|_1 \tau(p)^{-1} \tau(q) + \alpha \tau(q). 
$$

But since $p$ lies in $A$, we have $\|E_A(pyp)\|_1 = \|pE_A(y)p\|_1\leq \|E_A(y)\|_1$ for any $y\in M$. Also, 
since $1\leq i_1 < i_2 \leq \ell$ and $i_2-i_1 \leq k -1$, the element $y=x_{i_1} v_{i_1+1} \cdots v_{i_2} x_{i_2}$ 
belongs to $F_v^{m,n}$ with $m=i_2-i_1-1 \leq k-2$. Altogether, it follows that the last term in $(4)$ is majorized by 

$$
(1+\alpha) \|E_A(x_{i_1}
v_{i_1+1} \cdots v_{i_2} x_{i_2})\|_1 \tau(p)^{-1} \tau(q) + \alpha \tau(q) \tag 5 
$$
$$
\leq (1+\alpha) \varepsilon_{k-2} \tau(vv^*) \tau(p)^{-1} \tau(q) + \alpha \tau(q)  
$$
$$
\leq (1+\alpha) \varepsilon_{k-2}  \tau(q) + \alpha \tau(q), 
$$
where the last inequality is due to the fact that 
$\tau(vv^*)\leq 1/2$ implies $\tau(vv^*)/\tau(p)$ $ \leq 1$. If we now take $\alpha\leq \varepsilon_0/4$, from 
the first term of $(3)$ and last term of $(5)$, we get that for all $i=(i_1, ..., i_\ell)$ with $\ell\geq 2$ we have 

$$
\|E_A (z_{0,i} w_{i_1} z_{1, i} w_{i_2} z_{2,i} \dots w_{i_\ell}
z_{\ell, i}) \|_{_1} \leq 2 \varepsilon_{k-2} \tau(q). \tag 6 
$$

Since $2 \varepsilon_{k-2} \leq  \varepsilon_{k-1}$ and since there are at
most $\overset k \to{\underset i = 2 \to \Sigma} \left( \matrix k \\ i
\endmatrix \right) = 2^k - k - 1$ elements in the sum in $(1)$ for which
$\ell \geq 2$, from $(6)$ we get
$$
 \sum_{\ell \geq 2}  \sum_i \left\| E_A
\left( z_{0,i} \prod_{j=1}^\ell w_{i_j} z_{j,i} \right) \right\|_1 \leq (2^k - k -1) \varepsilon_{k-1} \tau(q) \tag 7
$$

Finally, from the sum on the right hand side of (1) we will now estimate
the terms with $\ell = 1$.  These are terms which are obtained from $x_0
v_1 x_1 v_2 x_2 \dots v_k x_k$ by replacing exactly one $v_i$ by $w_i$, so
they are of the form $z=z_{0,i} w_{i} z_{1,i} $, where $i=1, 2, ..., k$, 
$z_{0,i}=x_0v_1x_1... v_{i-1}x_{i-1}p$, $z_{1,i}=px_i v_{i+1} ... v_k x_k$ and $w_{i}=w^s$ if $v_i=v^s$. Note that there are $k$ of them.  

One should notice at this point that in the above estimates we
only used the fact that $w^*w = ww^* = q$ and that $A$ is a MASA, not the actual form of $w$,  
nor the fact that $A$ is singular.  We
will make the appropriate choice for $w$ now, to get the necessary
estimates for these last terms (in the process, we will also deal with the condition $|\tau(yv^k)| \leq \varepsilon \tau(vv^*)$, 
$\forall y\in Y$, which is required for $w$ to belong to $\mycal W$).  The singularity assumption on $A$ will play a crucial role, due to the following: 

\proclaim{4.6. Lemma} Let $A\subset M$  be a singular MASA. 
Let $Y_1=Y^*_1 \subset M \ominus A$, $Y_2\subset M$ be finite sets, $q\in A$ a nonzero projection. Given any $\beta > 0$ and $n\geq 1$ there exists 
a unitary element $w\in Aq$ such that $\|E_A(y_1w^i y_2)\|_1 < \beta $,  
$\|E_A(y_2w^i y_1)\|_1 < \beta $, $|\tau(y_2w^j)|\leq \beta$, for all $ y_1\in Y_1, y_2 \in Y_2$,   $0 < |i| \leq n$.  
\endproclaim
\noindent 
{\it Proof}. We may clearly assume $\|y_i\|\leq 1$, $\forall y_i\in Y_i$, $i=1, 2$. 
Let $\langle M, e_A \rangle$ be the Jones basic construction von Neumann algebra of the inclusion $A \subset^{e_A} M$, 
endowed with its canonical (semifinite) trace $Tr_{\langle M, e_A \rangle}$. Consider the semifinite von Neumann algebra $\Cal M = \langle M, e_A \rangle^{\oplus 2n} \oplus \Cal B(L^2M)^{\oplus 2n}$ 
and denote by $Tr$ the trace on $\Cal M$ defined by $Tr(x_1, x_2, ..., x_{2n}, y_1, y_2, ..., y_{2n})= \Sigma_j Tr_{\langle M, e_A \rangle}(x_j) + \Sigma_i Tr_{\Cal B(L^2M)}(y_i)$.  
Let $K_0\subset \Cal M$ denote the convex hull of the set 

$$\{ (w^j (\Sigma_{y_1\in Y_1} y_1^* e_A y_1) w^{-j}, w^j(\Sigma_{y_2\in Y_2} y_2e_{\Bbb C}y_2^*)w^{-j})_{1\leq |j| \leq n} \mid w\in \Cal U(Aq)\} 
\subset \Cal M.$$ 

One should notice right away that each $j$'th entry $z_j \in q\langle M, e_A \rangle q$, of the first $2n$ coordinates 
of an element $z=(z_j)_j\in K_0$ satisfies $z_j=  q z_j q$ and $e_A z_j =0= z_j e_A $ 
(the latter because $e_A w^j y_1 e_A = w^j e_A y_1 e_A =w^j E_A(y_1)  e_A=0$, $\forall y_1\in Y_1$; 
similarly $e_A y_1 w^j e_A=0$). 

Note further that $K_0$ is bounded both in the operator norm on $\Cal M$ (by $\Sigma_{y_1} \|y_1\|^2  + \Sigma_{y_2} \|y_2\|^2\leq |Y_1| + |Y_2|$) 
and in the Hilbert-norm $\|$ $\|_{2,Tr}$ on $\Cal M$ (by $\Sigma_{y_1} \|E_A(y_1^*y_1)\|_2^2 + \Sigma_{y_2} \tau(y_2y_2^*) \leq |Y_1| + |Y_2|$). Thus, its weak closure $K=\overline{K_0}^w$  
is a weakly compact bounded subset 
in both $\Cal M$ and $L^2(\Cal M, Tr)$. In particular, $K$ contains a unique element $b\in K$ with $\|b\|_{2,Tr}=\min \{\|z\|_{2,Tr} \mid z \in K\}$. 

Note also that 
the group $\Cal U(Aq)$ acts on $K_0$ by 

$$
\sigma_w((x_j, x'_j)_{1\leq |j|\leq n})=(w^jx_jw^{-j}, w^jx'_jw^{-j})_j, \forall w\in \Cal U(Aq),
$$ 
and that this action preserves the Hilbert norm $\|$ $\|_{2,Tr}$. 
Thus, $\sigma$ extends to an action of $\Cal U(Aq)$ on $K$, still denoted by $\sigma$. Since $\|\sigma_w(b)\|_{2,Tr}=\|b\|_{2,Tr}$, by the uniqueness of $b$ 
as the element of minimal norm in $K$, it follows that $\sigma_w(b)=b$, $\forall w\in \Cal U(Aq)$. 

Hence, if $b=(b_j, b'_j)_j$ are the $4n$ components of $b$, 
then for each $j$ with $1\leq |j|\leq n$, we have $w^jb_j=b_jw^j$, $w^jb'_j=b'_jw^j$ for all $w\in \Cal U(Aq)$. Since any unitary element in $Aq$ can be expressed as a $j$'th power 
of a unitary in $Aq$, it follows that $ub_j=b_ju$, $ub'_j=b'_ju$, $\forall u\in \Cal U(Aq)$. But since any element in $Aq$ is a linear combination of unitary elements 
in $Aq$, this implies $b_j\in Aq'\cap q \langle M, e_A \rangle q = A'\cap q\langle M, e_A \rangle q$ and $b'_j \in Aq'\cap q\Cal B(L^2M)q$. 
But by (1.4 in [P8]), the supremum of finite projections in $A' \cap q\langle M, e_A \rangle q$ 
is equal to the supremum of the projections $qve_Av^*q$ with $v\in \Cal N(A)$. Since $A$ is singular in $M$, this implies $b_j = e_A b_j e_A$. But $e_Ab_j=0$, 
and so $b_j = 0$, $\forall j$.  On the other hand, since $b'_j\in \Cal B(L^2M)$  are Hilbert-Schmidt (thus compact) and commute with the diffuse algebra $Aq$, it follows that $b'_j=0$, 
$\forall j$, as well. 

We have thus proved that $0_{\Cal M}=(0, ..., 0) \in K$. This implies that for any $\beta>0$ there exists $w\in \Cal U(Aq)$ such that 

$$
Tr_{\langle M, e_A \rangle}(w^j (y_1e_A y_1^*)w^{-j}(\Sigma_{y_2} y_2e_Ay_2^*)) \tag 1
$$
$$
 + Tr_{\Cal B(L^2M)}(w^j(\Sigma_{y_2}y_2e_{\Bbb C}y_2^*)w^{-j} e_{\Bbb C}) 
< \beta^2, 
$$ 
for all $y_1\in Y_1$, where the sums are taken over $y_2\in Y_2$. Indeed, for if not then  

$$
Tr_{\langle M, e_A \rangle}(w^j (\Sigma_{y_1} y_1 e_A{y_1}^*)w^{-j}(\Sigma_{y_2} y_2 e_A y_2^*)) 
$$
$$ + Tr_{\Cal B(L^2M)}(w^j(\Sigma_{y_2}y_2e_{\Bbb C}y_2^*)w^{-j} e_{\Bbb C}) 
\geq \beta^2,
$$ 
for all $w\in \Cal U(Aq)$. By taking convex combinations over $w\in \Cal U(Aq)$ 
and then weak limits, this would imply $Tr ((b_j \Sigma_{y_2} y_2e_Ay_2^*, b'_j e_{\Bbb C})_j)\geq \beta^2 $, $\forall b=(b_j,b'_j)_j \in K$, in particular for $b=0$,   
thus $0\geq  \beta^2  > 0$,  a contradiction. 

In particular, any $w\in \Cal U(Aq)$ that satisfies $(1)$,  will also satisfy 

$$
Tr_{\langle M, e_A \rangle} (w^j (y_1e_A{y_1}^*)w^{-j} y_2 e_Ay_2^*) 
< \beta^2 \tag 2
$$
$$
Tr_{\Cal B(L^2M)}(w^jy_2e_{\Bbb C}y_2^*w^{-j}e_{\Bbb C}) < \beta^2,
$$ 
for all $y_1\in Y_1, y_2 \in Y_2$ and all $j$ with $1\leq |j|\leq n$. 
The second set of these inequalities, translates into $|\tau(w^jy_2)|^2 \leq \beta^2$, $\forall y_2\in Y_2$, $1\leq |j|\leq n$. 
At the same time, by taking into account the definitions of $\|$ $\|_1$ and of $Tr_{\langle M, e_A \rangle}$, and  
by using  the Cauchy-Schwartz inequality in $(\langle M, e_A \rangle, Tr)$, the first set of the inequalities entails the estimates  

$$
\|E_A(y_1 w^j y_2)\|_1=\sup \{ |\tau(y_1w^jy_2 a)| \mid a\in  A, \|a\|\leq 1\}  \tag 3
$$
$$
= \sup \{ | Tr(e_Ay_1w^jy_2 e_A a)| \mid a\in A, \|a\| \leq 1\} 
$$
$$
\leq  Tr(e_Ay_2^*w^{-j}y_1^* e_Ay_1w^{j}y_2 e_A)^{1/2} Tr(e_A)^{1/2}
$$
$$
= Tr(w^{-j}y_1^* e_A y_1 w^j y_2 e_Ay_2^*)^{1/2} \leq \beta, 
$$
and similarly $\|E_A(y_2 w^j y_1)\|_1\leq \beta$, $\forall y_1\in Y_1=Y_1^*$, $y_2 \in Y_2$, $j=\pm1, \pm 2, ..., \pm n$.

\hfill 
$\square$ 

\vskip .1in

\noindent 
{\it End of proof of} 4.5: Denote by $Z$ the set of elements of the form 
$x_0v_1x_1... v_{i-1}x_{i-1}p$, or $px_i v_{i+1} ... v_k x_k$, for all possible choices arising from elements in $\overset
n \to{\underset k = 1 \to \cup} F_v^{k,n}$. 
By applying Lemma 4.6 to $\beta  =\varepsilon_{k-1}\tau(q)/2k$, $n\geq 1$ and $Y_2=Y\cup Z \cup Z^* \cup \{E_A(z) \mid z\in Z\cup Z^*\}$, 
$Y_1=\{y_2-E_A(y_2) \mid 
y_2 \in Y_2\}$,  it follows that there exists $w\in \Cal U(Aq)$ such that 

$$ 
\|E_A( ((x_0 v_1 x_1 \dots v_{j-1} x_{j-1} - E_A(x_0 v_1 x_1 \dots v_{j-1} x_{j-1}p)w_j x_j v_{j+1} \dots v_k x_k) \|_{_1} \tag 8 
$$
$$
\leq \varepsilon_{k-1} \tau(q)/2k,  
$$
$$
\|E_A(x_0 v_1 x_1 \dots v_{j-1} x_{j-1} w_j (x_j v_{j+1} \dots v_k x_k - E_A(px_j v_{j+1} \dots v_k x_k))) \|_{_1}  \tag 8' 
$$
$$
\leq \varepsilon_{k-1} \tau(q)/2k.  
$$
$$
|\tau(w^j y_2)| \leq \varepsilon \tau(q). \tag 8''
$$

From $(8)$ and $(8')$, it follows that for each element with $\ell=1$ in the summation $\Sigma_\ell \Sigma_i z_{0,i} \Pi_{j=1}^\ell
w_{i_j} z_{j,i}$ in $(1)$, i.e., of the form $x_0 v_1 x_1 \dots v_{j-1} x_{j-1} w_j x_j v_{j+1} \dots v_k x_k$, 
we have the estimate: 

$$
\|E_A(x_0 v_1 x_1 \dots v_{j-1} x_{j-1} w_j x_j v_{j+1} \dots v_k x_k)\|_1 \tag 9 
$$
$$
\leq 2 \varepsilon_{k-1} \tau(q)/ 2k +  \|E_A(x_0 v_1 x_1 \dots v_{j-1} x_{j-1})w_j E_A(x_j v_{j+1} \dots v_k x_k)\|_1
$$
$$
\leq \varepsilon_{k-1} \tau(q)/k + \gamma,   
$$
where $\gamma$ is the minimum of $\|E_A(x_0 v_1 x_1 \dots v_{j-1} x_{j-1})q\|_1$, $\|qE_A(x_j v_{j+1} \dots v_k x_k)\|_1$, which by  
the second inequality in $(2)$ is majorized by the minimum between 
$\|(1+\alpha)E_A(px_0 v_1 x_1 \dots v_{j-1} x_{j-1}p)\|_1 \tau(q)$ and  $(1+\alpha)\|E_A(px_j v_{j+1} \dots v_k x_kp)\|_1 \tau(q)$. 
Since $\|E_A(pyp)\|_1 = \|pE_A(y)p\|_1 \leq \|E_A(y)\|_1$, the latter is majorized by 
the minimum between $(1+\alpha) \|E_A(x_0 v_1 x_1 \dots v_{j-1} x_{j-1})\|_1 \tau(q)$, $(1+\alpha)\|E_A(x_j v_{j+1} \dots v_k x_k)\|_1 \tau(q)$.  
Both elements $x_0v_1 x_1 \dots v_{j-1} x_{j-1}$,  $x_j v_{j+1} \dots v_k x_k$ belong to some $F_v^{j,n}$ with $j\leq k-1$, and at least one of them 
with $j\neq 0$. Thus, by the properties of $v$ we have $\gamma \leq (1+\alpha) \varepsilon_{k-1} \tau(vv^*) \tau(q)$. Since $\alpha$ was taken 
$\leq \varepsilon_0/4 \leq 1/4$, one gets $\gamma \leq \varepsilon_{k-1}$. 

Hence, the last term in $(9)$ is majorized by $\varepsilon_{k-1}\tau(q)/k + 
\varepsilon_{k-1}\tau(q)$. Since there are $k$ terms with $\ell=1$, obtained by taking $j=1, ..., k$, by summing up over $j$ in $(9)$ 
and combining with $(7)$, we deduce by applying $E_A$ to $(1)$ the following final estimate: 

$$
\|E_A(x)\|_1 \leq  \|E_A(x_0 \Pi_{i=1}^k v_i x_i)\|_1  + \Sigma_\ell \Sigma_i \|E_A(z_{0,i} \Pi_{j=1}^\ell
w_{i_j} z_{j,i})\|_1 \tag 10
$$
$$
\leq \varepsilon_k \tau(vv^*)+ (2^k - k -1) \varepsilon_{k-1} \tau(q) + (k+1) \varepsilon_{k-1} \tau(q)
$$
$$
= \varepsilon_k \tau(vv^*) + \varepsilon_k\tau(ww^*) = \varepsilon_k \tau((v+w)(v+w)^*). 
$$

At the same time, from $(8'')$, we get  
$$|\tau(u^jy_2)|=|\tau(v^j+w^j)y_2)|\leq |\tau(v^jy_2)|+|\tau(w^jy_2)| \tag 10'
$$
$$
\leq \varepsilon \tau(vv^*) + \varepsilon \tau(q)=\varepsilon \tau(uu^*).
$$

Altogether, this shows that $u=v+w\in \mycal W$. Since $u \geq v$ and $u\neq v$, this contradicts the maximality of $v \in {\mycal W}$. 

We conclude that $\tau(v^*v) >\tau(f)/2$. If we now take $\delta\leq \varepsilon/2^{n^2}$, 
then $\varepsilon_n = 2^{1+2+ ... + n}\delta $ $< 2^{n^2} \delta \leq \varepsilon$ and the statement follows.  

\hfill $\square$.

\proclaim{Lemma 4.7} Let $A_n \subset M_n$, $\text{\bf A}\subset \text{\bf M}$ be as in $4.1$. Let $X\subset \text{\bf M}\ominus \text{\bf A}$, $Y\subset \text{\bf A}$  
be countable sets and $f$ a non-zero projection in $\text{\bf A}$. 
Then there exists a partial isometry $w$ in $\text{\bf A}f$ such that $\tau(ww^*) \geq \tau(f)/2$, $E_{\text{\bf A}}(x)=0$  
and $\tau(w^jy)=0$, for all $n\geq k \geq 1$,  all 
$x \in X_w^{k,n}$ and all $y\in Y$. 
\endproclaim 
\noindent
{\it Proof}. Let $X = \{ x^k \}_k$, $Y=\{y^k\}_k$ be enumerations of the sets and let $x^k = (x_n^k)_n$, $y^k=(y^k_n)_n$ be 
representations of $x^k\in  \Pi_\omega M_n$, $y^k\in \Pi_\omega A_n$, which we can take so that 
$x_n^k\in M_n$ satisfy $E_{A_n} (x_n^k) = 0$, for all
$k$.  Let also $f_n \in \Cal P(A_n)$ be so that $f=(f_n)_n$. By applying Lemma 4.5 for the inclusion $A_n \subset
M_n$, the projection $f_n \in A_n$, the positive element $\varepsilon = 2^{-n}$, the integer $n$ and the finite sets $X_n = \{ x_n^k \mid k \leq
n\}$,  $Y_n=\{y_n^j \mid k\leq n\}$, we get a partial isometry $w_n$ in $A_nf_n$ such that
$\tau(w_n^* w_n)\geq \tau(f_n)/2$  and

$$
\| E_{A_n} (x) \|_1 \leq 2^{-n}, \forall x \in \underset
k \leq n \to \cup (X_n)_{w_n}^{k,n},
$$
$$
|\tau(w_n^jy)|\leq 2^{-n}, \forall y\in Y_n, 1\leq |j|\leq n. 
$$

But then $w = (w_n)$ clearly satisfies  the required conditions.  

\hfill
$\square$

\vskip .05in
\noindent {\it Proof of} 4.1. Since  $4.1 (b)$ is an immediate consequence of part $4.1. (a)$, we only need to prove the latter. 
To do this, we construct recursively a sequence 
of partial isometries $v_1, v_2, .... \in \text{\bf A}$ such that 
\vskip .05in 
\noindent
$(i)$ $v_{j+1}v_j^*v_j = v_j$ and $\tau (v_jv_j^*) \geq 1 - 1/2^j$, $\forall j \geq 1$. 

 \vskip .05in 
\noindent
$(ii)$ $E_{\text{\bf A}}(x)=0, \forall n\geq k \geq 1$,  
$\forall x \in X_{v_j}^{k,n}$

 \vskip .05in 
\noindent
$(iii)$ $|\tau (v_j^i y)$, $\forall n\geq 1$, $y\in Y$. 

\vskip .05in

Assume we have constructed $v_j$ for $j=1, ..., m$. If $v_m$ is a unitary element, then we let $v_j=v_m$ for all $j\geq m$. 
If $v_m$ is not a
unitary element, then let $f = 1 - v_m^*v_m \in \text{\bf A}$. Note that $E_{\text{\bf A}} (x') = 0$, for all $x' \in X' \overset \text{\rm def}
\to = \underset k \leq n \to \cup X_{v_m}^{k, n}$. 

If we apply now Lemma 4.7 to $\text{\bf A}\subset \text{\bf M}$, the projection  $f\in \text{\bf A}$, 
and to the countable set 
$X' \subset \text{\bf M}$, then we get a partial isometry $u \in \text{\bf A}f$ 
such that $\tau(uu^*)\geq \tau(f)/2$ and $E_{\text{\bf A}f} (x) = 0$ for all $x \in \underset k \leq n \to
\cup (X')_u^{k,n}$.  But then $v_{m+1}=v_m + u$ will satisfy  both $(i)$ and $(ii)$ for $j=m+1$. 

It follows now from $(i)$ that the sequence $v_j$ converges in the norm $\|$ $\|_2$ to a unitary element $v\in \text{\bf A}$, which due to $(ii)$ and $(iii)$ will satisfy the conditions 
required in part $(a)$ of $4.1$.  

\hfill $\square$

\vskip .05in

To deduce 4.2.1$^\circ$ from 4.1, we'll need the following:

\proclaim{Lemma 4.8} Let $u, v$ be unitary elements in a finite von Neumann algebra $M$ such that  
$\tau(u)=0$ and $\tau(v^j)=0$, for any non-zero integer $j$ with $ |j| \leq n-1$, for some $n \geq 2$. 
Assume  $\tau(x_0y_1x_1y_2 .... y_k x_k)=0$, for any $k\geq 1$ and any choice of $y_i \in \{u, u^*\}$,  $x_1, ..., x_{k-1} \in \{v^j \mid 1 \leq |j| \leq k-1\}$ and 
$x_0, x_k \in \{v^j \mid 1 \leq |j| \leq k-1\}\cup \{1\}$. Then we have: 

\vskip .05in 
\noindent
$(a)$ $\{u^* v^j u v^{-j} \mid j=1, 2, ..., n-1 \}$ are freely independent Haar unitaries in $M$. 
\vskip .05in 
\noindent 
$(b)$  $\|\Sigma_{j=0}^{n-1} v^j u v^{-j}\| \leq \sqrt{n-1} + 1$. 
\endproclaim 
\noindent
{\it Proof}.  Part $(a)$ is easy to check and we leave it as an exercise (see e.g. [AO] for similar calculations). 

To deduce part  $(b)$, recall that by a well known result  of Kesten ([Ke]),  if $w_1, ..., w_m$ are freely independent Haar unitaries 
in a finite von Neumann algebra $M$, then $\|\Sigma_{i=1}^m w_i\|=\sqrt{m}$. Thus, by $(a)$ we get: 
$$
\|\Sigma_{j=0}^{n-1} v^j u v^{-j} \| = \|u(1+ \Sigma_{j=1}^{n-1} u^* v^j u v^{-j})\|
$$
$$
= \|1+ \Sigma_{j=1}^{n-1} v^* v^j u v^{-j}\| \leq 1 + \|\Sigma_{j=1}^{n-1}u^*  v^j u v^{-j}\| = 1+\sqrt{n-1}. 
$$

\hfill 
$\square$ 

\noindent
{\it Proof of} 4.2. By Theorem 4.1 there exists a diffuse von Neumann subalgebra $B_0\subset \text{\rm A}=A^\omega$ such that 
any word with alternating letters from $\{u, u^*\}$, $B_0 \ominus \Bbb C1$, has trace $0$. Let $v\in B_0$ be a unitary element such that 
$\tau(v^j)=$ for $j=1, 2, .., n-1$ and $v^n=1$. Thus, if $\lambda\in \Bbb C$ is a primitive $n$'th root of $1$ then $v=\Sigma_{k=0}^{n-1} \lambda^{k}e_{k+1}$,  
where $e_k\in B_0$ are spectral projections of $v$ with $\tau(e_k)=1/n$, $\forall k$. An easy 
calculation shows that $n^{-1}\Sigma_{j=0}^{n-1} v^j u v^{-j} = \Sigma_{k=1}^n e_k u e_k$.  But then 4.2.1$^\circ$  
follows from $4.8$ $ (b)$. 

To prove 4.2.2$^\circ$ let $B_0\subset A^\omega$ be free with respect to the von Neumann algebra $\Bbb C e + \Bbb C(1-e)$. 
Then the calculation of norms in  [V2] shows that if $q$ is any projection 
of trace $1/n$ in $B_0$ with $1/n\leq \tau(e)$, then $\|qeq-\tau(e)q\|\leq 2/\sqrt{n}$. By applying again [V2] for $p\in B_0$ with $\tau(p)=1/2$, we get 
$\|e(p-\tau(p)1)e\| = (\tau(e)(1-\tau(e)))^{1/2}$ and thus $\|pep\|\leq \tau(p) + (\tau(e)(1-\tau(e)))^{1/2}=1/2 + (\tau(e)(1-\tau(e)))^{1/2}$. Since $\tau(1-p)=1/2$ as well, 
we get similarly $\|(1-p)e(1-p)\|\leq 1/2 + (\tau(e)(1-\tau(e)))^{1/2}$

\hfill
$\square$ 

\noindent
{\it Proof of} 4.3. By Proposition 2.2, in order to prove Corollary $4.3$, it is sufficient to prove that any projection 
$e\in \text{\bf M}$ whose expectation on $\text{\bf A}$ is a scalar multiple of some projection $f\in \text{\bf A}$, can be paved. But 
this is indeed the case, because $\text{\bf A}f\subset f\text{\bf M}f$ is itself an ultraproduct of singular inclusions, for which 4.2 applies.

\hfill 
$\square$

\heading 5. Final remarks
\endheading

\noindent
{\bf 5.1. Examples of singular MASAs}.  Dixmier's first examples of singular MASAs $A$ in II$_1$ factors $M$  ([D1]), were constructed  
from group-subgroup situations, $H \subset G$, as $A=L(H)\subset L(G)=M$, with $G$ infinite conjugacy class (ICC) and $H\subset G$ an 
abelian subgroup satisfying certain conditions. These conditions are met for instance by wreath product inclusions groups $H \subset G = K \wr H$, with $H$ infinite abelian 
and $K$ non-trivial and by  the the inclusions $L(\Bbb Z) \subset L(\Bbb Z * \Gamma_0)$, for any non-trivial group $\Gamma_0$. 
Another criterion for singularity of MASAs in factors 
was found in [P2]. It can be used to recover the previous examples, as well as others. 
It shows for instance that $A=L^\infty([0,1])$ is singular in $A * N$ for any finite von Neumann algebra $N$. It also shows that  
the group algebra $A=L(H)$ is singular in any crossed product II$_1$ factor 
$M=B^{\otimes H} \rtimes H$, arising from a Bernoulli action $H \curvearrowright B^{\otimes H}$, 
for any non-trivial finite von Neumann ``base''-algebra $B$. In fact, by (3.1 in [P11]), all these MASAs $A$ are singular in the following stronger sense: If $u\in \Cal U(M)$ 
is so that $uAu^* \cap A$ is diffuse, then $u\in A$. This {\it absorption} phenomenon from ([P11]) is actually valid for any 
inclusion $L(H)\subset M=N \rtimes H$, arising from a mixing action of $H$ on a finite von Neumann algebra $N$. 

Another strengthening of the notion of singularity for a MASA $A \subset M$ was emphasized in [P3] and it requires that the only automorphisms of $M$ 
that normalize $A$ are the inner automorphisms Ad$(u)$ with $u\in \Cal U(A)$. Such MASAs were called ultrasingular in [P3], 
but we will call them {\it supersingular} from now on, because they have the property that any two automorphisms of $M$ that coincide on $A$ 
must differ by some Ad$(u)$, with $u\in \Cal U(A)$. Equivalently, embeddings with same range of $M$ into another algebra are uniquely 
determined by their values on $A$. It was shown in [P3] that any II$_1$ factor $M$ whose outer automorphim group 
Out$(M)$ is countable (e.g. if $M$ has property T,  by [C2]), do have supersingular MASAs.  

We note here that results from (Section 4 and 5 of [P12]) show in particular 
that if one reduces the singular $MASA$, $A=L(\Bbb Z) \subset L([0, 1]^\Bbb Z)\rtimes \Bbb Z=R$, by a 
projection $p\in A$ which is not fixed by any ``rotation'' by a character $\gamma \in \hat{\Bbb Z}$, then $Ap\subset pRp\simeq R$ is 
supersingular. Moreover, if $p, q\in A$ are not conjugate by such a rotation, then $Ap, Aq$ are distinct singular MASAs in $R$. More 
precisely, we have:

\proclaim{5.1.1. Theorem [P12]}  Let $H$ be a torsion free abelian group $($such as $H=\Bbb Z)$ and $H \curvearrowright X=X_0^H$ a Bernoulli $H$-action. 
Let $R=L^\infty(X) \rtimes H$ and denote $A=L(H)$. If $p, q\in A$ are non-zero projections and  
$\theta: pRp \simeq qRq$ is an isomorphism carrying $Ap$ onto $Aq$,  then there exists a character $\gamma \in \hat{H}$ such that $\theta$ is 
the restriction of $\theta_\gamma\in \text{\rm Aut}(R)$ to $pRp$. Moreover, 
the only automorphisms of $pRp\simeq R$ that normalize $Ap$ are the restrictions of the automorphisms $\theta_\gamma$ 
that satisfy $\gamma(Y)=Y$ (a.e.), where $Y\subset \Bbb T$ is the subset with characteristic function $\chi_Y=p$. In particular, if $\{p_t  \mid t\in (0, 1]\}$ is a family 
of projections in $L(H)$ with $\tau(p_t)=t$, then  $Ap_t\subset  p_tRp_t\simeq R$ provide a  
family of distinct singular MASAs  in the hyperfinite $\text{\rm II}_1$ factor, which are supersingular for $t\not\in \Bbb Q$. 
\endproclaim
\noindent

\vskip .05in 
\noindent
{\bf 5.2. Characterizations of singularity for MASAs}. Another strengthening of singularity for MASAs was discovered in ([P4]), 
where it is shown that if $A$ is a diffuse abelian von Neumann algebra and $N$ is any finite 
von Neumann algebra, then $A$ is maximal amenable (equivalently, maximal injective) in $A * N$.  

We notice in 5.2.1 below an immediate consequence of Theorem 4.1, showing that for any singular MASA 
$A\subset M$, the ultrapower $A^\omega$ is maximal amenable in $M^\omega$, i.e., if $A^\omega \subset P \subset M^\omega$ for some amenable 
von Neumann algebra $P$, then $P=A^\omega$.  Moreover, any $P\subset M^\omega$ that contains $A^\omega$ and has countable dimension both 
as a left and right Hilbert module, must coincide with $A^\omega$. 

We also provide an alternative characterization of 
singularity for MASAs in terms of moments, as those MASAs that contain Haar unitaries which are asymptotically free 
with respect to sets perpendicular to it. For this to happen, asymptotic $4$-independence is in fact sufficient. This should be compared to Theorem 3.9 
where it was shown that asymptotic 2-independence occurs for any MASA, and to 5.3.1 below, which shows that in fact in 
arbitrary MASAs asymptotic 3-independence occurs as well. 

\proclaim{5.2.1. Theorem} Let $A_n\subset M_n$ be a sequence of MASAs in finite von Neumann algebras $M_n$ and denote 
$\text{\bf A}=\Pi_\omega A_n$, $\text{\bf M}=\Pi_\omega M_n$. The following are  equivalent: 
\vskip .05in
$1^\circ$  There exists a sequence of projections $p_n \in A_n$ such that $\underset n \rightarrow \omega \to \lim \tau(p_n)=1$ 
and $A_np_n$ is singular in $p_nM_np_n$, $\forall n$. 
\vskip .05in 
$2^\circ$ $\text{\bf A}$ is singular in $\text{\bf M}$; 
\vskip .05in 
$3^\circ$ $\text{\bf A}$ is maximal amenable in  $\text{\bf M}$; 
\vskip .05in 
$4^\circ$  If $\Cal H\subset L^2\text{\bf M}$ is a Hilbert $\text{\bf A}$-bimodule with $\dim \Cal H_{\text{\bf A}}, \dim _{\text{\bf A}}\Cal H\leq \aleph_0$  
 $($i.e., $\exists X\subset \Cal H$ countable such that ${\text{\rm sp}}\text{\bf A} X$   
and $\text{\rm sp} X \text{\bf A}$ are dense in $\Cal H)$, 
then $\Cal H$ is of the form $L^2\text{\bf A}p$, for some $p\in \Cal P(\text{\bf A})$. 
In particular, $\text{\bf A}$ is maximal among subalgebras $ P\subset \text{\bf M}$ that contain $\text{\bf A}$ and 
have the property that $L^2P$ is countably generated both as a left and right Hilbert $\text{\bf A}$-module. 
\vskip .05in 
$5^\circ$ Given any countable set $X\subset \text{\bf M}\ominus \text{\bf A}$, there exists $B_0\subset \text{\bf A}$ 
diffuse such that $B_0$, $X$  are free independent relative to $\text{\bf A}$.
\vskip .05in
$6^\circ$ Given any countable set  $X\subset \text{\bf M}\ominus \text{\bf A}$, there exists $B_0\subset \text{\bf A}$ 
diffuse such that $B_0$, $X$  are $4$-independent.
\vskip .05in
$7^\circ$ Given any selfadjoint element $x\in \text{\bf M}\ominus \text{\bf A}$, there exists $B_0\subset \text{\bf A}$ 
diffuse such that $B_0$, $\{x\}$  are $4$-independent.
\vskip .05in 
Moreover, if $A_n \subset M_n$ are all equal to the same MASA $A\subset M$ then the above are also equivalent to $A\subset M$ 
being singular and in $6^\circ, 7^\circ$ 
it is sufficient to take $X$, resp. $\{x\}$ inside $M\ominus A$. 
\endproclaim 
\noindent
{\it Proof}. If  $u_n  \in \Cal N_{M_n}(A_n)$, then $u=(u_n)_n$ normalizes $\text{\bf A}$ as well, 
acting non-trivially iff $u \not\in \text{\bf A}$. Also, if $e_n$ is the maximal projection in $A_n$ with the 
property that $u_ne_n\in A_n$, then $u$ acts nontrivially on $\text{\bf A}$ iff $\underset n \rightarrow \omega \to \lim \tau(e_n)=0$. 
This shows that $2^\circ \implies 1^\circ$. The converse is implicit in [P10] (due to Remark 5.2 in [P3]). Indeed, for  
if $u\in \Cal N_{\text{\bf M}}(\text{\bf A})$ is not in $\text{\bf A}$, then there exists a non-zero projection $q\in \text{\bf A}$ such that 
$uqu^* q =0$ and $u$, $q$ can be represented by sequences $u=(u_n)_n$, $q=(q_n)_n$, with $u_n$ unitaries in $M_n$, 
$q_n$ projections in $A_n$, such that $u_n q_n u_n^* q_n =0$. Moreover, since $\underset n \rightarrow \omega \to \lim \tau(p_n)=1$, we may 
assume $q_n \leq p_n$ and $u_nq_n u_n^*\leq p_n$. 
But by ([P10]), by the singularity of $A_np_n \subset p_nM_np_n$, for each $n$ there exists a 
unitary element $v_n \in q_nM_nq_n$ such that $\|E_{A_n}(u_nv_nu_n^*)\|_2 \leq \|q_n\|_2/n$. Thus, $v=(v_n)_n\in \text{\bf A}$ satisfies $uvu^* \perp \text{\bf A}$, 
a contradiction. Thus $1^\circ, 2^\circ$ are equivalent. 

The implication $1^\circ \implies 5^\circ$ is shown in Theorem $4.1. (a)$, and $5^\circ \implies 6^\circ \implies 7^\circ$ are trivial. 
To see that $6^\circ \implies 1^\circ$, assume there exist  $v_n\in M_n$ partial isometries  
such that $v_nv_n^*, v_n^*v_n$ are mutually orthogonal projections in $A_n$  and $v_nA_nv_n^*=A_nv_nv_n^*$. If we denote $v=(v_n)_n$ 
and $u\in \text{\bf A}$ would be 
a Haar unitary that's $4$-independent with respect to $X=\{v, v^*\}$, then the equality $vuv^*u^*=u^*vuv^*$ (due to abelianess of $A^\omega$) 
implies $0\neq \tau(vv^*)=\tau(vu^*v^*uvuv^*u^*)=0$, a contradiction. Taking $X=\{v+v^*\}$, this actually proves $7^\circ \implies 1^\circ$ as well. 

The implication $3^\circ \implies 2^\circ$ is trivial. To prove the converse, note that if $N$ is any von Neumann algebra that strictly contains $\text{\bf A}$, then there exist  
two orthogonal projections $p_1, p_2 \in \text{\bf A}$ that are equivalent via some partial isometry $v$ in $N$ (exercise!). If $q=2^{-1}(p_1 + p_2 + v+v^*)$, 
then $q$ is a projection in $N$ such that $E_{\text{\bf A}}(q)=2^{-1}p$, where $p=p_1+p_2$. By Theorem 4.1, 
there exists a diffuse von Neumann subalgebra $B_0\subset \text{\bf A}$ such that any alternating word in $B_0\ominus \Bbb C$ 
and $q-2^{-1}p$ has trace $0$. Thus, the algebras $B_0p$ and $\Bbb C q + \Bbb Cp$ are free independent, implying that if $u\in B_0p$ is a Haar unitary 
then $u$ and $(v+v^*)u(v+v^*)$  generate a copy of the free group factor $L(\Bbb F_2)$ inside $pNp$. Since amenability is a hereditary property (by [S]), 
this shows that $N$ cannot be amenable.  

$4^\circ \implies 2^\circ$ is trivial (because if $\text{\bf A}$ is not singular and $u\in \Cal N_{\text{\bf M}}(\text{\bf A})\setminus \Cal U(\text{\bf A})$, 
then   the von Neumann algebra $P$ generated by $u$ and $\text{\bf A}$ has the countable set $X=\{u^n \mid n\in \Bbb Z\}$ satisfying sp$X \text{\bf A}$, 
sp$\text{\bf A}X$ dense in $P\neq \text{\bf A}$. Finally, to prove $2^\circ \implies 4^\circ$, assume $X\subset \Cal H\ominus \text{\bf A}$ is a separable subspace 
such that the span of both $X \text{\bf A}$ and $\text{\bf A} X$ are $\| \ \|_2$-dense in $\Cal H\ominus \text{\bf A}$. By $4.1. (a)$, there exists $B_0 \subset \text{\bf A}$ 
diffuse such that $B_0$ is free independent to $X$ relative to $\text{\bf A}$. In particular, given any Haar unitary $u\in B_0$, we have $E_{\text{\bf A}}(x_1^*ux_2)=0$, 
for all $x_1, x_2 \in X$. Thus, $uX \perp \text{\rm sp}X \text{\bf A}=\Cal H \ominus \text{\bf A}$, a contradiction.  

The last part of the statement, when all $A_n\subset M_n$ are assumed to be equal, is now trivial. 
\hfill 
$\square$ 

\vskip .05in 
\noindent
{\bf 5.3. Controlling moments through incremental patching}. In Theorem 3.9, we have proved 
that if $A \subset M$ is an arbitrary MASA, then for any countable $X \subset M^\omega \ominus A^\omega$, there exists a diffuse 
von Neumann subalgebra  $B_0\subset A^\omega$ such that $B_0$ is $2$-independent with respect to $X$. 
We chose to prove this  through a ``global'' construction 
of finite dimensional approximations of such a $2$-independent $B_0$. But we can also prove 
this result differently, through 
the method used in the proofs of the previous section, and which consists in  
controlling the moments  incrementally, by patching ``infinitesimal  pieces'' of an appropriate Haar unitary. This method does use  
a technical result from Section 3, namely property $(d')$ of Theorem 3.6, but which was   
already known since [P1] (see also A.1 in [P5]): If $M$ is a finite von Neumann algebra, $A \subset M$ a MASA and 
$X\subset M\ominus A$ a finite set, then given any $\varepsilon > 0$, there exists a non-zero projection $q\in A$ such that $\|qxq\|_1 \leq \varepsilon \tau(q)$, 
$\forall x\in X$. 

In fact, as shown in 5.3.1 below, the ``incremental patching'' method can be used to obtain a slightly stronger result  for arbitrary MASAs $A\subset M$, showing that   
one can construct separable, diffuse von Neumann subalgebras 
$B_0\subset A^\omega$ that are $3$-independent with respect to any given countable set 
$X\subset M^\omega\ominus A^\omega$.  As we saw in Theorem 5.2.1, this is the best one can do for an arbitrary MASA, 
as existence of a $B_0$ that's $4$-independent with respect to any given countable set $X\subset M\ominus A$ forces $A$ to be singular (in which case   
$B_0$ can even be chosen free independent with respect to the given $X$).

\proclaim{5.3.1. Theorem}  Let $M_n$ be a sequence of finite 
factors with $\text{\rm dim} M_n \rightarrow \infty$ and for each $n$, let $A_n\subset M_n$ be a MASA. 
Denote by $\text{\bf A}=\Pi_\omega A_n \subset \Pi_\omega M_n$. 
Let  $Q \subset \Pi_\omega M_n$ be an arbitrary separable von Neumann subalgebra such that $E_{\text{\bf A}}(Q) = \text{\bf A}\cap Q$, i.e. $Q$ and 
$\text{\bf A}\subset \Pi_\omega M_n$ make a commuting square, and denote $B_1=\text{\bf A}\cap Q$. There exists a diffuse von Neumann 
subalgebra $B_0\subset \text{\bf A}$ such that $B_0$ is $3$-independent to $Q\ominus B_1$, more preciseley:  
$\tau(xa)=0$, $\forall x\in Q, a\in B_0\ominus \Bbb C$; 
$\tau(x_1a_1x_2a_2)=0$, $\tau(x_1a_1x_2a_2x_3a_3)=0$, $\forall x_i \in Q\ominus B_1$, $a_i \in B_0 \ominus \Bbb C$ 
$($N.B.: the odd level independence relations follow from the even ones$)$. 
\endproclaim
\noindent
{\it Proof}. We proceed along the lines of the proofs of Lemmas 4.5, 4.7 and Theorem 4.1, from the previous section. 
If $F$ is a subset in a von Neumann algebra and $v$ a partial isometry with $vv^*=v^*v$, then we denote 

$$
F^k_{v,n}=\{\Pi_{j=1}^k v^{i_j} x_j  \mid x_j \in F, 1\leq |i_j| \leq n, 1\leq j \leq k\} \}. 
$$
We first prove the following: 
\vskip .05in
{\it Fact}. Let $M$ be a finite von Neumann algebra and 
$A\subset M$ a MASA. Given any finite set $F\subset M\ominus A$, with $\|x\| \leq 1$, $\forall x\in F$, any $n\geq 1$ and any $\delta > 0$, there exists a 
Haar unitary $v\in A$ such that  $|\tau(x)|\leq \delta$, $\forall x\in \overset
3 \to{\underset k = 1 \to \cup} F^k_{v,n}$. 

\vskip .05in 
To prove this, denote by ${\mycal W}=\{ v \in A \mid vv^*\in \Cal P(A),  |\tau(x)| \leq \delta \tau(v^*v), 
\forall x \in \cup_{k=1}^3 F^k_{v,n}, \tau(v^m)=0, \forall m\neq 0 \}$.  Endow ${\mycal W}$ with the order $\leq$ in
which $w_1 \leq w_2$ iff $w_1 = w_2 w_1^*w_1$. $({\mycal W} , \leq)$
is then clearly inductively ordered.  Let $v$ be a maximal element
in ${\mycal W}$. Assume $\tau(v^*v) < 1$ and denote $p = 1 - v^*v$.  If $w$ is a
partial isometry in $Ap$  and  $u = v + w$, 
then by using that $u^{i_j}=v^{i_j}+ w^{i_j}$ and expanding  $x =u^{i_1}x_1 u^{i_2}x_2 ... u^{i_k}x_k \in F^k_{u,n}$, $k = 1, 2, 3$,  
as a binomial product, we have $\tau(x) = \tau(\Pi_{j=1}^k v^{i_j} x_j) + \Sigma \tau( ... x_{j-1}w^{i_j}x_{j} .... )$, and thus 
$$
|\tau(x)| \leq  |\tau(\Pi_{j=1}^k v^{i_j} x_j)| + \Sigma |\tau( ... x_{j-1}w^{i_j}x_{j} .... )|, \tag 1
$$
where the sum is taken over all terms that have at least one occurrence of $w^{i_j}$. Since $v\in {\mycal W}$, we have 
$|\tau( \Pi_{j=1}^k v^{i_j} x_j)|\leq \delta \tau(vv^*)$. We will prove that we can choose $w\neq 0$ 
so that the summation on the right hand side of $(1)$ is majorized by $\delta \tau(ww^*)$, giving $|\tau(x)|\leq \delta \tau(vv^*) + \delta \tau(ww^*)=\delta \tau(uu^*)$. 
This will contradict the maximality of $v$, thus showing that $vv^*=1$, i.e $v$ is a Haar unitary in $A$. 
We construct $w$ by first making an appropriate choice for its 
support projection $q=ww^*$, then choosing $w$ as an appropriate Haar unitary in $Aq$. 

In order to estimate the summation $ \Sigma |\tau( ... x_{j-1}w^{i_j}x_{j} .... )|$ in $(1)$, 
note the following: in case $k=1$ the sum has just one member, being of the form  $|\tau(w^jy)|$, for some $1\leq |j|\leq n$, $y\in X$; 
in case $k=2$, the sum has three terms, being of the form  
$$
|\tau(w^{j_1}x_1v^{j_2}x_2)| +  |\tau(v^{j_1} x_1 w^{j_2}x_2)|  + |\tau(w^{j_1}x_1w^{j_2}x_2)|;  \tag 2
$$
in case $k=3$, the sum has seven terms, being of the form 
$$
|\tau(w^{j_1}x_1 v^{j_2} x_2 v^{j_3}x_3)| + |\tau(v^{j_1}x_1w^{j_2}x_2v^{j_3}x_3)| + |\tau(v^{j_1}x_1v^{j_2}x_2w^{j_3}x_3)| \tag 3
$$
$$
+ |\tau(w^{j_1}x_1 w^{j_2} x_2 v^{j_3}x_3)| + |\tau(v^{j_1}x_1w^{j_2}x_2w^{j_3}x_3)| + |\tau(w^{j_1}x_1v^{j_2}x_2w^{j_3}x_3)| 
$$
$$
+| \tau(w^{j_1}x_1 w^{j_2} x_2 w^{j_3}x_3)|. 
$$

Now note that for each summand for which we have $2$ or $3$ appearances of non-zero powers of $w$ in the above sums (one term for $k=2$ and 
four terms for $k=3$), such appearances must be consecutive, i.e. they will be of the form $|\tau(....w^iyw^j... )|$, for some $i, j \neq 0$, $y\in F\subset M\ominus A$ 
(for one of the terms, one uses the equality $\tau(w^{j_1}x_1v^{j_2}x_2w^{j_3}x_3)=\tau(x_1v^{j_2}x_2w^{j_3}x_3w^{j_1})$). 
If $q=ww^*$, then for each one of these terms  we have $|\tau(....w^iyw^j... )| \leq \|qyq\|_1$. By (2.1 in [P1]), or (A.1 in [P5]), or by using $3.6 (d')$ in this paper, 
applied to the MASA $Ap\subset pMp$ 
and the set $pFp\subset pMp\ominus Ap$, one can choose the  
non-zero projection  $q\in Ap$ such that $\|q y q\|_1 \leq 2^{-3}\delta \tau(q)$, $\forall y \in pFp$. It thus follows that the sum of terms having two or more 
appearances of powers of $w$ are majorized by $2^{-1}\delta \tau(q)$ (because there is one such term when $k=2$ and four when $k=3$). 

All remaining terms 
and the case $k=1$ have just one occurrence of $w^j$, $j\neq 0$, i.e are of the form $|\tau(y_1w^jy_2)|=|\tau(w^jE_A(qy_2y_1q))|$, for some $y_1, y_2 \in M$, 
$1\leq |j| \leq n$. There are $k$ many such terms for each $k= 1, 2,  3$. 
Let's denote by $Y_0$ the set of all $y_1, y_2$ which appear this way, and note that this is a finite set in $M$. Thus  $Y=E_A(qY_0\cdot Y_0q)$ 
is finite as well. It is sufficient  to find now a Haar unitary in $Aq$ such that $|\tau(w^jy)|\leq 2^{-8}\delta \tau(q)$, $\forall y\in Y$, $1\leq |j|\leq n$, because 
then the sum of the $k$ terms will be majorized by $2^{-1}\delta \tau(q)$ which added up to the quantity $2^{-1}\delta \tau(q)$ that majorizes 
the terms with at least $2$ occurrences of powers of $w$ gives that for all $x\in \cup_{k=1}^3 F^k_{u,n}$, we have $|\tau(x)|\leq \delta \tau(uu^*)$. 
Since $Aq$ is diffuse, it contains a separable diffuse subalgebra $A_0\subset Aq$, which is isomorphic to $L^\infty(\Bbb T)$ with the Lebesgue 
measure corresponding to $\tau(q)^{-1}\tau_{|A_0}$. Let then $w_0\in A_0$ be a Haar unitary generating $A_0$. Since $\{w_0^m\}_m$ tends to $0$ 
in the weak operator topology and $Y\subset A$ is a finite set, there exists $n_0\geq n$ such that $|\tau(w_0^m y)|\leq 2^{-4}\delta \tau(q)$, 
for all $y\in Y$ and $|m| \geq n_0$. But then $w=w_0^{n_0}$ is still a Haar unitary and it satisfies all the required conditions. 

This ends the proof of the {\it Fact}. 

\vskip .05in 
With this in hand, we proceed as follows: 
Let $X_0\subset Q$ be $\| \ \|_2$-dense countable subset and denote by $X=\{ y/\|y\| \mid y=x-E_{\text{\bf A}}(x), x\in X_0\setminus \text{\bf A} \}$. 
Note that $X$ is a countable subset of $\text{\bf M}\ominus \text{\bf A}$ and each element in $X$ has operator norm equal to $1$. Write $X$ as a sequence $\{x_n\}_n$. 
For each $n$ we now apply the above {\it Fact} to the set $F_n=\{x_1, ..., x_n \}$ and $\delta = 1/n$, to get a Haar unitary $v_n\in \text{\bf A}$ such that 
$$
|\tau(\Pi_{j=1}^k v_n^{i_j} x_{t_j})| < 1/n, \forall  |i_j|, t_j \in  \{1, ..., n\}, k= 1, 2, 3. \tag 4
$$

Let $x_n=(x_{n,m})_m$, $v_n=(v_{n,m})_m$ be representations of 
the $x_n$'s and $v_n$'s with $x_{n,m}\in M_m$, $\|x_{n,m}\|\leq 1$, $v_{n,m}\in \Cal U(A_m)$. Thus, $(4)$ and the fact that $v_n$ are 
Haar unitaries, translates into  
$$
\underset m \rightarrow \omega \to \lim \tau(v^k_{n,m})=0, \forall k\neq 0; 
\underset m \rightarrow \omega \to \lim |\tau(\Pi_{j=1}^k v_{n,m}^{i_j} x_{t_j, m})| < 1/n, \forall  |i_j|, t_j \in  \{1, ..., n\}. \tag 5
$$

Let $V_n$ of $\omega$ denote the set of all $m \in \Bbb N$ with the property that  
$$
|\tau(v_{n,m}^k)| < 1/n, 1\leq |k|\leq n; |\tau(\Pi_{j=1}^k v_{n, m}^{i_j} x_{t_j, m})| < 1/n, \forall  |i_j|, t_j \in  \{1, ..., n\}. \tag 6
$$
From $(5)$, it follows that $V_n$ corresponds to a closed-open neighborhood of 
$\omega$ in $\Omega$, under the identification $\ell^\infty\Bbb N=C(\Omega)$. With this in mind, define recursively $W_0=\Bbb N$, 
$W_{n+1}=W_n\cap V_{n+1}\cap \{ n \in \Bbb N \mid n > \min W_n\}$ and note that $\{W_n\}_n$ this way defined is a strictly decreasing 
sequence of neighborhoods of $\omega$ satisfying $W_n \subset \cap_{j\leq n} V_j$. 

Finally, let $u=(u_m)_m \in \text{\bf A}$ be defined by $u_m=v_{n,m}$, for $m\in W_n \setminus W_{n-1}$, $n\geq 1$. 
It is then immediate to check that $u$ is a Haar unitary element in $\text{\bf A}$ and that the von Neumann algebra $B_0$ it generates 
satisfies the required $3$-independence conditions. 

\hfill 
$\square$

\vskip .05in 
\noindent
{\bf 5.4. Exact paving size for ultraproducts of singular MASAs}. The order of magnitude of the paving size 
in Corollary 4.3 should be $\varepsilon^{-2}$, for any $x$, not only for $x=v$ unitary element with $E_{\text{\bf A}}(v)=0$ and projections that 
expect on scalars (i.e., the cases covered by Cor 4.2). We pose here two questions which, 
if  answered in the affirmative, would imply this fact: 
\vskip .05in 
$(a)$ Can any $x$ with  $E_{\text{\bf A}}(x)=0$ and norm $\leq 1/2$ (or of norm $\leq c$ for an even smaller universal constant $c>0$) 
be written as a convex combination of unitaries having $0$-expectation on $\text{\bf A}$? If so, then 4.2.1$^\circ$ would imply that 
$\text{\rm n}(x, \varepsilon)$ is majorized by a constant multiple of $\varepsilon^{-2}$ for any $x\in M^\omega$. 
\vskip .05in 
$(b)$ Is it true that if $M$ is a II$_1$ factor and $x=x^*\in M\ominus \Bbb C$, 
$u\in \Cal U(M)$ a Haar unitary, such that $\tau(u^{i_1}x u^{i_2} x ...)=0$, for any alternating word with $i_j\neq 0$, then 
$\|\Sigma_{i=1}^n u^i x u^{-i}\|$ has order of magnitude $\sqrt{n}$ ? Again, if this would hold true in this generality, then we would not need 
Proposition 2.3 at the end of the proof of Theorem 4.3, the result following directly from 4.1$(a)$, with the estimate $\varepsilon^{-2}$ for the 
order of magnitude of the paving size.

\vskip .05in 
\noindent
{\bf 5.5. A conjecture generalizing  Kadison-Singer}. While we have not been able to settle the classic Kadison-Singer 
problem in its equivalent formulations of Theorem 2.2, i.e., 
by proving that one can pave all elements in $R^\omega$ (resp. in $\text{\bf M}=\Pi_\omega M_{n \times n}(\Bbb C)$) 
over its MASA $D^\omega$ (resp. $\text{\bf D}=\Pi_\omega D_n$),  we believe this is true and  that in fact the following more general conjecture holds true: 

\vskip .05in 
\noindent
{\it 5.5.1. Conjecture}: Given any sequence of MASAs in finite factors, $A_n \subset M_n$, the ultraproduct inclusion $\Pi_\omega A_n \subset \Pi_\omega M_n$ 
has the Kadison-Singer (equivalently, the paving) property.

\head  References \endhead

\item{[Ak]} C. Akemann: {\it The dual space of an operator algebra}, Trans. Amer. 
Math. Soc. {\bf 126}  (1967), 286-302. 

\item{[AkA]} C. Akemann, J. Anderson: ``Lyapunov theorems for 
operator algebras'', Mem. AMS {\bf 94} (1991). 

\item{[AkO]} C. Akemann, Ostrand: {\it Computing norms in group $C^*$-algebras}, Amer. J. Math. {\bf 98} (1976), 1015-1047. 

\item{[A1]} J. Anderson: {\it Extensions, restrictions and
representations of states on}
C$^*$-{\it alge-\newline bras}, Trans. Amer.
Math. Soc. {\bf 249} (1979), 303-329. 

\item{[A2]} J. Anderson: {\it Extreme points in sets of positive linear maps on $\Cal B(\Cal H)$}, Jour. Functional Analysis 
{\bf 31} (1979), 195-217. 

\item{[BeHKW]} K. Berman, H. Halpern, V. Kaftal, G. Weiss: {\it Matrix norm inequalities and the relative Dixmier property} 
Integral Equations and Op. Theory {\bf 11} (1988), 28-49. 

\item{[BT]} J. Bourgain, L. Tzafriri: {\it On a problem of
Kadison and Singer}, J. Reine Angew. Math. {\bf 420} (1991), 1-43.

\item{[CaFTW]} P. Casazza, M. Fickus, J. Tremain, E. Weber: {\it 
The Kadison-Singer Problem in Mathematics and Engineering: a detailed account}, in ``Operator theory, operator algebras, and applications'' 
299-355, Contemp. Math., Vol. {\bf 414}, Amer. Math. Soc., Providence, RI, 2006. 

\item{[C1]} A. Connes: {\it Classification of injective factors},
Ann. of Math., {\bf 104} (1976), 73-115.

\item{[C2]} A. Connes: {\it A type} II$_1$ {\it factor with countable fundamental group}, 
J. Operator Theory {\bf 4} (1980), 151-153. 

\item{[CFW]} A. Connes, J. Feldman, B. Weiss:
{\it An amenable equivalence relation is generated by a single
transformation}, Erg. Theory Dyn. Sys.  {\bf 1} (1981),
431-450.

\item{[D1]} J. Dixmier: {\it Sous-anneaux ab\'eliens maximaux dans les facteurs de type fini}, Ann. of Math. {\bf 59} (1954), 279-286. 

\item{[D2]} J. Dixmier: ``Les alg\'ebres d'op\'erateurs dans l'espace hilbertien'', Gauthier-Vill-\newline ars, Paris 1957, 1969. 

\item{[F]} J. Feldman: {\it Nonseparability of certain finite factors},  Proc. Amer. Math. Soc. {\bf 7} (1956), 23-26. 

\item{[KS]} R.V. Kadison, I.M. Singer: {\it Extensions of pure
states}, Amer. J. Math. {\bf 81} (1959), 383-400.

\item{[KR]} R.V. Kadison, J. Ringrose: ``Fundamentals of the 
theory of operator algebras'', Pure and Applied Mathematics, {\bf 100}, 
Academic Press, Inc. New York, 1983

\item{[Ke]} H. Kesten: {\it Symmetric random walks on groups}, Trans. Amer. Math. Soc. {\bf 92} (1959), 336-354. 

\item{[MvN1]} F. Murray, J. von Neumann:
{\it On rings of operators}, Ann. Math. {\bf 37} (1936), 116-229.

\item{[MvN2]} F. Murray, J. von Neumann: {\it Rings of operators
IV}, Ann. Math. {\bf 44} (1943), 716-808.

\item{[vN]} J. von Neumann: {\it Einige satze uber messbare abbildungen}, Ann. of Math. {\bf 33} (1932), 574-586. 

\item{[OW]} D. Ornstein, B. Weiss: {\it Ergodic theory of
amenable group actions I. The Rohlin Lemma} Bull. A.M.S. (1) {\bf 2}
(1980), 161-164.

\item{[PiP]} M. Pimsner, S. Popa: {\it Entropy and index for subfactors}, Annales Scient. Ecole Norm. Sup. {\bf 19} (1986), 57-106. 

\item{[P1]} S. Popa: {\it On a problem of R.V. Kadison on maximal
abelian *-subalgebras in factors}, Invent. Math., {\bf 65} (1981),
269-281.

\item{[P2]} S. Popa: {\it Orthogonal pairs of *-subalgebras in
finite von Neumann algebras}, J. Operator Theory, {\bf 9} (1983),
253-268.

\item{[P3]} S. Popa:  {\it Singular maximal abelian *-subalgebras in
continuous von Neumann algebras}, J. Funct. Analysis, {\bf 50}
(1983), 151-166. 

\item{[P4]} S. Popa: {\it Maximal injective subalgebras in factors
associated with free groups}, Advances in Math., {\bf 50} (1983),
27-48.

\item{[P5]} S. Popa: {\it Classification of amenable subfactors of
type} II, Acta Mathematica, {\bf 172} (1994), 163-255.

\item{[P6]} S. Popa: {\it Free independent sequences in type} II$_1$ {\it factors
and related problems}, Asterisque, {\bf 232} (1995), 187-202.

\item{[P7]} S. Popa: {\it The relative Dixmier property for inclusions
of von Neumann algebras of finite index}, Ann. Sci. Ec. Norm. Sup.
{\bf 32} (1999), 743-767.

\item{[P8]} S. Popa: {\it On a class of type} II$_1$ {\it factors with
Betti numbers invariants}, Ann. of Math {\bf 163} (2006), 809-899
(math.OA/0209310).

\item{[P9]} S. Popa: {\it On the distance between MASA's in type}
II$_1$ {\it factors}, Fields Institute Communications, {\bf 30}
(2001), 321-324.

\item{[P10]} S. Popa: {\it Strong Rigidity of} II$_1$ {\it Factors
Arising from Malleable Actions of $w$-Rigid Groups} I, Invent. Math.,
{\bf 165} (2006), 369-408. (math.OA/0305306).

\item{[P11]} S. Popa: {\it Strong Rigidity of} II$_1$ {\it Factors
Arising from Malleable Actions of $w$-Rigid Groups} II, Invent. Math.,
{\bf 165} (2006), 409-453. (math.OA/0407137).

\item{[S]} J. Schwartz: {\it Two finite, non-hyperfinite,
non-isomorphic factors}, Comm. Pure Appl. Math. {\bf 16} (1963),
19-26.

\item{[V1]} D. Voiculescu:  {\it Symmetries of some reduced free product 
$C^*$-algebras},  In: ``Operator algebras and their connections with topology
and ergodic theory'', Lect. Notes in Math. Vol. {\bf 1132}, 566-588 (1985).

\item{[V2]} D. Voiculescu: {\it Multiplication of certain noncommuting random variables.}, J. Operator Theory {\bf 18} (1987), 223-235. 

\item{[Wa]} S. Wasserman: {\it On tensor product of certain group C$^*$-algebras}, J. Funct. Anal. {\bf 23} (1976), 239-254. 

\item{[We]} N. Weaver: {\it The Kadison-Singer problem in discrepancy theory}, Discrete  Math. {\bf 278} (2004), 227-239. 

\item{[W]} F. B. Wright:  {\it A reduction for algebras of finite type},  Ann. of Math. {\bf 60} (1954), 560 - 570.

\enddocument